\documentclass[smallcondensed]{svjour3}     
\smartqed  
\usepackage{graphicx}
\usepackage{commath}
\usepackage{ wrapfig, subcaption, setspace, booktabs}
\graphicspath{ {./figures/} }
\usepackage[T1]{fontenc}
\usepackage[font=small, labelfont=bf]{caption}
\usepackage[english]{babel}
\usepackage{url, lipsum}
\usepackage{float}
\usepackage{hyperref}
 \usepackage{booktabs}
 \usepackage{multirow}
 \usepackage[table,xcdraw]{xcolor}
\usepackage[square,sort,comma,compress,numbers]{natbib}
\makeatletter
\DeclareRobustCommand\bigop[2][1]{%
  \mathop{\vphantom{\sum}\mathpalette\bigop@{{#1}{#2}}}\slimits@
}
\newcommand{\bigop@}[2]{\bigop@@#1#2}
\newcommand{\bigop@@}[3]{%
  \vcenter{%
    \sbox\z@{$#1\sum$}%
    \hbox{\resizebox{\ifx#1\displaystyle#2\fi\dimexpr\ht\z@+\dp\z@}{!}{$\m@th#3$}}%
  }%
}
\makeatother

\newcommand{\bigA}{\DOTSB\bigop[0.92]{\mathrm{A}}}

\newcommand{\ary}[1]{\boldsymbol{\mathsf{#1}}}

\hypersetup{
colorlinks=true,
linkcolor=black
}

\usepackage[utf8]{inputenc}
\usepackage[english]{babel}
\usepackage{mathrsfs}
 \usepackage{amsfonts}

\begin{document}

\title{Assessment of an Isogeometric Approach with Catmull-Clark Subdivision Surfaces using the Laplace-Beltrami Problems
}



\titlerunning{Assessment of Catmull-Clark Subdivision IGA}        

\author{Zhaowei Liu$^{1*}$ \and Andrew McBride$^1$ \and Prashant Saxena$^1$ \and Paul Steinmann$^{1,2}$}

\institute{* \email{zhaowei.liu@glasgow.ac.uk} \\
		  $1$. Glasgow Computational Engineering Centre, University of Glasgow, Glasgow, G12 8LT, United Kingdom,\\
           $2$. Chair of Applied Mechanics, Friedrich-Alexander Universit\"at Erlangen-N\"urnberg, Paul-Gordan-Str. 3, D-91052, Erlangen, Germany
}

\date{Received: date / Accepted: date}

\maketitle

\begin{abstract}
An isogeometric approach for solving the Laplace-Beltrami equation on a two-dimensional manifold embedded in three-dimensional space using a Galerkin method based on Catmull-Clark subdivision surfaces is presented and assessed. {The scalar-valued Laplace-Beltrami equation requires only $C^0$ continuity and 
is adopted to elucidate key features and properties of the isogeometric method using Catmull-Clark subdivision surfaces.} Catmull-Clark subdivision bases are used to discretise both the geometry and the physical field. A fitting method generates control meshes to approximate any given geometry with Catmull-Clark subdivision surfaces. The performance of the Catmull-Clark subdivision method is compared to the conventional finite element method. Subdivision surfaces without extraordinary vertices show the optimal convergence rate. However, extraordinary vertices introduce error,  which decreases the convergence rate. A comparative study shows the effect of the number and valences of the extraordinary vertices on accuracy and convergence. An adaptive quadrature scheme is shown to reduce the error.
\end{abstract}

\section{Introduction}
\label{intro}
\citet{hughes2005isogeometric} proposed the concept of isogeometric analysis (IGA) in 2005. The early works on IGA~\cite{cottrell2009isogeometric,reali2006isogeometric,bazilevs2008isogeometric} focussed on geometries modelled using Non-Uniform Rational B-Splines (NURBS) as these are widely used in computer aided design (CAD). NURBS can be used to model free-form, two-dimensional curves. However, a NURBS surface is a tensor product surface generated by two NURBS curves, thereby imposing limitations for modelling complex geometries with arbitrary topologies. Complex CAD models are always composed of a number of NURBS patches. These patches are often poorly connected in the design stage. When such models are used for analysis, the unmatched patches must be treated carefully to ensure the geometries are watertight. Furthermore, because NURBS can not be locally refined, adaptive mesh refinement method cannot be employed. A number of alternative CAD techniques were developed and adopted in IGA to overcome these limitations, including Hierarchical B-splines~\cite{forsey1988hierarchical,vuong2011hierarchical}, T-splines~\cite{sederberg2003t,bazilevscalo2010}, PHT-splines~\cite{deng2008polynomial,nguyen2011rotation}, { THB-splines~\cite{giannelli2012thb,buchegger2016adaptively} }and LR B-splines~\cite{dokken2013polynomial,johannessen2014isogeometric}. Some of these recent techniques are being adopted by the engineering design market. However, the majority are the subject of academic research and not widely used in the CAD community. Moreover, computing the basis functions for analysis using these alternative approaches can be expensive.  Catmull and Clark~\cite{catmull1978recursively} developed a bicubic B-spline patch subdivision algorithm for describing smooth three dimensional objects. The use of Catmull-Clark subdivision surfaces to model complex geometries in the animation and gaming industries dates back to 1978. Catmull-Clark subdivision surfaces can be considered as uniform bi-cubic splines which can be efficiently evaluated using polynomials.

{In CAD, distortion of regular parametrizations are inevitable and indeed vital when modelling complex geometries. Allowing `extraordinary vertices' ensures that Catmull-Clark subdivision surfaces can be used for modelling complex geometries with arbitrary topology. \citet{cirakortiz2000} implemented Loop subdivision surfaces for solving the Kirchhoff-Love shell formulation. This was the first application of subdivision surfaces to engineering problems. Subdivision surfaces have subsequently been used in electromagnetics~\cite{Dault2015}, shape optimisation~\cite{bandara2016shape, BANDARA201862} , acoustics~\cite{liu2018isogeometric, chen2020acoustic} and lattice-skin structures~\cite{xiao2019interrogation}.

{
Catmull-Clark subdivision surfaces face a number of challenges when used for analysis. Many of these have been discussed in the literature, however a unified assessment is lacking. This manuscript provides a clear and concise discussion of the challenges and limitations of Catmull-Clark subdivision surfaces.
 
Engineering designs often require exact geometries including circles, spheres, tori and cones. However, subdivision surfaces can not capture these geometries exactly. Moreover, there are always offsets between the control meshes and the surfaces. Fitting subdivision surfaces~\cite{litke2001fitting} aim to overcome this limitation. Although the fitting subdivision surfaces still can not model arbitrary geometries exactly as they are interpolated using cubic splines, they can approximate the given geometries closely through least-square fitting.  Another challenge of subdivision surfaces is that they can model smooth closed manifolds easily but require special treatment to model manifolds with boundaries. A common solution is to introduce `ghost' control vertices to provide bases for interpolating. From the perspective of analysis, the shape functions will span into `ghost' elements~\cite{cirakortiz2000}. In addition, the spline basis functions do not possess an interpolating property. Thus it is difficult to directly impose Dirichlet boundary conditions. Meshless methods and extended finite element methods have developed strategies to overcome this problem~\cite{fernandez2004imposing,moes2006imposing}. A common strategy is to modify the weak form of the governing equation. {Methods} include the Lagrangian multiplier method~\cite{babuvska1973finite}, the penalty method~\cite{arnold1982interior} and Nitsche's method~\cite{nitsche1971variationsprinzip,hansbo2002unfitted}. 

Conventional Catmull-Clark subdivision surfaces can not be locally refined. Truncated hierarchical Catmull-Clark subdivision surfaces (THCCS), developed by~\citet{wei2015truncated}, overcome this limitation. They generalise truncated hierarchical B-splines (THB-splines) to meshes with arbitrary topology. \citet{wei2016extended} subsequently improved their method using a new basis function insertion scheme and thereby enhanced the efficiency of local refinement. The extraordinary vertices introduce singularities in the parametrisation~\cite{takacs2012h2, nguyen2014comparative}. Catmull-Clark subdivision surfaces have $C^2$ continuity everywhere except at the surface points related to extraordinary vertices where, as demonstrated by~\citet{peters1998analysis}, they possess $C^1$ continuity. \citet{stam1998exact} developed a method to evaluate Catmull-Clark subdivision surfaces directly without explicitly subdividing, thus allowing one to evaluate elements containing extraordinary vertices. Although the surface gradients can not be evaluated at the extraordinary vertices, they can be evaluated at nearby quadrature points. Thus, subdivision surfaces can be used as $C^1$ elements as required, for example, in thin shell theory~\cite{cirakortiz2000}. Nevertheless, the evaluation of points around extraordinary vertices of Catmull-Clark surfaces introduces error. The conventional evaluation method repeatedly subdivides the element patch until the target point fall into a regular patch allowing a uniform bi-cubic B-spline patch to be mapped the subdivided element patch. The extraordinary vertex also introduces approximation errors because of the singular parametrisations at extraordinary vertices~\cite{peters1991parametrizing,neamtu1994degenerate}. Stam's natural parametrisation only can achieve $C^0$ continuity at extraordinary vertices. Recently \citet{wawrzinek2016integration} introduced a characteristic subdivision finite element scheme that adopted a characteristic reparameterisation for elements with extraordinary vertices. The evaluated limiting surface is at least $C^1$ everywhere and the numerical accuracy is improved. \citet{zhang2018subdivision} optimised the subdivision scheme to improve its approximation properties when used for thin-shell theory. }
}

{Using the finite element method to solve the partial differential equations (PDEs) on surfaces dates back to the seminal work by~\citet{dziuk1988finite}, which developed a variational formulation to approximate the solution of the Laplace-Beltrami problems on two dimensional surfaces. This method was extended to solve nonlinear and higher-order equations on surfaces by~\citet{dziuk2007surface}. \citet{dziuk2013finite} also provided a thorough review on finite element methods for approximating the solution of PDEs on surfaces. \citet{dedner2013analysis} proposed a discontinuous Galerkin (DG) method for solving a elliptic problem with the Laplace-Beltrami operator on surfaces. Adaptive DG~\cite{dedner2016adaptive} and high-order DG~\cite{antonietti2014high} methods were also developed for solving PDEs on surfaces. However, the accuracy of these methods depends on the approximation of the mean curvatures of the surfaces. The geometrical error is dominant when conventional Lagrangian discretisation is used to approximate solutions on complex surfaces. Isogeometric discretisation maintains the exact geometry and overcomes this limitation. \citet{dede2015isogeometric} proposed an isogeometric approach for approximating several surface PDEs involving the Laplace-Beltrami operator on NURBS surfaces. \citet{bartezzaghi2015isogeometric} solved PDEs with high order Laplace-Beltrami operators on surfaces using NURBS based isogeometric Galerkin method. More accurate results are obtained using an IGA approach over the conventional finite element method. \citet{langer2015multipatch} present an isogeometric DG method with non-matching NURBS patches allowing the approximation of PDEs on more complex surfaces. }

{This work presents a thorough and unified discussion of several major issues related to isogeometric Galerkin formulation based on Catmull-Clark subdivision surfaces. The difficulties associated with imposing Dirichlet boundary conditions, the reduction of the approximation power around extraordinary vertices, and
the problem of sufficient numerical integration in the element with extraordinary vertices will be examined and discussed.}
Previous studies~\cite{cirakortiz2000,Cirak:2001aa} on Catmull-Clark subdivision surfaces for analysis introduce ghost degrees of freedoms for constructing basis functions in elements at boundaries. We propose a method which {modifies} the basis functions at boundaries to ensure they are only associated with given control vertices. No additional ghost degrees of freedom are involved.  A penalty method is employed to impose Dirichlet boundary conditions. This does not change the size or symmetry of the system matrix and is straightforward to implement. An adaptive quadrature scheme inspired by~\cite{juttler2016numerical} is presented to increase the integration accuracy for elements with extraordinary vertices. The proposed method can perform isogeometric analysis on complex geometries using Catmull-Clark subdivision discretisations. A test for approximating Poisson's problem on a square plate is conducted to demonstrate the properties of the method  in a simplified setting so as to distill the key features. The approach is also used for solving the Laplace-Beltrami equation which is a benchmark problem for curved manifolds~\cite{juttler2016numerical,nguyen2014comparative}. A comparative convergence study is conducted between the Catmull-Clark subdivision method and the conventional finite element method. The effects of the extraordinary vertices and {modified bases at boundaries} on convergence are examined.
Catmull-Clark subdivision surfaces are limiting surfaces generated by successively subdividing given control meshes. {They are identical to uniform bi-cubic B-splines. Thus,  they have difficulty to represent desired geometries exactly. Here, a least-squares fitting method is used to fit any given geometry using Catmull-Clark subdivision surfaces.}

{This manuscript first summarises the subdivision algorithm and the evaluation method for Catmull-Clark subdivision surfaces.} Then, techniques for using Catmull-Clark for numerical analysis and improving accuracy are presented in Section~\ref{sec:techniques}. Section~\ref{sec:Laplace-Beltrami} presents the Laplace-Beltrami problem and Section~\ref{sec:FE} shows a Galerkin method with Catmull-Clark subdivision surface bases.  Section~\ref{sec:numerical} showcases the numerical results. 
\section{Catmull-Clark subdivision surfaces}
\label{sec:CC_curves_surfaces}
\begin{figure}
\centering
	\includegraphics[width=0.7\linewidth]{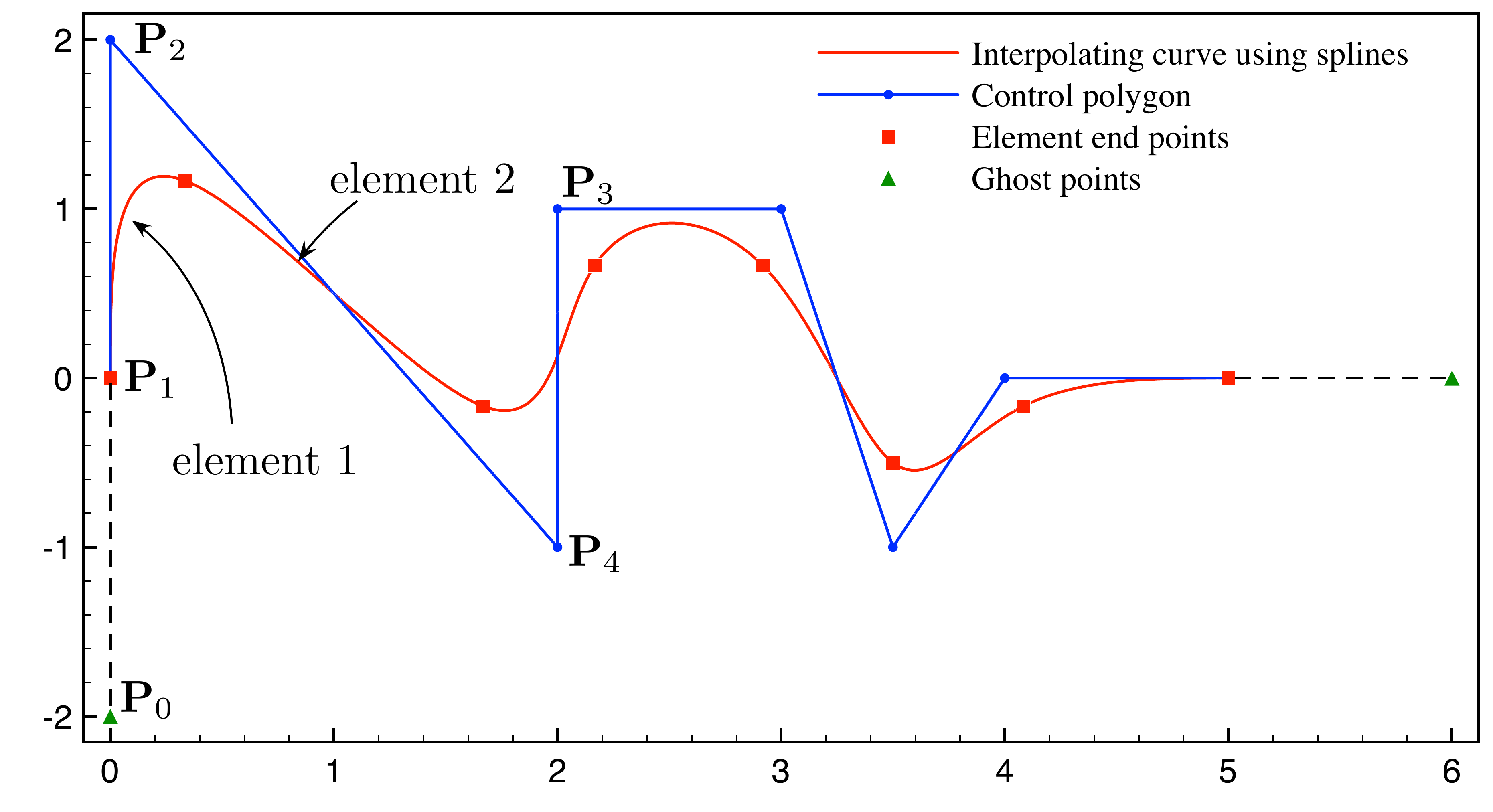}
	\caption{A subdivision curve is interpolated using basis splines and its control polygon.}
	\label{fig:interpolating_curve}
\end{figure}

\begin{figure}
\centering
\begin{subfigure}[b]{0.45\linewidth}
\centering
	\includegraphics[width=\linewidth]{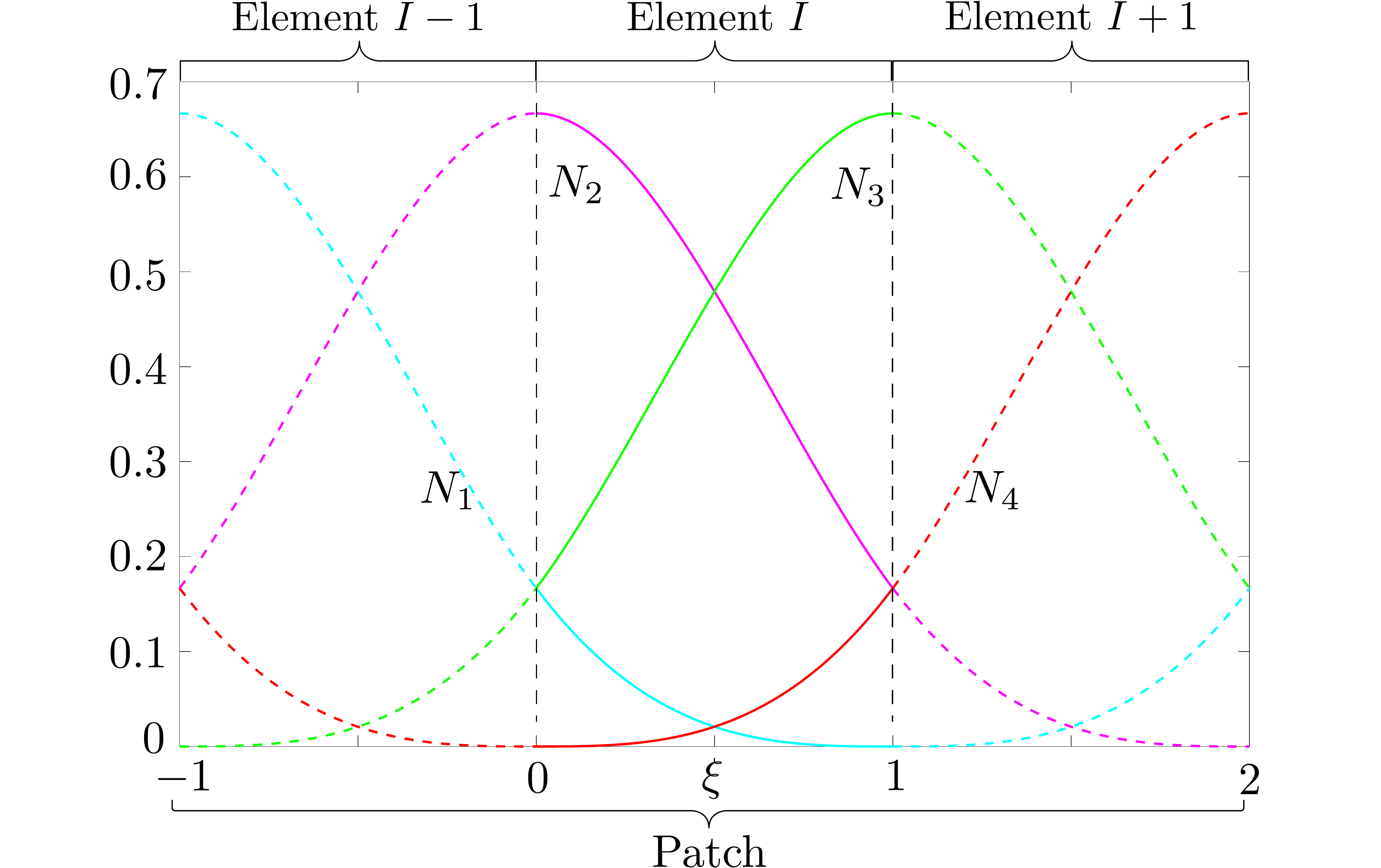}
	\caption{}
	\label{fig:basis}
\end{subfigure}
\begin{subfigure}[b]{0.54\linewidth}
\centering
	\includegraphics[width=\linewidth]{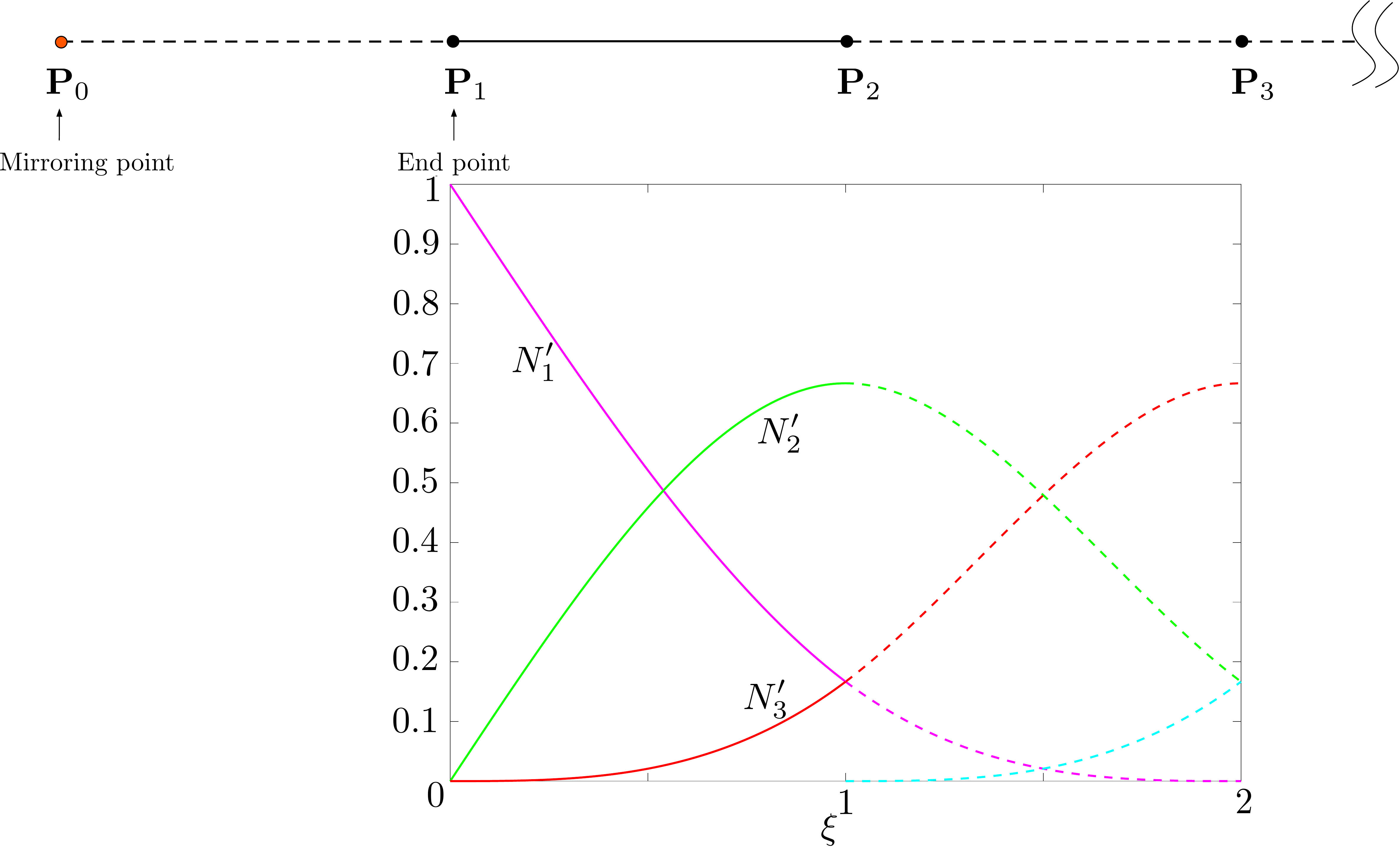}
	\caption{}
	\label{fig:basis_end}
\end{subfigure}
\begin{subfigure}[b]{0.6\linewidth}
\centering
	\includegraphics[width=\linewidth]{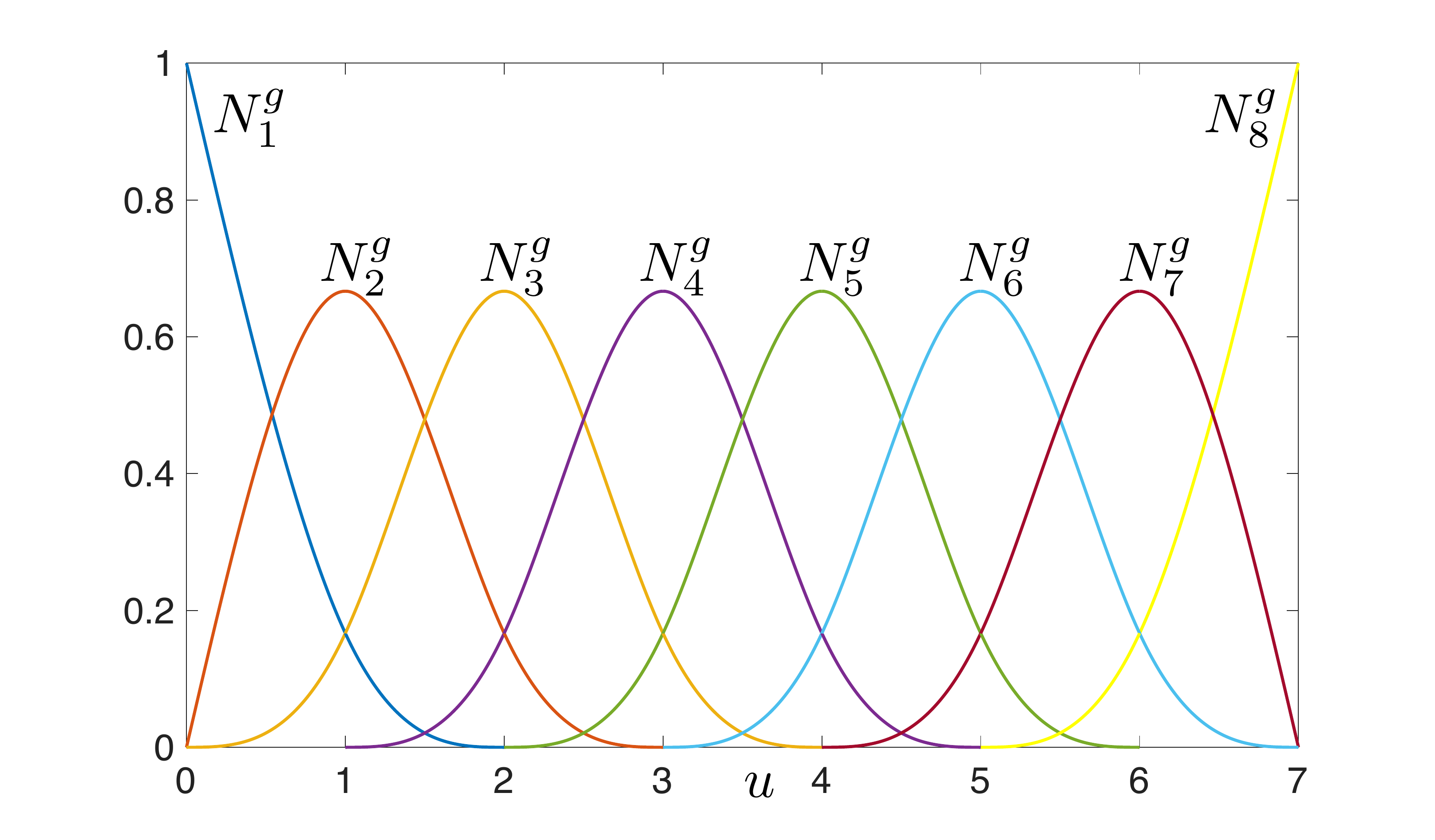}
	\caption{}
	\label{fig:bases}
\end{subfigure}
	\caption{(a) Basis splines for interpolating element $I$ as a Catmull-Clark curve. (b) The construction of mirroring ghost point to maintain the location of the end point. Basis splines are reconstructed for interpolating an end-element in a Catmull-Clark curve. (c) Global basis functions for interpolating a curve.}
\end{figure}

There exist a variety of subdivision schemes, but the basic idea is to use a subdivision scheme to generate a smooth surface through a limiting procedure of repeated refinement steps starting from an initial polygonal grid.  The Catmull-Clark algorithm can generate curves and surfaces which are identical to cubic B-splines. The algorithms for curves and surfaces are shown in Appendix~\ref{ap:A.1} and~\ref{ap:A.2}, respectively. {This section will briefly introduce the methods for interpolating and evaluating curves and surfaces using the Catmull-Clark subdivision algorithm.\\

\subsection{Curve interpolation and evaluation based on  the subdivision algorithm}}
Figure~\ref{fig:interpolating_curve} shows a curve generated using a subdivision algorithm. {The interpolated curve is identical to a cubic B-spline curve.} The limiting curve can be interpolated using cubic basis splines and associated control points. With a control polygon containing $n$ control points, the curve is naturally divided into $n-1$ elements. 
{Each element in the curve} is associated with one segment of the control polygon. To interpolate on the target element, four control points including the neighbouring control points are required. For example, if one aims to evaluate the geometry of element $2$ in Figure~\ref{fig:interpolating_curve}, the four control points $\mathbf{P}_1$,$\mathbf{P}_2$,$\mathbf{P}_3$ and $\mathbf{P}_4$ are required and the curve point is evaluated as
\begin{equation}
\mathbf{x}(\xi) = \sum_{A = 1}^4 N_A(\xi) \mathbf{P}_A,
\end{equation}
where $\xi \in [0,1]$ is the parametric coordinate within an element. The basis functions for element $2$ are defined by
\begin{align}
N_1(\xi) &= \frac{1}{6}[1-3\xi+3\xi^2-\xi^3], & N_2(\xi) &= \frac{1}{6}[4-6\xi^2+3\xi^3], \nonumber \\
N_3(\xi) &= \frac{1}{6}[1+3\xi+3\xi^2-3\xi^3], & N_4(\xi) &= \frac{1}{6}\xi^3.
\label{eq:basis}
\end{align}
The bases are visualised in Figure~\ref{fig:basis}. They are $C^2$ continuous across element boundaries.
Element $1$ in Figure~\ref{fig:interpolating_curve} contains the end of the curve, which has an end curve point that coincides with the control point. In order to evaluate this curve, one needs to mirror the point $\mathbf{P}_2$ to $\mathbf{P}_0$ as
\begin{equation}
\mathbf{P}_0 = 2\mathbf{P}_1 - \mathbf{P}_2.
\label{eq:mirroring}
\end{equation}
The curve point can now be evaluated using basis splines with a set of control points shown in Figure~\ref{fig:basis_end}. However, if one adopts {a spline discretisation} for analysis, this strategy of end element treatment will introduce additional `ghost-like' degrees of freedom. To avoid this problem, the expression for $\mathbf{P}_0$~\eqref{eq:mirroring} is substituted into the interpolating equation yielding
\begin{equation}
\mathbf{x}(\xi) = \sum_{A=0}^3 N_{A+1}(\xi) \mathbf{P}_A = \sum_{B = 1 }^3 N'_B(\xi)\mathbf{P}_B.
\label{eq:boundary_basis}
\end{equation}
Hence only three control points are required to evaluate a curve point and the modified basis functions for interpolating end elements are defined by
\begin{align}
N'_1(\xi) &= \frac{1}{6}[6-6\xi+\xi^3], & N'_2(\xi) &= \frac{1}{6}[6\xi - 2\xi^3], & N'_3(\xi) &= \frac{1}{6}\xi^3.
\label{eq:basis_boundary}
\end{align}
Figure~\ref{fig:basis_end} illustrates the modified basis functions. {It achieves the same basis functions as the cubic B-Spline with $p+1$ multiple knots at the two end points.} The new basis functions do not possess the Kronecker delta property but do have the interpolating property at the boundary. {The performance of modified bases in analysis will be discussed in Section~\ref{sec:patch}}.

The global basis functions for interpolating the curve in Figure~\ref{fig:interpolating_curve} are shown in Figure~\ref{fig:bases}. {It is worth noting that this subdivision curve is a cubic B-spline curve and represents a special case of Lane-Riesenfeld subdivision it can not model {conical shapes} exactly.} This property is significantly different to NURBS and motivates Section~\ref{sec:geometry_fitting} on geometry fitting.\\

\subsection{Interpolating and evaluating Catmull-Clark subdivision surfaces}
One defines the number of elements connected with the vertex as the valence. A regular vertex in a Catmull-Clark surface mesh has a valence of 4.  A vertex with a valence not equal to 4 is called an extraordinary vertex. This allows subdivision surfaces to handle arbitrary topologies. In their seminal paper~\cite{catmull1978recursively}, Catmull and Clark proposed a way to modify the weight distributions for extraordinary vertices in order to describe complex geometries. With this simple solution, Catmull-Clark surfaces can use a single mesh to present surfaces of arbitrary geometries while other spline-based CAD tools, such as NURBS surfaces, need to link multiple patches. {The limiting surface of the Catmull-Clark subdivision algorithm has $C^2$ continuity over the surface except at the extraordinary vertices where they have $C^1$ continuity as proven by~\citet{peters1998analysis}. } This section will illustrate the methods of interpolating and evaluating both regular element and element with an extraordinary vertex in Catmull-Clark subdivision surfaces.
\subsubsection*{Element in a regular patch}
\begin{figure}
\centering
\begin{subfigure}[b]{0.51\linewidth}
\centering
	\includegraphics[width=0.9\linewidth]{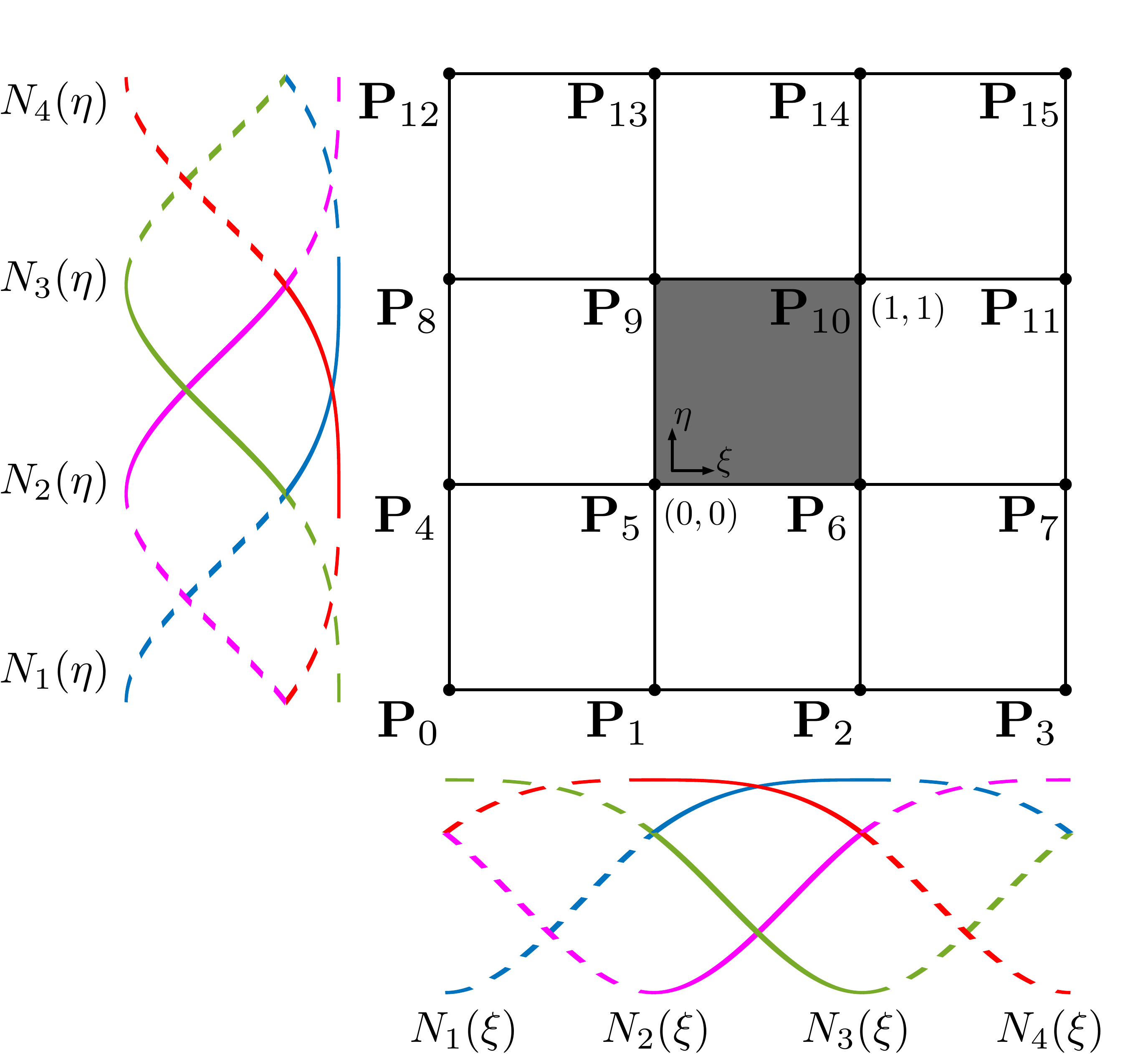}
	\caption{}
	\label{fig:CC_element_patch}
\end{subfigure}
\begin{subfigure}[b]{0.49\linewidth}
\centering
	\includegraphics[width=\linewidth]{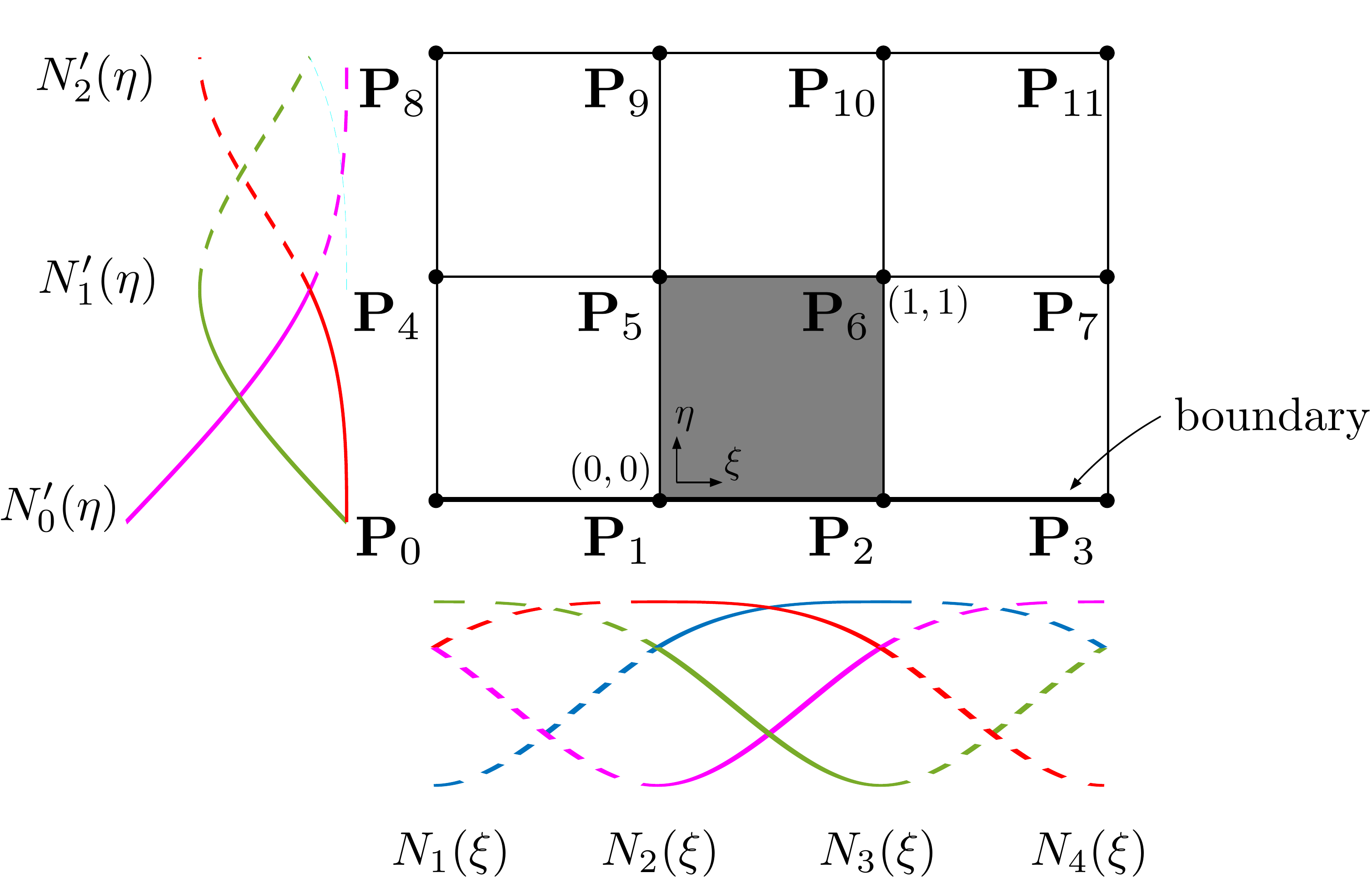}
	\caption{}
	\label{fig:CC_element_patch_boundary}
\end{subfigure}
\begin{subfigure}[b]{0.49\linewidth}
\centering
	\includegraphics[width=0.7\linewidth]{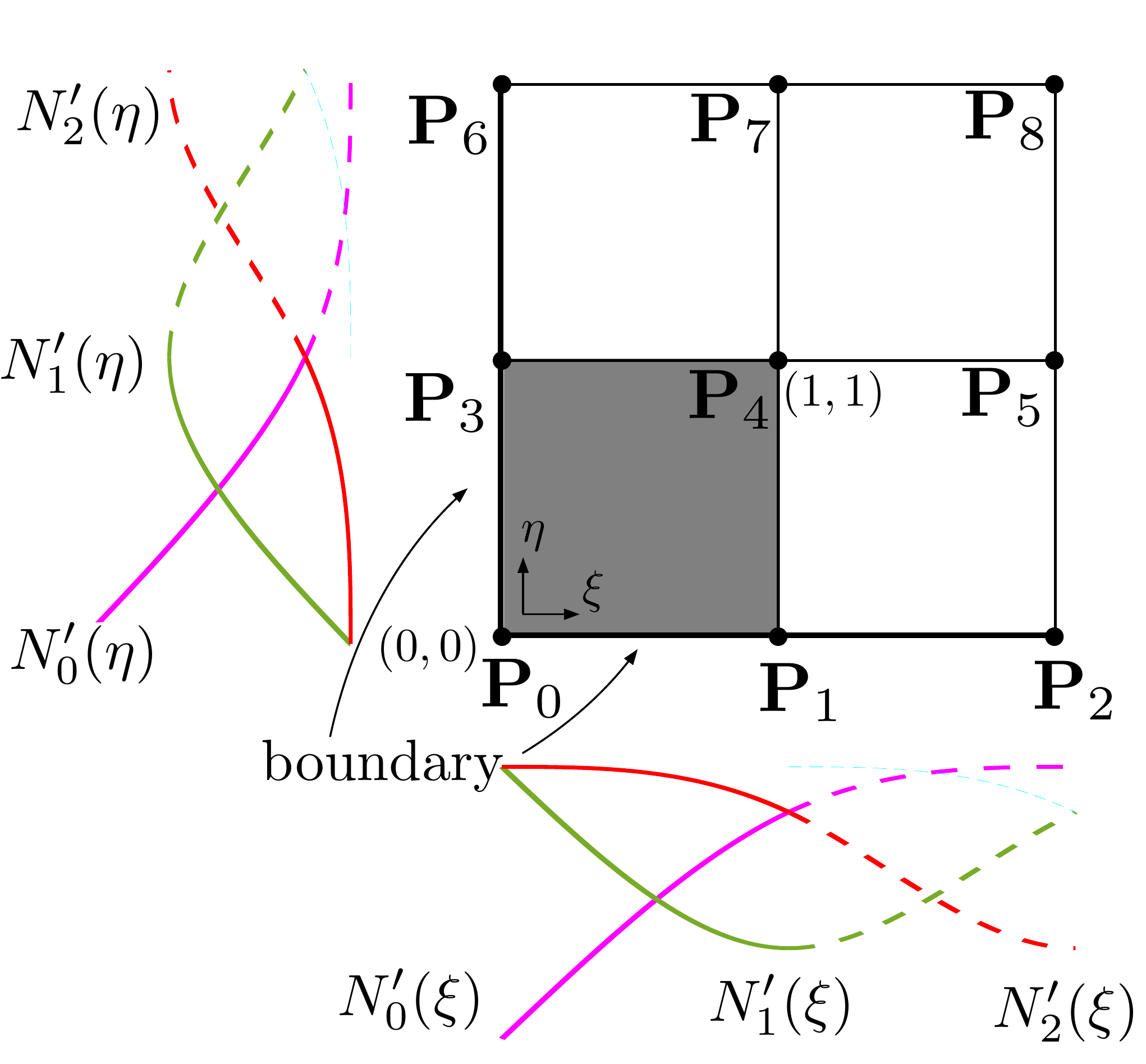}
	\caption{}
	\label{fig:CC_element_patch_corner}
\end{subfigure}
\caption{Element patches for evaluating a Catmull-Clark subdivision element. (a) A regular element. (b) An element with a face on the boundary. (c) An element with two faces on the boundary.}
	\label{fig:CC_element_patches}
\end{figure}

Figure~\ref{fig:CC_element_patch} shows a subdivision surface element~(dashed) which does not contain an extraordinary vertex. In order to evaluate a point in this Catmull-Clark element, an element patch must be formed. The patch consists of the element itself and the elements which share vertices with it. A regular element patch has $9$ elements with $16$ control vertices. The surface point can be evaluated using the $16$ basis functions associated with these control points as
\begin{equation}
\mathbf{x}(\boldsymbol{\xi}) = \sum_{A = 0}^{15} N_A(\boldsymbol{\xi}) \mathbf{P}_A,
\end{equation}
where $\boldsymbol{\xi}:=(\xi,\eta)$ is the parametric coordinate of a Catmull-Clark subdivision surface element. A Catmull-Clark surface is obtained as the tensor product of two Catmull-Clark curves. The basis functions are defined by
\begin{equation}
N_i(\boldsymbol{\xi}) = N_{i\%4}(\xi)N_{\lfloor i/4 \rfloor}(\eta), \,  i =0,1,\dots,15,
\label{eq:cc_bases_tp}
\end{equation}
where $N(\xi)$ or $N(\eta)$ are the basis functions defined in Equation~\eqref{eq:basis} and presented in Figure~\ref{fig:CC_element_patch}. $\lfloor \bullet \rfloor$ is the modulus operator and $\%$ denotes the remainder operator which gives the remainder of the integer division.

Figure~\ref{fig:CC_element_patch_boundary} shows the element patch of a subdivision surface element~(shaded) that has an edge on the physical boundary. This type of element has only $5$ neighbour elements so that it belongs to an element patch which has $12$ control vertices. To evaluate this element, a common solution is to generate a set of `ghost' vertices outside the domain to form a full element patch~\cite{cirakortiz2000}. However, this method involves additional degrees of freedom in numerical analysis. Instead, the curve basis functions in Equation~\eqref{eq:basis_boundary} are adapted to deal with the element on the boundary. The same strategy is used for elements which have two edges on the physical boundary as shown in Figure~\ref{fig:CC_element_patch_corner}.

\subsubsection*{Element in a patch with an extraordinary vertex}
\begin{figure}
\centering
\begin{subfigure}[b]{0.49\linewidth}
\centering
	\includegraphics[width=0.8\linewidth]{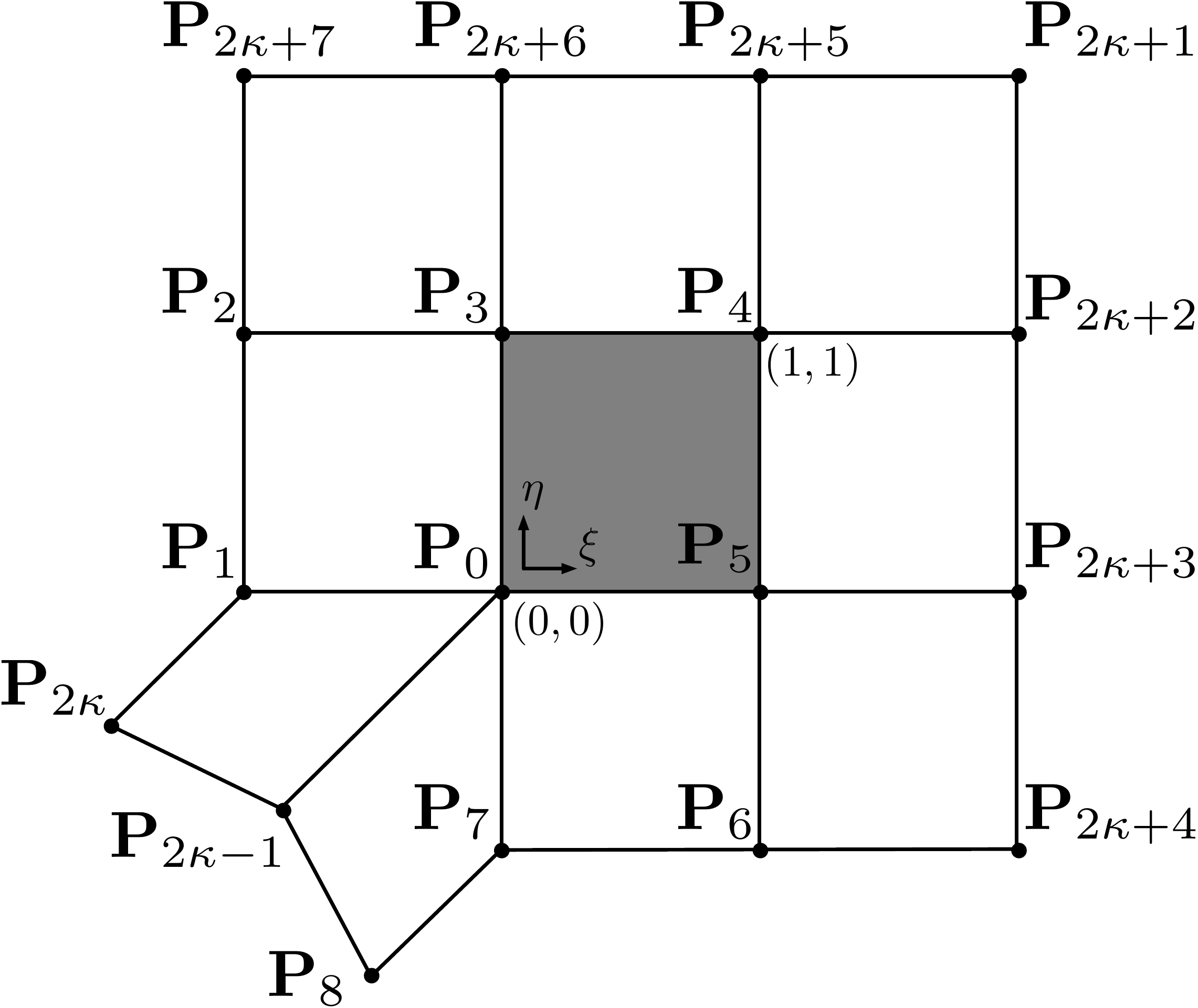}
	\caption{}
	\label{fig:Catmull_clark_ex_mesh}
\end{subfigure}
\begin{subfigure}[b]{0.49\linewidth}
\centering
	\includegraphics[width=0.65\linewidth]{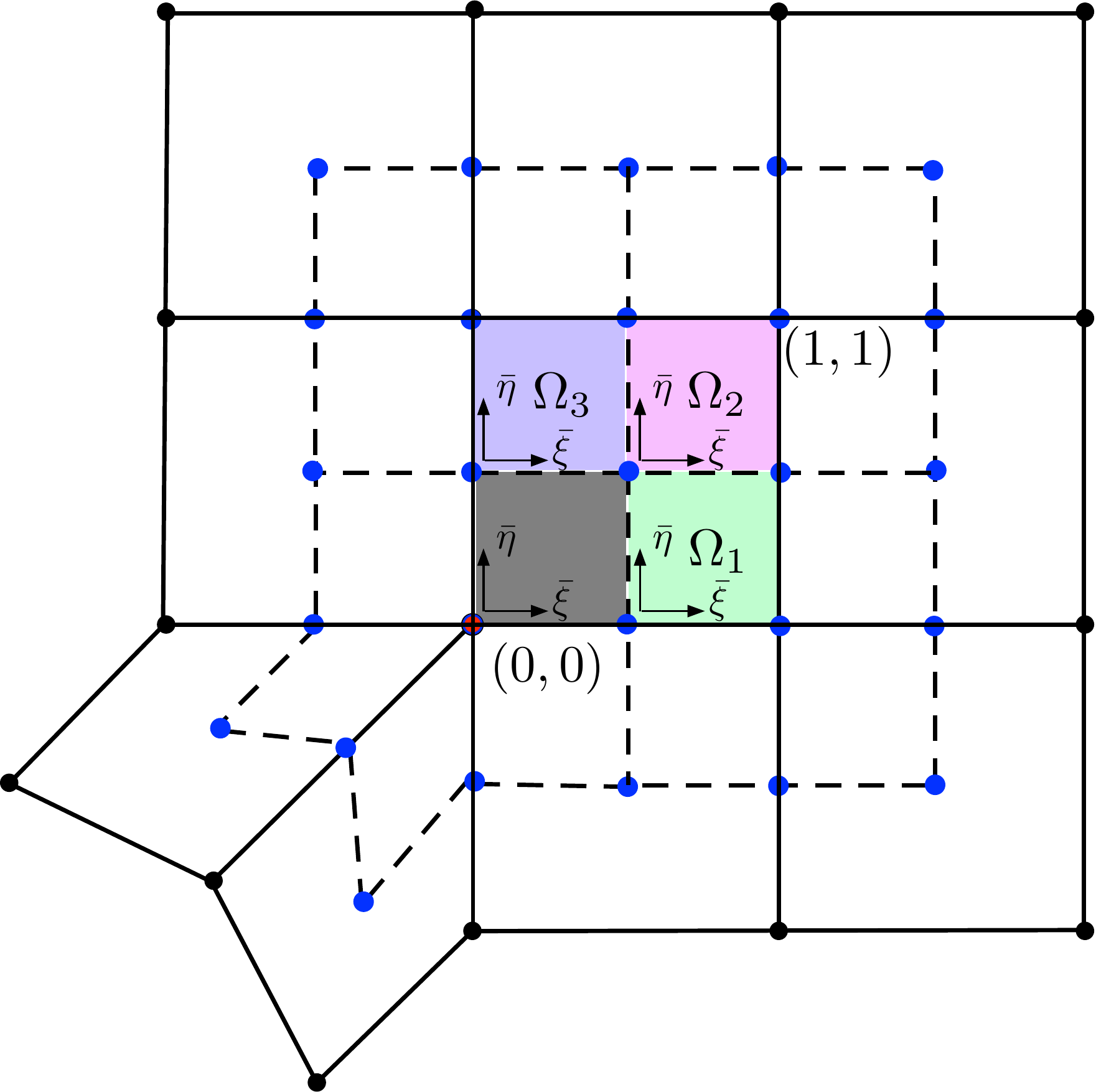}
	\caption{}
	\label{fig:catmull_clark_ex}
\end{subfigure}
\begin{subfigure}[b]{0.49\linewidth}
\centering
	\includegraphics[width=0.63\linewidth]{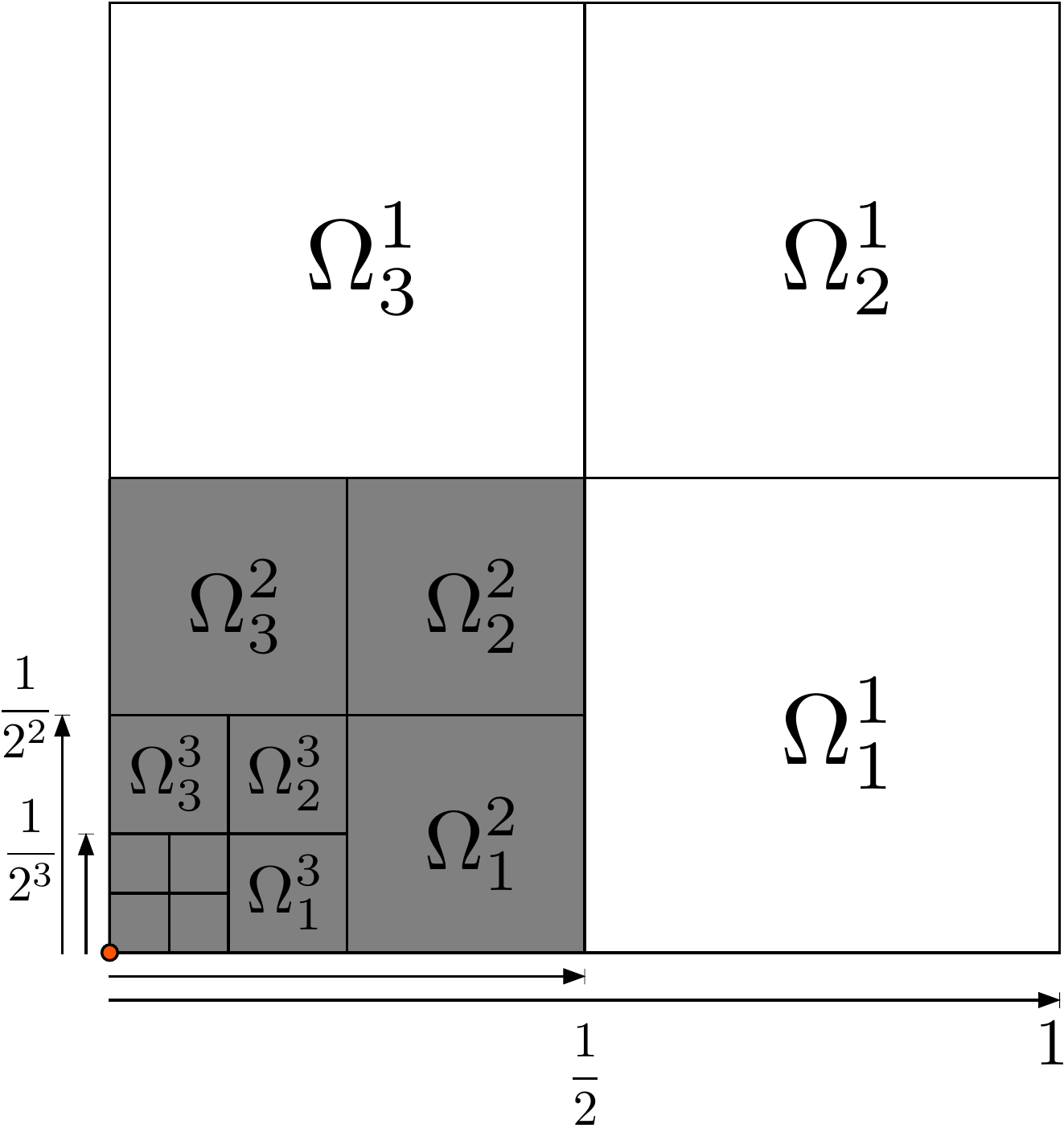}
	\caption{ }
	\label{fig:CC_domain_subdivision}
\end{subfigure}
\begin{subfigure}[b]{0.49\linewidth}
\centering
	\includegraphics[width=\linewidth]{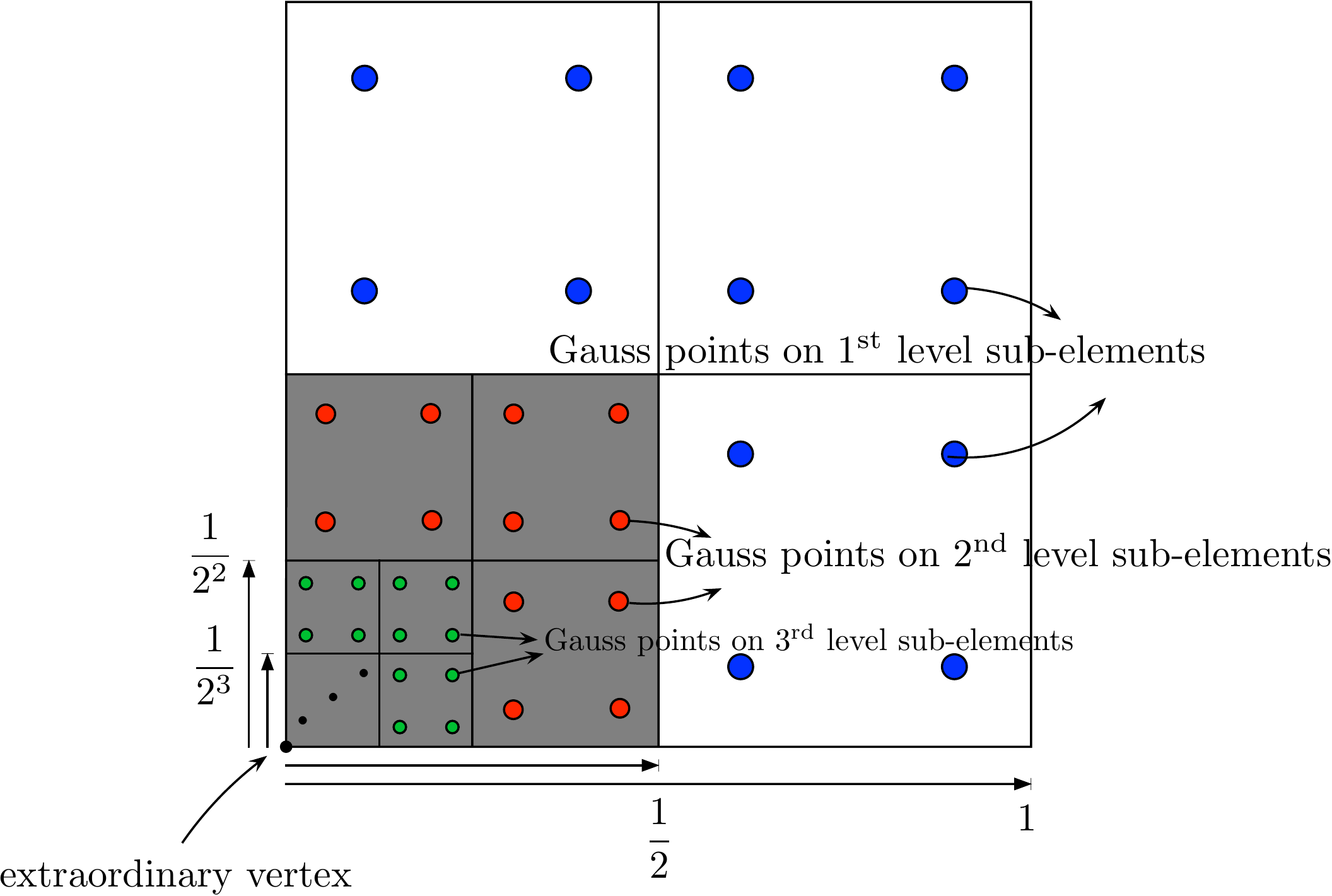}
	\caption{}
	\label{fig:CC_domain_adaptive_quadrature}
\end{subfigure}
\caption{(a) An irregular patch for a Catmull-Clark subdivision element with an extraordinary vertex.  (b) One level of subdivision of the element patch divides the element into four sub-elements; three of them are in regular patches. (c) Successive subdivisions of the element until the evaluated point falls into a sub-element with a regular patch. (d) Adaptive Gauss quadrature scheme for the element with an extraordinary vertex.}
\end{figure}

Extraordinary vertices are a key advantage of Catmull-Clark subdivision surfaces which allows them to model complex geometries with arbitrary topologies. However, it increases the difficulty of evaluating the surfaces. Figure~\ref{fig:Catmull_clark_ex_mesh} shows a Catmull-Clark subdivision element which contains one extraordinary vertex.

In order to evaluate this element, one needs to re-numerate the control points as shown in Figure~\ref{fig:Catmull_clark_ex_mesh}. After applying one level of subdivision, new control points are generated and this element is subdivided into four sub-elements, as shown in Figure~\ref{fig:catmull_clark_ex}. The sub-elements $\Omega_1$, $\Omega_2$ and $\Omega_3$ are now in a regular patch. However, the last sub-element (grey) still has an extraordinary vertex. If the target point to be evaluated is in this region, we must repeatedly subdivide the element until the point falls into a sub-element with a regular patch. Then, the point can be evaluated within the sub-element with the new set of control points $\mathscr P_{n,k}$, where $n$ is the number of subdivision required and $k = 1, 2, 3$ is the sub-element index shown in Figure~\ref{fig:catmull_clark_ex}. The new control point set is computed as
\begin{equation}
\mathscr P_{n,k} = {\ary D_k \ary A \bar{\ary A}^{n-1}} \mathscr P_0,
\end{equation}
where $\ary D_k$ is a selection operator to pick control points for the sub-elements. $\ary A$ and $\bar{\ary A}$ are two types of subdivision operators. $\mathscr  P_0$ is the initial set of control points.  The detailed approach is given in~\cite{stam1998exact} and also can be found in Appendix~\ref{ap:bases_ev}. $\mathscr  P_{n,k}$ contains 16 control points. Then, a surface point in the element with an extraordinary vertex can be computed as
\begin{equation}
\mathbf{x}(\boldsymbol{\xi}) =  \sum_{A=0}^{15} N_A(\bar {\boldsymbol{ \xi}})\, \mathbf{P}^{n,k}_{A},
\label{eq:subd_interpolating}
\end{equation}
where $\bar{ \boldsymbol{ \xi}}$ is the parametric coordinates of the evaluated point in the sub-element, which can be mapped from $\boldsymbol{ \xi}$ as
\begin{equation}
\bar{ \boldsymbol{ \xi}} = 
\left\{
\begin{split}
&(2^n\xi - 1, 2^n\eta)\quad &\text{if} \quad k =1 \\
&(2^n\xi - 1, 2^n\eta-1)\quad & \text{if} \quad k =2\\
&(2^n\xi , 2^n\eta-1)\quad  &\text{if} \quad k =3
\end{split}
\right.
.
\label{eq:mapping_2}
\end{equation}
Equation~\eqref{eq:subd_interpolating} can thus be rewritten as 
\begin{equation}
\mathbf{x}(\boldsymbol{\xi}) =  \sum_{A=0}^{2\kappa+7} \hat{N}_A( {\boldsymbol{ \xi}})\, \mathbf{P}^{0}_{A},
\end{equation}
where $\hat{N}$ is the Catmull-Clark subdivision surfaces basis function. Define $\hat{\mathbf{N}}$ as a set of $2\kappa+8$ basis functions in an element with an extraordinary vertex and $\mathbf{N}$ is a set of $16$ regular basis functions defined in Equation~\eqref{eq:cc_bases_tp}. $\hat{\mathbf{N}}$ can be calculated in a vector form as
\begin{equation}
 \hat{\mathbf{N}}(\boldsymbol{ \xi}) = {[\ary D_k \ary A \bar{\ary A}^{n-1}]^T \mathbf{N}}(\bar{ \boldsymbol{ \xi}}).
\label{eq:subd_basis}
\end{equation}
The derivatives of the Catmull-Clark subdivision surfaces basis functions for elements containing extraordinary vertices are expressed as
\begin{equation}
\frac{\partial \hat{\mathbf{N}}(\boldsymbol{ \xi}) }{\partial \mathbf{\boldsymbol{\xi}}}=\left[
\begin{array}{cc}
\frac{\partial \hat{N}_0}{\partial \xi} & \frac{\partial \hat{N}_0}{\partial \eta} \\
\frac{\partial \hat{N}_1}{\partial \xi} & \frac{\partial \hat{N}_1}{\partial \eta} \\
\vdots & \vdots \\
\frac{\partial \hat{N}_{2\kappa+7}}{\partial \xi} & \frac{\partial \hat{N}_{2\kappa+7}}{\partial \eta} \\
\end{array}
\right],
\end{equation}
 and can be computed by
\begin{equation}
\frac{\partial \hat{\mathbf{N}}(\boldsymbol{ \xi}) }{\partial \mathbf{\boldsymbol{\xi}}} =[\ary D_k \ary A \bar{\ary A}^{n-1}]^T \frac{\partial {{\mathbf{N}}}(\bar{\boldsymbol{ \xi}})}{\partial \bar{ \boldsymbol{ \xi}}}	\frac{\partial \bar{\boldsymbol{\xi}}}{\partial \boldsymbol{\xi}},
\end{equation}
where $\frac{\partial \bar{\boldsymbol{\xi}}}{\partial \boldsymbol{\xi}}$ can be considered as a mapping matrix defined by
\begin{equation}
\frac{\partial \bar{\boldsymbol{\xi}}}{\partial \boldsymbol{\xi}} = \left[
\begin{array}{cc}
2^n & 0\\
0 & 2^n
\end{array}
\right].
\label{eq:singular_problem}
\end{equation}

\begin{remark}
The calculation of the basis functions $\hat{\mathbf{N}}$ at a physical point $\mathbf{x}$ involves two mappings. The first is from the physical domain to the parametric domain of an element with an irregular patch, $\mathbf{x} \mapsto \boldsymbol{\xi}$. Because the irregular patch does not have the tensor-product nature, $n$ levels of subdivisions are required and the point is mapped to the parametric domain of a sub-element, $\boldsymbol{\xi} \mapsto \bar{\boldsymbol{\xi}}$. This second mapping is defined in Equation~\eqref{eq:mapping_2}. The value of $n$ approaches positive infinity when $\boldsymbol{\xi}$ approaches the extraordinary vertex which has the parametric coordinate $(0,0)$. Hence the diagonal terms in the mapping matrix~\eqref{eq:singular_problem} tend to positive infinity as $n \to \infty$. This results in the basis functions $\hat{\mathbf{N}}$ not being differentiable at $\boldsymbol{\xi} = \mathbf{0}$. This problem is termed singular configuration in~\cite{ juttler2016numerical}, and singular parameterisation in~\cite{takacs2012h2,nguyen2014comparative}.
 \label{rm:singular}
\end{remark} 
\section{Techniques for analysis and improving accuracy}
This section presents three techniques which are essential for using Catmull-Clark subdivision surfaces in numerical analysis.  A geometry fitting method using Catmull-Clark surfaces is introduced in Section~\ref{sec:geometry_fitting}. Section~\ref{sec:aq_ev} illustrates an adaptive quadrature scheme  for integrating element with an extraordinary vertex to improve accuracy. Section~\ref{sec:penalty} introduces the penalty method for applying essential boundary conditions. 
\label{sec:techniques}
\subsection{Geometry fitting}
\label{sec:geometry_fitting}
\begin{figure}
\centering
\begin{subfigure}{0.75\linewidth}
\includegraphics[width=\linewidth]{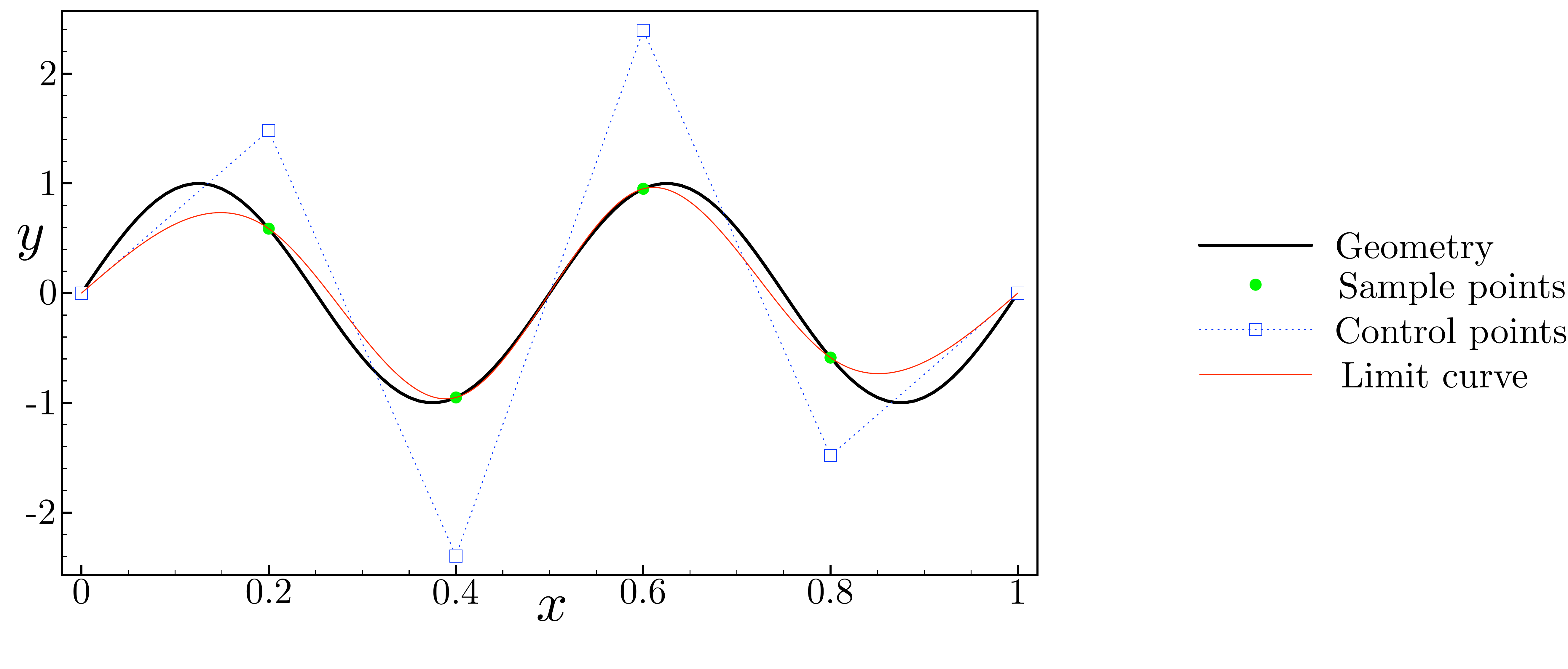}
\caption{Fitting curve with $n_s=6$ sampling points.}
\label{fig:fitting_1}
\end{subfigure}
\begin{subfigure}{0.49\linewidth}
\includegraphics[width=\linewidth]{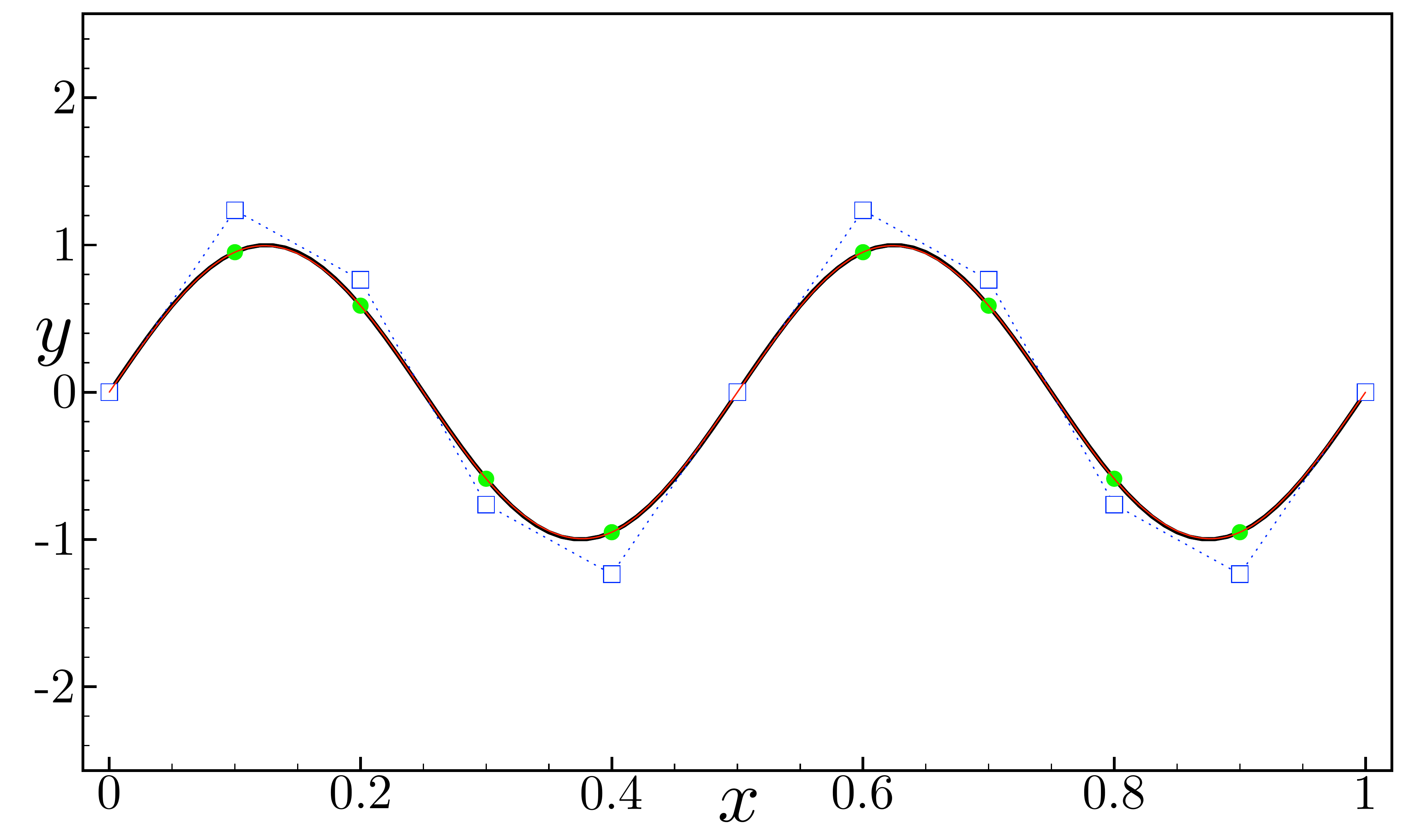}
\caption{Fitting curve with $n_s =11$ sampling points.}
\label{fig:fitting_2}
\end{subfigure}
\begin{subfigure}{0.49\linewidth}
\includegraphics[width=\linewidth]{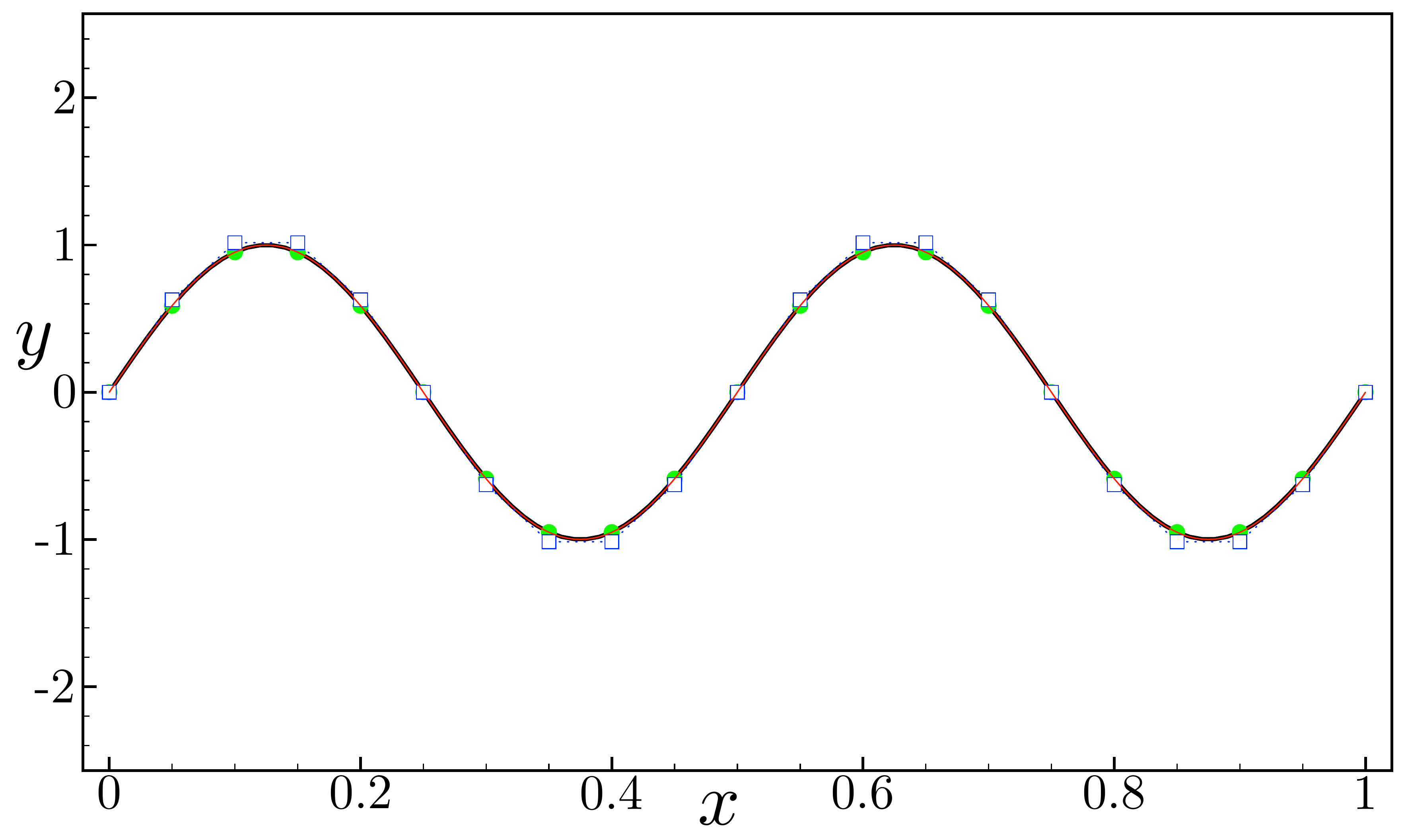}
\caption{Fitting curve with $n_s = 21$ sampling points.}
\label{fig:fitting_3}
\end{subfigure}
\caption{{The process of constructing cubic B-spline curves to approximate a given curve. }}
\label{fig:fitting}
\end{figure}

{Catmull-Clark subdivision surfaces are CAD tools which can construct limiting surfaces from control polygons and meshes. However, in a number of engineering problems, the geometry is given as an industry design and a limit surface that is a "best approximation" of this desired geometry required. \citet{litke2001fitting} introduced a method for fitting a Catmull-Clark subdivision surface to a given shape. They employed, both a least-squares fitting method and a quasi-interpolation method to determine a set of control points for a given surface.} The least-square fitting method is used here. One first chooses a set of sample points $\mathscr{S} = \{ \mathbf{s}_1, \mathbf{s}_2, \dots, \mathbf{s}_{n_s} \} \in \Gamma$, where $\Gamma$ is the geometry, $n_s$ is the number of sample points. Each sample point should be evaluated using Catmull-Clark subdivision bases with control points as
\begin{equation}
\mathbf{s}(\boldsymbol{ \xi}) =\sum_{A = 1}^{n_b} N_{A}(\boldsymbol{ \xi})\mathbf{P}_{A},
\end{equation}
where $n_b = 2\kappa+8$ is the number of local basis functions. Then the set of sample points can be evaluated as
\begin{equation}
\mathscr{S} = \ary{L} \mathscr{P},
\end{equation}
where $\mathscr{P} = \{ \mathbf{P}_1, \mathbf{P}_2, \cdots, \mathbf{P}_{n_c} \}$ is a set of control points with $n_c$ control points. $\ary L$ is an evaluation operator of Catmull-Clark curves or surfaces. Set $\boldsymbol{ \xi} = (0,0)$ to ensure the sample points correspond to the control vertices and $n_s \equiv n_c$, then $\ary L$ is a square matrix. The control points can be calculated as
\begin{equation}
\mathscr P = \ary L^{-1} \mathscr S.
\label{eq:fitting_inv}
\end{equation}
If more sampling points $n_s$ are chosen than the required number of control points $n_c$, then $\ary L$ is invertible, a least-squares method is used to obtain a set of control points $\hat{\mathscr P}$ that minimises $\norm{\mathscr S - \ary L \mathscr P}^2$ as
\begin{equation}
\hat{\mathscr P} = [\ary {L}^T \ary L]^{-1}\ary {L}^T\mathscr{S}.
\end{equation}
{Figure~\ref{fig:fitting} shows an example of fitting a geometry using cubic B-spline curves based on the Catmull-Clark subdivision algorithm.} The given curve is defined as $y = \sin(4\pi x)$. Figure~\ref{fig:fitting_1} shows that 6 sample points are chosen from the given curve and one assembles the evaluation operator for these sampling points. The control points can be obtained by solving~\eqref{eq:fitting_inv}. Using these control points, the limit curve can be interpolated. Since 6 sample points is not sufficient to capture the given curve, the limit curve is significantly different to the given curve. Figures~\ref{fig:fitting_2} and~\ref{fig:fitting_3} show the curve fitting with 11 and 21 sample points, respectively. Increasing the number of samples points, the limit curve converges to the given curve.

\subsection{Adaptive quadrature rule for element with an extraordinary vertex}
\label{sec:aq_ev}
In numerical analysis, a Gauss quadrature rule is applied to integrate over Catmull-Clark subdivision elements. A one dimensional quadrature rule with $n_q$ Gauss points can exactly evaluate the integrals for polynomials of degree up to $2n_q-1$. The polynomial degree of a cubic B-spline function is $3$. Because the basis functions of a Catmull-Clark subdivision element in regular element patch are generated as the tensor product of two cubic splines, $2 \times 2$ Gauss points can be used in this case. However, if a Catmull-Clark subdivision element has an extraordinary vertex, the basis functions are generated by Equation~\eqref{eq:subd_basis}. In this case, basis functions are not polynomials and the derivatives of the basis functions suffer from the singular parametrisation problem, see Remark~\ref{rm:singular}. Thus, the standard Gauss quadrature can not be used to evaluate the element integral. Inspired by~\cite{juttler2016numerical}, an adaptive quadrature rule, well suited to Catmull-Clark subdivision surfaces is adopted by integration at a number of levels of subdivided elements. With $n_d$ levels of subdivisions, the element is subdivided into $3n_d+1$ sub-elements as shown in Figure~\ref{fig:CC_domain_adaptive_quadrature}. The sub-elements can be evaluated using cubic B-splines with new control vertices except for the ones having an extraordinary vertex. Thus the Gauss quadrature rule can be used to evaluate the integrals in $3n_d$ sub-elements. With a number of subdivisions, the integration error can be reduced. In this work, $n_d = 7$ is chosen in order to obtain sufficiently accurate values of the integrals.

\subsection{Penalty method for applying boundary condition}
\label{sec:penalty}
 The basis functions do not have the Kronecker delta and interpolating properties, so boundary conditions can not be directly applied using conventional methods. The method used here is a penalty method which uses a penalty parameter and boundary mass matrix to apply the boundary conditions approximately. It preserves the symmetry of the system matrix and does not increase its size. However, the penalty parameter should be carefully selected. If fine meshes with more degrees of freedom are adopted, a larger penalty parameter must be chosen. The Dirichlet boundary condition is defined as
\begin{equation}
u = \bar u \quad \text{on} \quad \partial \Gamma.
\label{eq:Poisson_bc}
\end{equation}
{An $L_2$ projection is used for applying} the Dirichlet boundary condition, where for test function $v$, one obtains
\begin{equation}
\int _{\partial \Gamma} v u \,\mathrm{d}L = \int_{\partial \Gamma} v \bar u \,\mathrm{d}L, \quad \forall v \in H^1_0(\partial \Gamma).
\label{eq:boundary_weak_form}
\end{equation}
{Using the cubic B-spline functions} in Equation~\eqref{eq:basis} to discretise $u$ and $v$ and the same strategy for formulating the system matrix, one introduces a boundary mass matrix as
\begin{equation}
\ary M_b = \bigA_{e= 1}^{n^b_{e}}\sum_{i = 1}^{n_q} \mathbf{G}^e(\boldsymbol{\xi}^i)  |J^j(\mathbf{x}({\xi}^{i}))| w^{i},
\end{equation}
where $n^b_e$ is the number of boundary elements, and
\begin{equation}
\mathbf{G}^e({\xi}) =
\left[
\begin{array}{cccc}
 N_1 ({\xi})   N_1 ({\xi})  &  N_1 ({\xi})   N_2 ({\xi}) & \cdots &  N_1 ({\xi})   N_{n_b} ({\xi})\\
  N_2 ({\xi})   N_1 ({\xi})  &  N_2 ({\xi})   N_2 ({\xi}) & \cdots  &  N_2 ({\xi})   N_{n_b} ({\xi})\\
\vdots & \vdots & \ddots & \vdots\\
 N_{n_b}  ({\xi})   N_1 ({\xi})  &  N_{n_b}  ({\xi})   N_2 ({\xi}) & \cdots &  N_{n_b} ({\xi})   N_{n_b} ({\xi})\\
\end{array}
\right].
\end{equation} 
The right hand side vector for applying the boundary conditions is thus
\begin{equation}
\ary f_b  = \bigA_{j= 1}^{n^b_{e}}\sum_{i = 1}^{n_q} \bar{u}(\mathbf{x}(\boldsymbol{\xi}_{i})) \mathbf{N}^j(\boldsymbol{\xi}_{i})  |J^j(\mathbf{x}(\boldsymbol{\xi}_{i}))| w_{i}.
\end{equation}
Then the discrete system of equations arising from~\eqref{eq:boundary_weak_form} is
\begin{equation}
\ary M_b \ary u = \ary f_b.
\label{eq:system_equation_bc}
\end{equation}
We note that the elements for applying boundary conditions are the discretisation of the surface boundary which are {one dimensional cubic B-spline curves} and only one-dimensional Gauss quadrature rule is used for integration. However, one uses the global degrees of freedom indices to assemble $\ary M_b$ and $\ary f_b$, so that they have the same size as the system matrix and global right-hand side vector, respectively. Assume the system of equations is expressed as
$\ary K \ary u = \ary f$, where $\ary K$ is the system matrix, $\ary u$ is the global coefficients vector to be solved for, and $\ary f$ is global right-hand side vector. Then, we scale $\ary M_b$ and $\ary f_b$ using a penalty factor $\beta$ and combine them with the systems of equations as
\begin{equation}
[\ary K + \beta\ary M_b] \ary u = \ary f + \beta \ary f_b. 
\end{equation}
The Dirichlet boundary condition~\eqref{eq:Poisson_bc} is here weakly applied to the system of equations. {A relatively large penalty factor $\beta = 10^8$ is selected for all numerical examples. It is sufficiently large to ensure good satisfaction of the constraint but not too large so as to significantly impact the conditioning of the system.}
\section{Laplace-Beltrami problem}
\label{sec:Laplace-Beltrami}
The governing partial differential equation which we want to solve to illustrate fundamental features of  subdivision surfaces is given by
\begin{equation}
-\Delta_{\Gamma} u = f \quad \text{on} \quad \Gamma,
\label{eq:Poisson}
\end{equation}
where $\Gamma$ is a two dimensional manifold (with outward unit normal vector $\mathbf{n}$) in three dimensional space $\mathbb{R}^3$ and $\Delta_{\Gamma}(\bullet)$ is the Laplace-Beltrami operator (also called surface Laplacian operator). The Dirichlet boundary condition is expressed in~\eqref{eq:Poisson_bc}.
We will use a manufactured solution to compute against the approximate solution. The Laplace-Beltrami operator is defined by
\begin{equation}
\Delta_{\Gamma} (\bullet) = \nabla_{\Gamma}\cdot\nabla_{\Gamma}(\bullet),
\end{equation}
where $\nabla_{\Gamma}(\bullet)$ is the surface gradient operator defined by
\begin{equation}
\nabla_{\Gamma}(\bullet) = [\mathbf{I} - \mathbf{n} \otimes \mathbf{n}] \cdot \nabla(\bullet).
\end{equation}
{
Hence the surface gradient of a scalar function $v$ can be calculated as the spatial gradient subtracted by its normal part as
\begin{equation}
\nabla_{\Gamma} v = \nabla v - \mathbf{n}[\mathbf{n} \cdot \nabla v],
\end{equation}
where $\nabla(\bullet)$ is the spatial gradient operator. Hence the surface Laplacian of $v$ is given by
\begin{equation}
\Delta_{\Gamma} v = \Delta v - \mathbf{n}\cdot[\nabla^2 v \cdot \mathbf{n}] - [\mathbf{n} \cdot \nabla v] \left[\nabla \cdot \mathbf{n} - \mathbf{n}\cdot [\nabla\mathbf{n} \cdot \mathbf n]\right],
\label{eq:Laplacian-Beltrami}
\end{equation}
where $ \nabla^2 v$ is the Hessian matrix of $v$, and $\nabla\mathbf{n}$ is the gradient of the normal vector, which is arranged in a matrix as
\begin{equation}
\nabla \mathbf n =
\left[ 
\begin{array}{ccc}
\frac{\partial n_1}{\partial x_1} &\frac{\partial n_2}{\partial x_1} & \frac{\partial n_3}{\partial x_1}\\
 \frac{\partial n_1}{\partial x_2}& \frac{\partial n_2}{\partial x_2}& \frac{\partial n_3}{\partial x_2}\\
 \frac{\partial n_1}{\partial x_3}&\frac{\partial n_2}{\partial x_3} &\frac{\partial n_3}{\partial x_3}
\end{array}
\right].
\end{equation}
We define the total curvature at a surface point $\mathbf{x} \in \Gamma$ as the surface divergence of the normal, that is $c(\mathbf{x}) := \nabla_\Gamma\cdot \mathbf{n}$.
For a given manufactured solution $u^m$, the right hand side of Equation~\eqref{eq:Poisson} can thus be computed as
\begin{equation}
   f(\mathbf{x}) = - \Delta u^m(\mathbf{x}) + \mathbf{n}(\mathbf{x})\cdot[\nabla^2 u^m(\mathbf{x})\cdot \,\mathbf{n}(\mathbf{x})] + c(\mathbf{x}) [\mathbf{n}(\mathbf{x}) \cdot \nabla u^m(\mathbf{x})], \quad \mathbf{x} \in \Gamma.
\label{eq:bt_rhs}
\end{equation}
}
\section{Galerkin formulation}
\label{sec:FE}
The weak formulation of problem~\eqref{eq:Poisson} is
\begin{equation}
\int_{\Gamma} \nabla_{\Gamma} u \cdot  \nabla_{\Gamma} v \, \mathrm{d} \Gamma = \int_{\Gamma} f v  \,\mathrm{d} \Gamma, \quad \forall v \in H^1_0(\Gamma),
\label{eq:weak_form}
\end{equation}
where $v$ is an admissible test function. The weak formulation is partitioned into $n_e$ number of elements, as
\begin{equation}
\sum_{k= 1}^{n_{e}}\int_{\Gamma_{e}} \nabla_{\Gamma} u \cdot \nabla_{\Gamma} v  \,\mathrm{d} \Gamma = \sum_{k= 1}^{n_{e}}\int_{\Gamma_{e}} f v  \,\mathrm{d} \Gamma.
\end{equation}
Discretising $v$, $\nabla u$ and $\nabla v$ using the Catmull-Clark basis functions ${\mathbf{N}}$  given in Equation~\eqref{eq:cc_bases_tp} produces
\begin{align}
v &= \sum_{A = 1}^{n_b} N_{A} v_{A}, \nonumber \\
\nabla_{\Gamma} u &= \sum_{A = 1}^{n_b} \nabla_{\Gamma} N_A u_A = \sum_{A = 1}^{n_b} \frac{\partial N_A}{\partial \boldsymbol{\xi}}\cdot\mathbf{J}^{-1} u_{A}, \nonumber \\	
\nabla_{\Gamma} v &= \sum_{A = 1}^{n_b} \nabla_{\Gamma} N_A v_A =  \sum_{A = 1}^{n_b} \frac{\partial N_A}{\partial \boldsymbol{\xi}}\cdot\mathbf{J}^{-1} v_{A},
\label{eq:discretisation_variables}
\end{align}
where $\mathbf{J}$ is the surface Jacobian for the manifold,  given in a matrix form as
\begin{equation}
\mathbf{J} = \frac{\partial \mathbf{x}}{\partial \boldsymbol{\xi}} =
\left[
\begin{array}{cc}
\frac{\partial x_1}{\partial \xi} & \frac{\partial x_1}{\partial \eta} \\
\frac{\partial x_2}{\partial \xi} & \frac{\partial x_2}{\partial \eta} \\
\frac{\partial x_3}{\partial \xi} & \frac{\partial x_3}{\partial \eta}
\end{array}
\right].
\end{equation}
For details on the computation of $\mathbf{J}^{-1}$ see~\cite{peterson2005mapped} and for a discussion of superficial tensors such as $\mathbf{J}$ in the context of Laplace-Beltrami equation, see~\cite{gurtin2002interface}.  If the element contains an extraordinary vertex, the shape functions $N_A$ are replaced by $\hat{N}_A$ in Equation~\eqref{eq:subd_basis}. The surface gradient of the shape functions is computed as
\begin{equation}
\nabla_\Gamma \hat{N}_A = \frac{\partial  \hat{N}_A}{\partial \bar{\boldsymbol{\xi}}} \cdot \frac{\partial \bar{\boldsymbol{\xi}}}{\partial \boldsymbol{\xi}}\cdot \mathbf{J}^{-1}
\end{equation}
and
\begin{equation}
\mathbf{J} =\frac{\partial \mathbf{x}}{\partial \bar{\boldsymbol{\xi}}} \cdot \frac{\partial \bar{\boldsymbol{\xi}}}{\partial \boldsymbol{\xi}} =  \sum_{A=0}^{2\kappa+7} \frac{\partial \hat{N}_A}{\partial \bar{\boldsymbol{\xi}}} \, \mathbf{P}_{A} \cdot \frac{\partial \bar{\boldsymbol{\xi}}}{\partial \boldsymbol{\xi}}.
\end{equation}
Integrating the discrete problem using Gauss quadrature, the system of equations~\ref{eq:weak_form} becomes
\begin{equation}
\underbrace{\left[\bigA_{e= 1}^{n_{e}}\sum_{i = 1}^{n_q}\mathbf{D}^e(\mathbf{\boldsymbol{\xi}}^i) |\mathbf{J}(\mathbf{x}(\boldsymbol{\xi}^{i}))| w^{i}\right]}_{\ary K}\ \ary u=\underbrace{ \bigA_{e= 1}^{n_{e}}\sum_{i = 1}^{n_q} f(\mathbf{x}(\boldsymbol{\xi}^{i})) \mathbf{N}^e(\boldsymbol{\xi}^{i})  |\mathbf{J}(\mathbf{x}(\boldsymbol{\xi}^{i}))| w^{i}}_{\ary f} ,
\label{eq:fem_formulation_discretisation}
\end{equation}
where $\bigA$ is the assembly operator and
\begin{equation}
\mathbf{D}^e(\boldsymbol{\xi}) =
\left[
\begin{array}{cccc}
\nabla_\Gamma N_1 (\boldsymbol{\xi}) \cdot \nabla_\Gamma N_1 (\boldsymbol{\xi})  & \nabla_\Gamma N_1 (\boldsymbol{\xi}) \cdot \nabla_\Gamma N_2 (\boldsymbol{\xi}) & \cdots & \nabla_\Gamma N_1 (\boldsymbol{\xi}) \cdot \nabla_\Gamma N_{n_b} (\boldsymbol{\xi})\\
\nabla_\Gamma N_2 (\boldsymbol{\xi}) \cdot \nabla_\Gamma N_1 (\boldsymbol{\xi})  & \nabla_\Gamma N_2 (\boldsymbol{\xi}) \cdot \nabla_\Gamma N_2 (\boldsymbol{\xi}) & \cdots & \nabla_\Gamma N_2 (\boldsymbol{\xi}) \cdot \nabla_\Gamma N_{n_b} (\boldsymbol{\xi})\\
\vdots & \vdots & \ddots & \vdots\\
\nabla_\Gamma N_{n_b} (\boldsymbol{\xi}) \cdot \nabla_\Gamma N_1 (\boldsymbol{\xi})  & \nabla_\Gamma N_{n_b} (\boldsymbol{\xi}) \cdot \nabla_\Gamma N_2 (\boldsymbol{\xi}) & \cdots & \nabla_\Gamma N_{n_b} (\boldsymbol{\xi}) \cdot \nabla_\Gamma N_{n_b} (\boldsymbol{\xi})\\
\end{array}
\right].
\end{equation}
$n_q$ is the number of quadrature points in each element, $w_i$ is the weight for $i^\text{th}$ quadrature point, $n_e$ is the number of elements and $n_b$ is the number of basis functions of the element. The basis functions $\boldsymbol{N}^e$ are replaced by $\hat{\boldsymbol{N}}^e$ if the element $e$ contains an extraordinary vertex. In this case, the basis functions are not differentiable and their derivatives approach positive infinity when points are close to the extraordinary vertex~(see Remark~\ref{rm:singular}). Thus $|\mathbf{J}|$ approaches positive infinity at extraordinary vertices. Errors result if quadrature is adopted to integrate the contributions from element containing extraordinary vertices.

The discrete system of equations to solve is thus given by
\begin{equation}
\ary K \ary u = \ary f.
\label{eq:system_equation}
\end{equation}

\section{Numerical results}
\label{sec:numerical}
{
A `patch test'~\cite{zienkiewicz1997finite} on a two-dimensional plate is first presented to assess the consistency and stability of the proposed formulation in a simplified setting.  Then, the Laplace-Beltrami equation is solved on both cylindrical and hemispherical surfaces. Convergence studies are conducted. The influence of extraordinary vertices is also investigated. All numerical results are computed using an open source finite element library: deal.II~\cite{dealII91,BangerthHartmannKanschat2007}.
}
\subsection{`Patch test'}
\label{sec:patch}
\begin{figure}
\centering
\begin{subfigure}[]{0.51\linewidth}
\centering
  \includegraphics[width=0.75\linewidth]{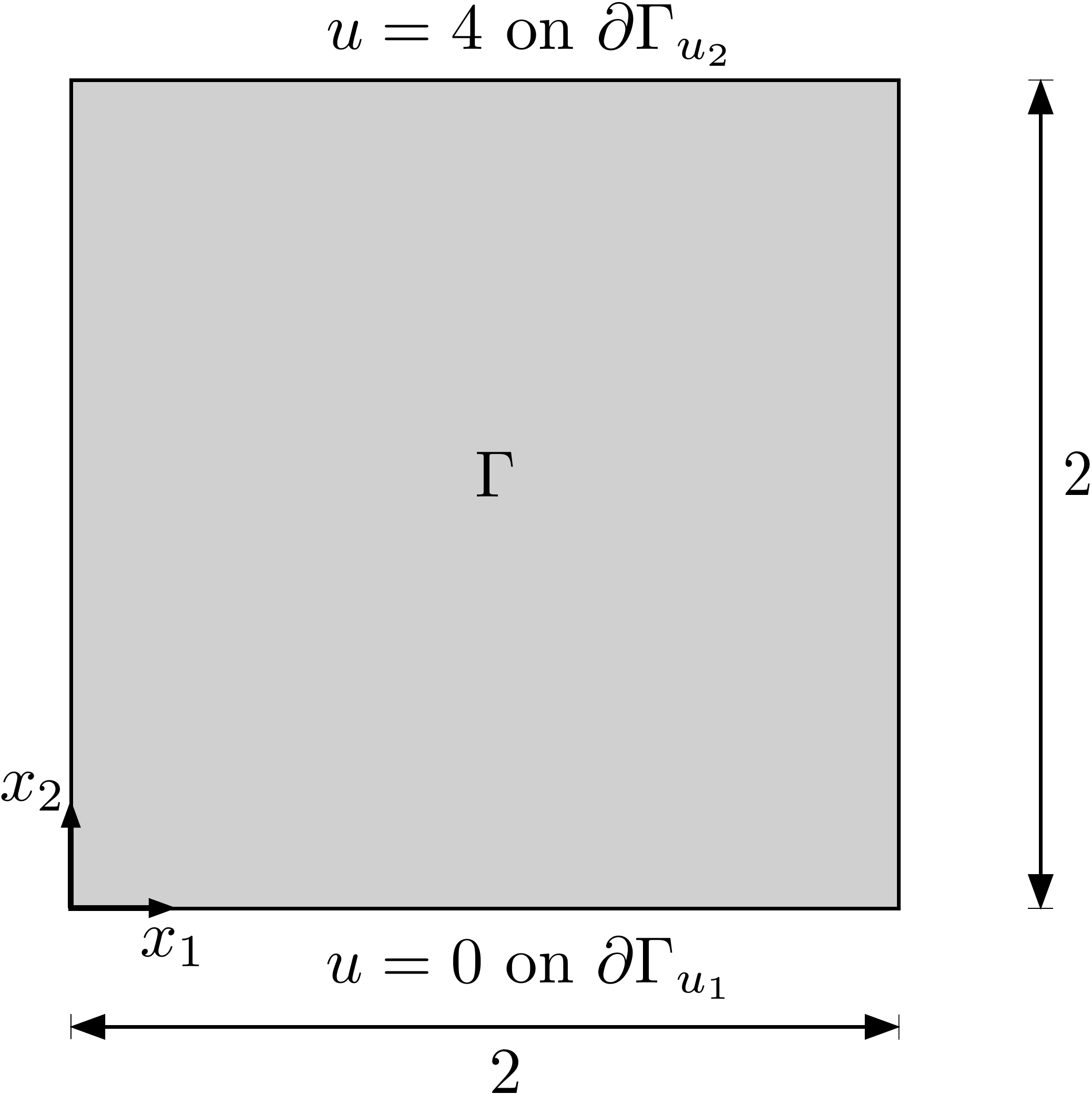}
	\caption{Plate problem setup.}
	\label{fig:plate_problem}
\end{subfigure}

\begin{subfigure}[]{0.47\linewidth}
\centering
  \includegraphics[width=\linewidth]{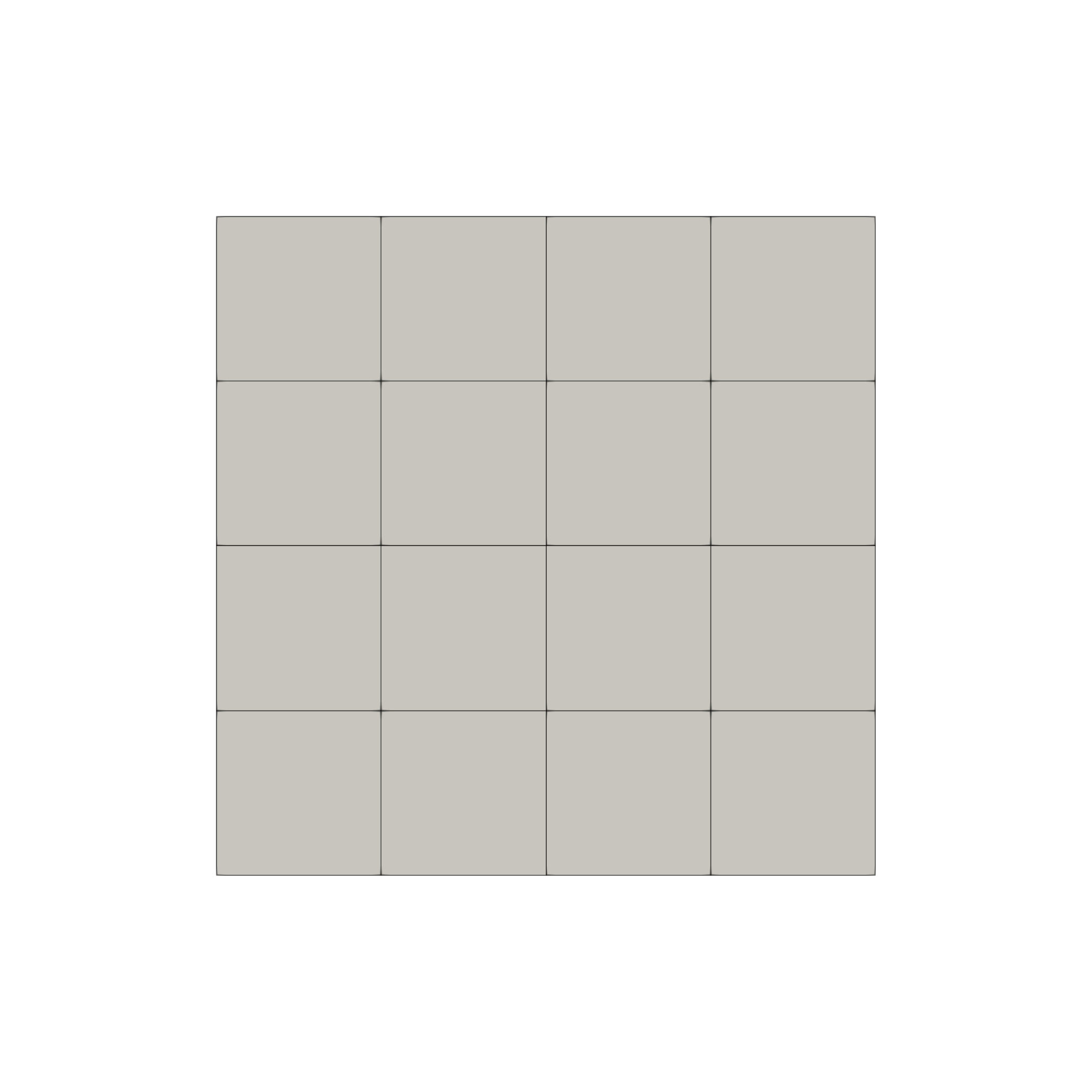}
	\caption{A regular mesh (Mesh 1).}
	\label{fig:plate_reg_grid}
\end{subfigure}
\begin{subfigure}[]{0.47\linewidth}
\centering
  \includegraphics[width=\linewidth]{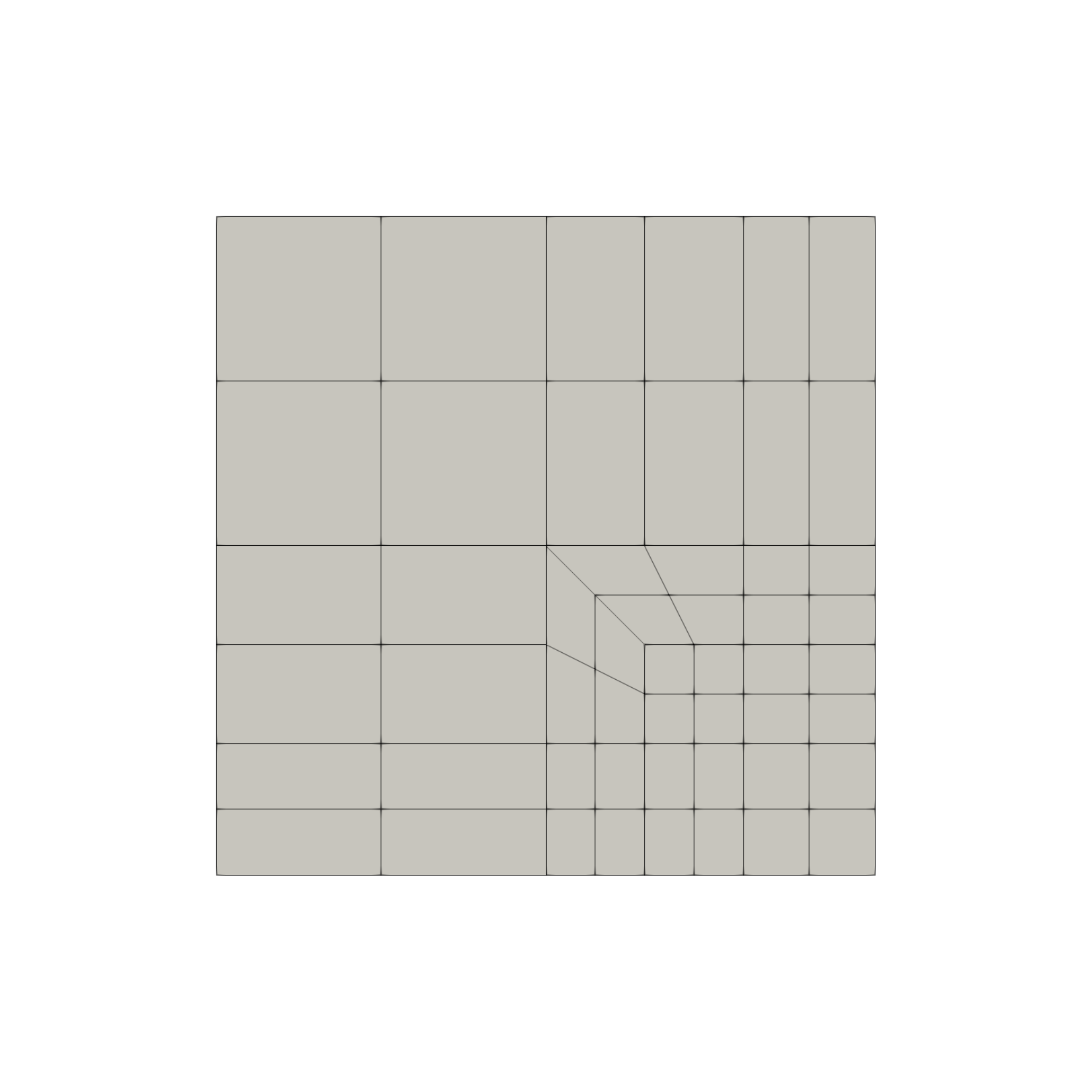}
	\caption{A mesh with extraordinary vertices (Mesh 2).}
	\label{fig:plate_irreg_grid}
\end{subfigure}
\caption{Schematic of the patch test on a plate.}
\label{fig:plate}
\end{figure}

The `patch test' is performed on a  two dimensional flat plate where the Laplace-Beltrami operator reduces to the Laplace operator.  The problem proposed in Section~\ref{sec:Laplace-Beltrami} reduces to the Poisson problem expressed given by
\begin{equation}
-\Delta u = f \quad \text{on} \quad \Gamma \in \mathbb{R}^2.
\end{equation}
This partial differential equation is solved on the square plate shown in Figure~\ref{fig:plate_problem} with the essential boundary conditions
\begin{equation}
\left\{
\begin{split}
u = 0 \quad \text{for} \quad x_2 = 0 \quad (\partial \Gamma_{u_1})\\
u = 4 \quad \text{for} \quad x_2 = 2 \quad (\partial \Gamma_{u_2})
\end{split}
\right..
\end{equation} 
{The essential boundary conditions are imposed using the penalty method. Natural homogeneous boundary conditions are applied on the remaining two edges of the plate.} Four different manufactured functions for $f$ are used.  The functions, analytical solutions for $u$ and their gradients $\frac{\partial u}{\partial x_2}$ are given in Table~\ref{tab:test_functions}. We investigate both a regular and an irregular mesh. The regular mesh is a $4 \times 4$ element patch without extraordinary vertices as shown in Figure~\ref{fig:plate_reg_grid}. {In all of the tests, a geometry error is absent. }
\begin{table}[]
\centering
\begin{tabular}{@{}llll@{}}
\toprule
Test No. &$ f $                            & $u$     & $\frac{\partial u}{\partial x_2} $  \\ \midrule
1                & $0$                            & $2x_2$ & $2$        \\
2                & $1$                            & $-\frac{1}{2}{x_2}^2 + 3x_2$     & $-x_2+3$       \\
3                & $x_2$                         & $-\frac{1}{6}{x_2}^3 + \frac{8}{3} x_2$     & $-\frac{1}{2}{x_2}^2 + \frac{8}{3}$        \\
4                & $\pi\sin(\pi x_2)$ &  $\frac{1}{\pi}\sin(\pi x_2)+2x_2$    & $\cos(\pi  x_2)+2$    \\
\bottomrule
\end{tabular}
\caption{Four test case functions for the plate problem. {
Test 1 has no right-hand side term, thus the analytical solution $u$ is linear and its gradient is a constant. The analytical solutions for Tests 2 and 3 are quadratic and cubic, respectively, and their gradients are linear and quadratic, respectively. Test 4 has a sine function as the right-hand side term which gives a cosine function as the gradient of the analytical solution.}}
	\label{tab:test_functions}
\end{table}

\begin{figure}
\centering
\begin{subfigure}{0.4\linewidth}
\centering
\includegraphics[width=	\linewidth]{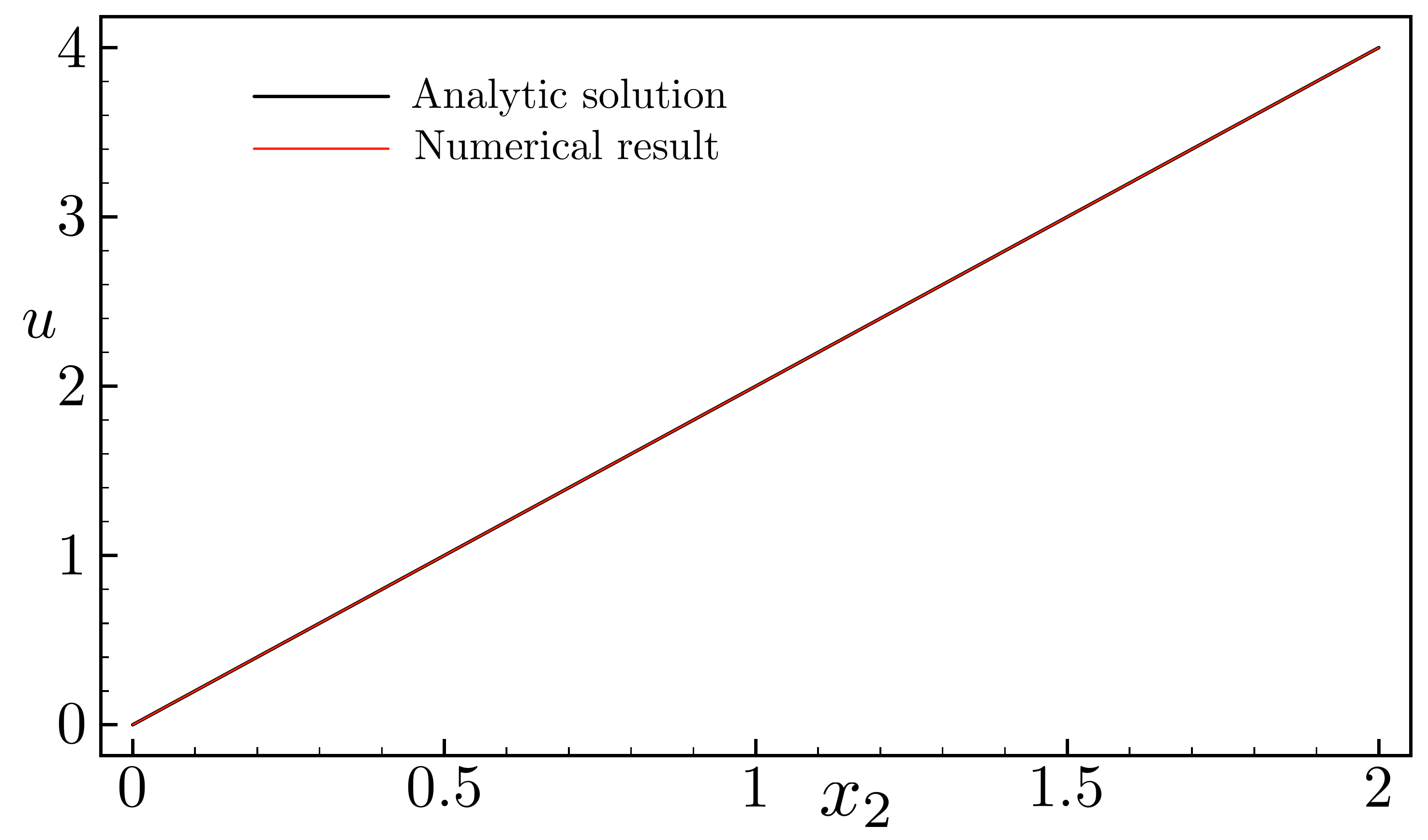}
\caption{Solution for Test 1.}
\label{fig:plate_f=0_u}
\end{subfigure}
\begin{subfigure}{0.4\linewidth}
\centering
\includegraphics[width=	\linewidth]{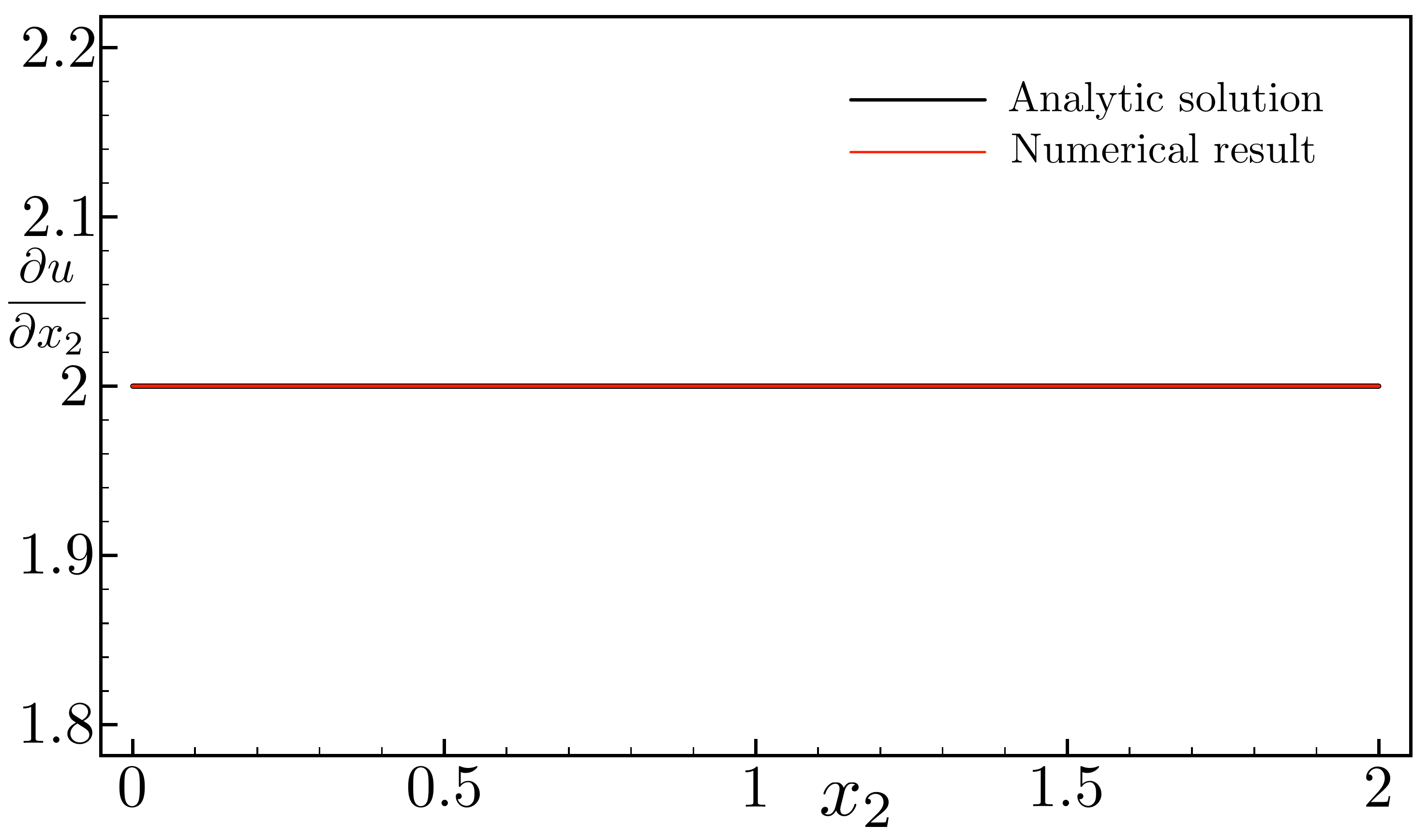}
\caption{Gradient for Test 1.}
\label{fig:plate_f=0_plot}
\end{subfigure}
\begin{subfigure}{0.4\linewidth}
\centering
\includegraphics[width=	\linewidth]{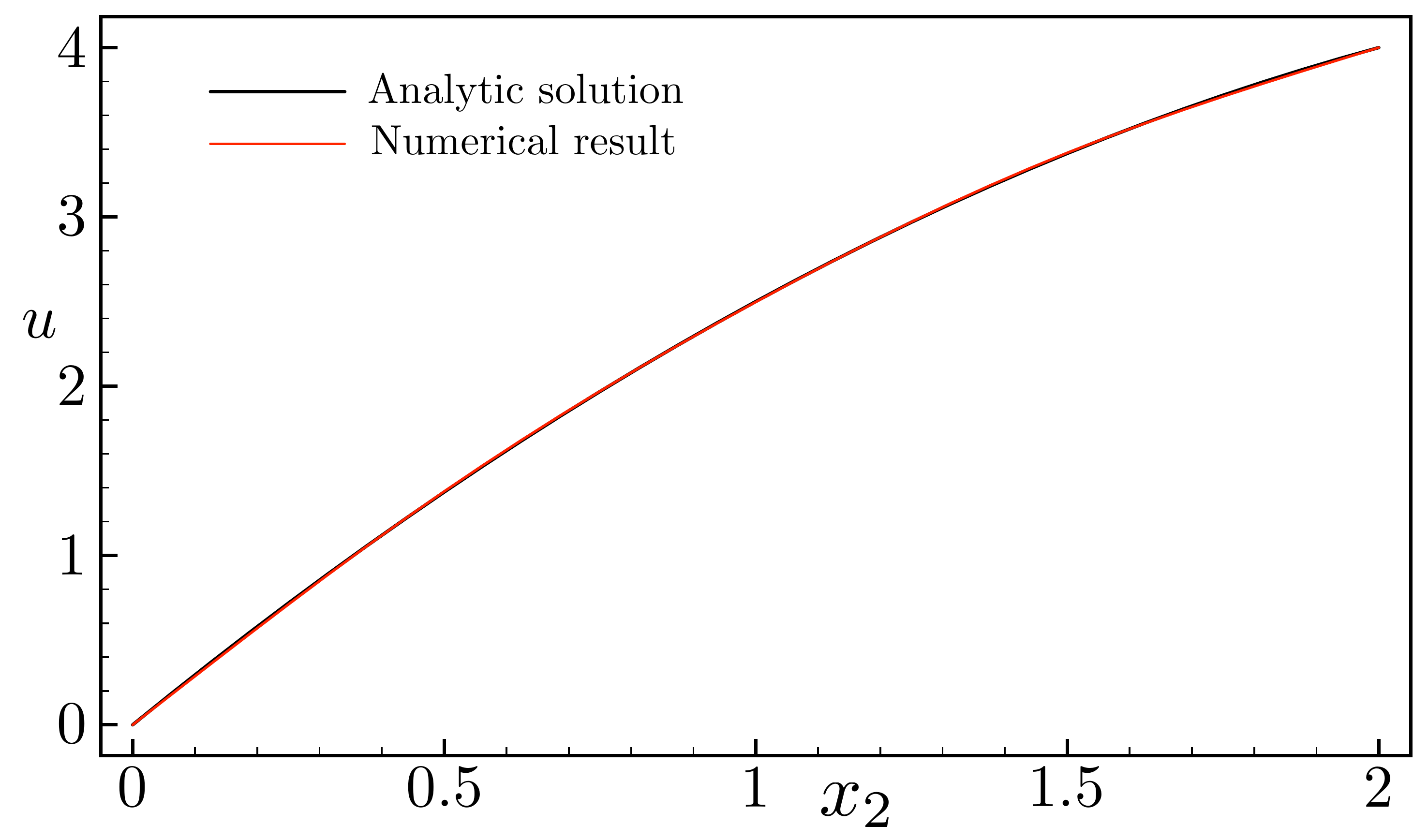}
\caption{Solution for Test 2.}
\label{fig:plate_f=1_u}
\end{subfigure}
\begin{subfigure}{0.4\linewidth}
\centering
\includegraphics[width=	\linewidth]{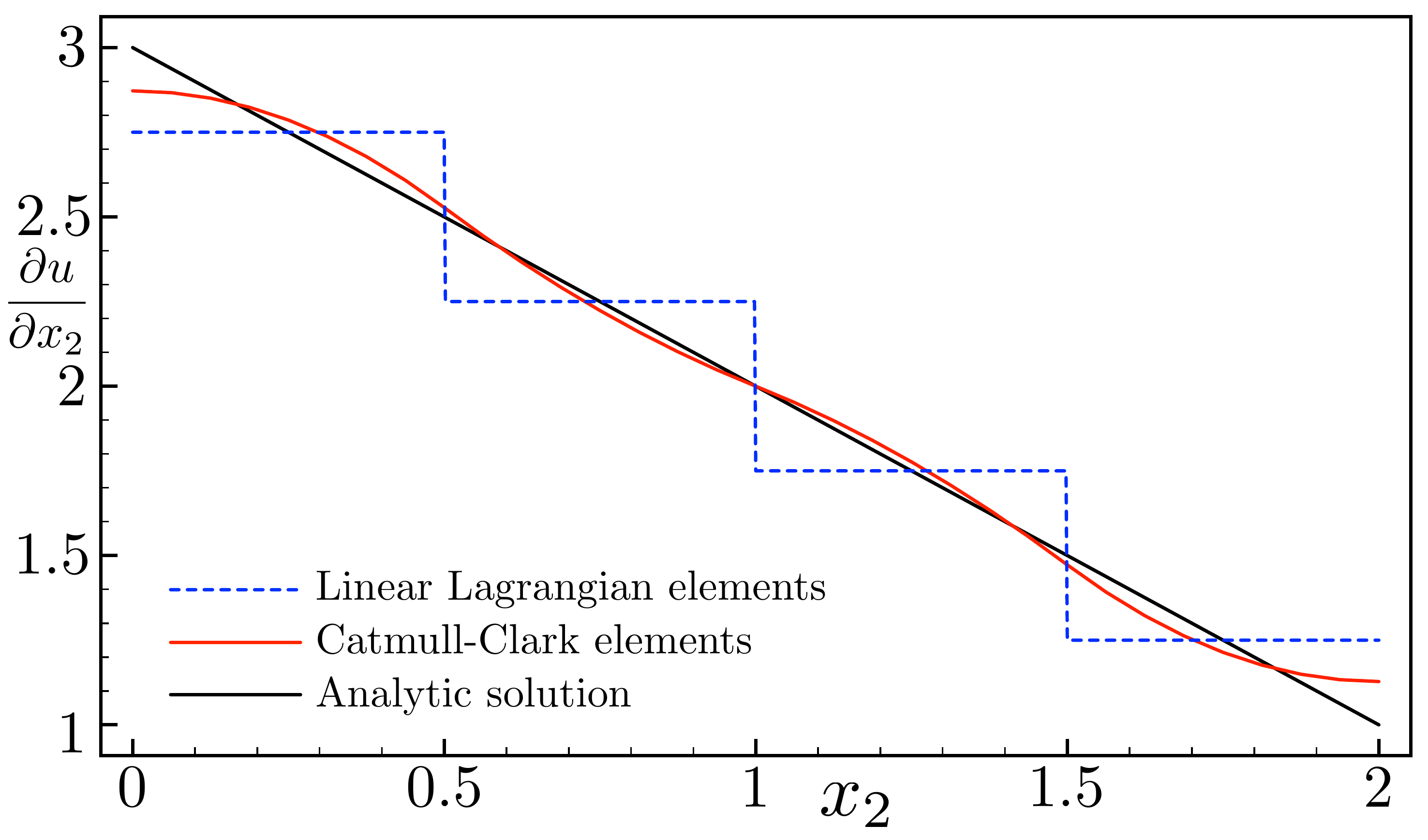}
\caption{Gradient for Test 2.}
\label{fig:plate_f=1_plot}
\end{subfigure}
\begin{subfigure}{0.4\linewidth}
\centering
\includegraphics[width=	\linewidth]{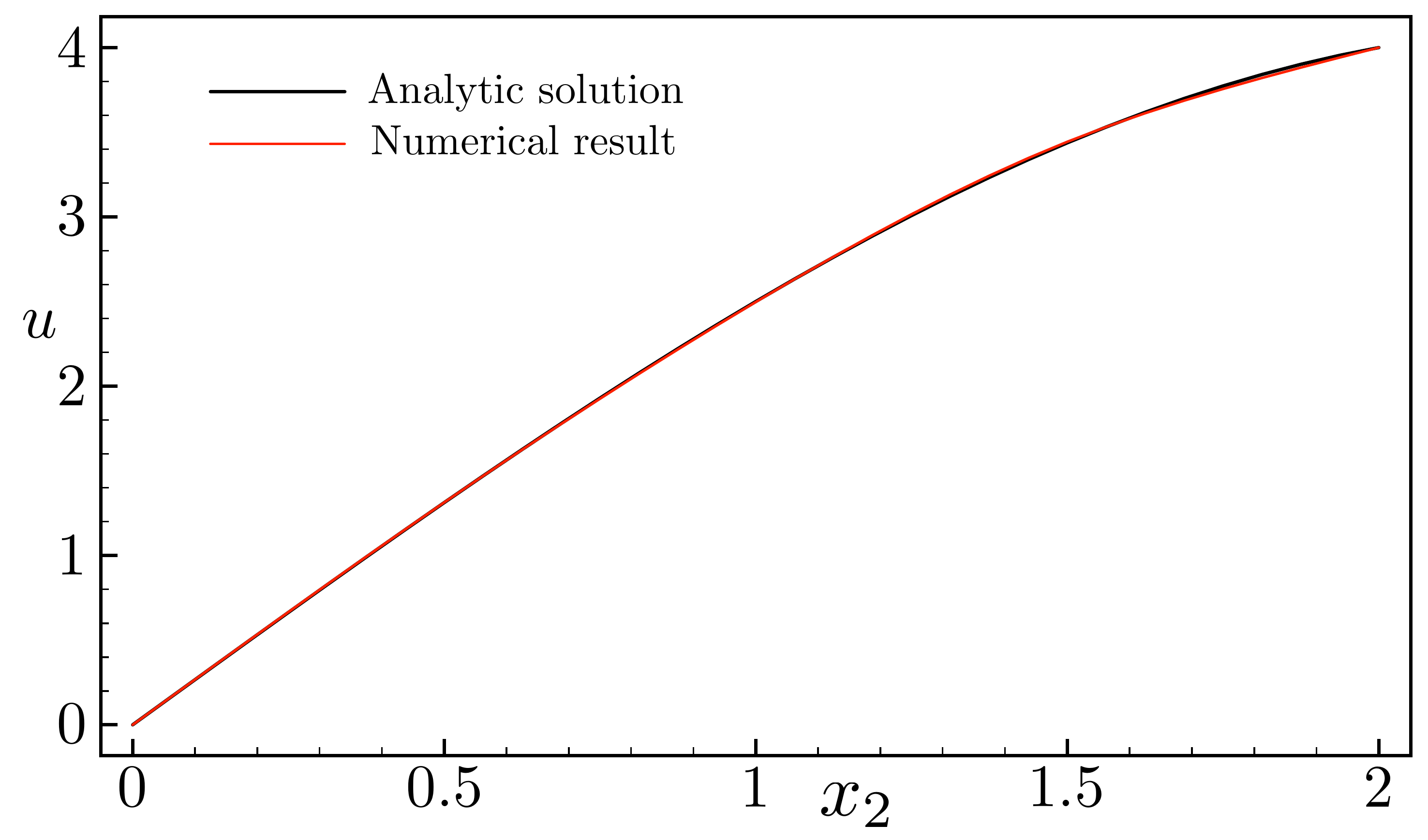}
\caption{Solution for Test 3.}
\label{fig:plate_f=y_u}
\end{subfigure}
\begin{subfigure}{0.4\linewidth}
\centering
\includegraphics[width=	\linewidth]{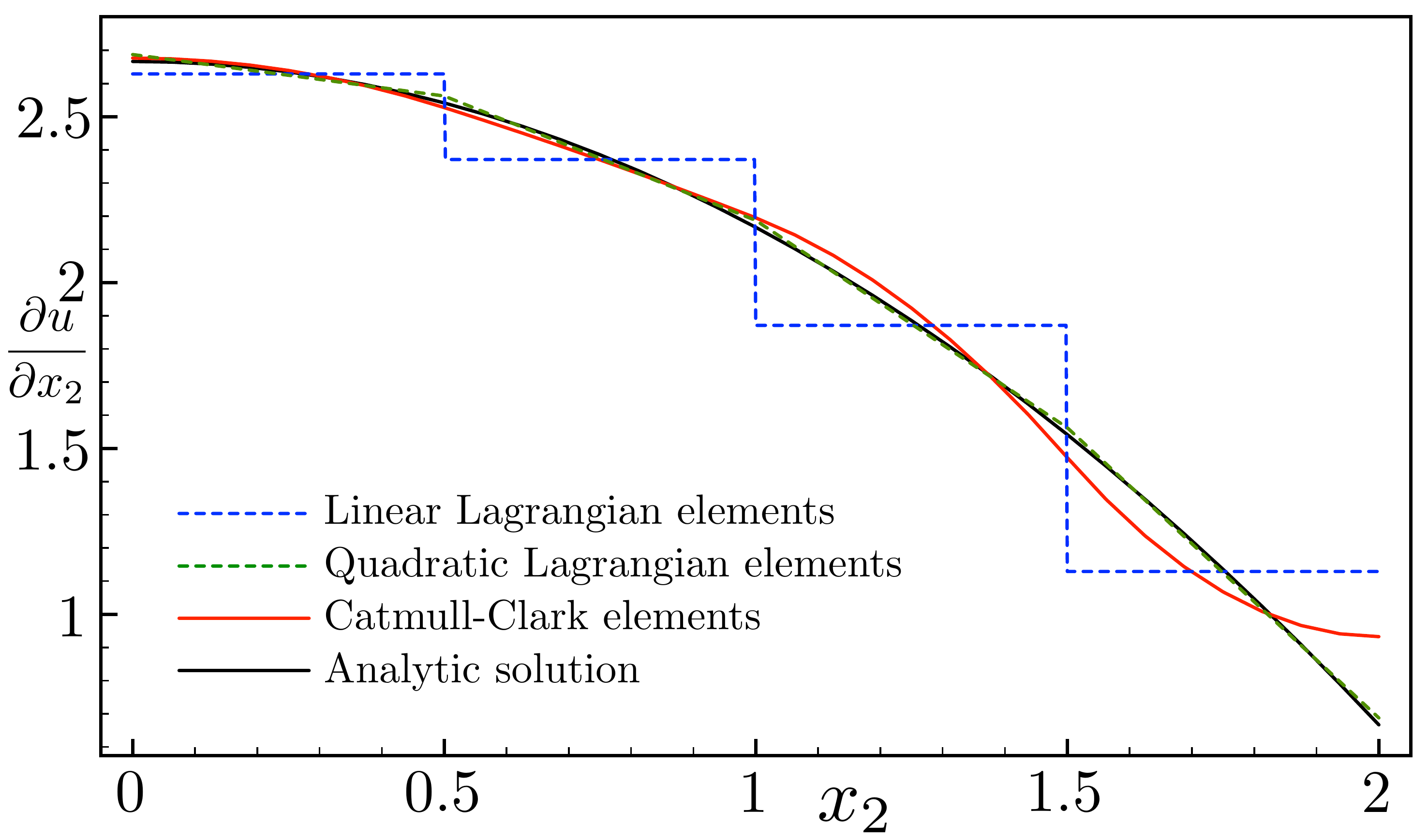}
\caption{Gradient for Test 3.}
\label{fig:plate_f=y_plot}
\end{subfigure}
\begin{subfigure}{0.4\linewidth}
\centering
\includegraphics[width=	\linewidth]{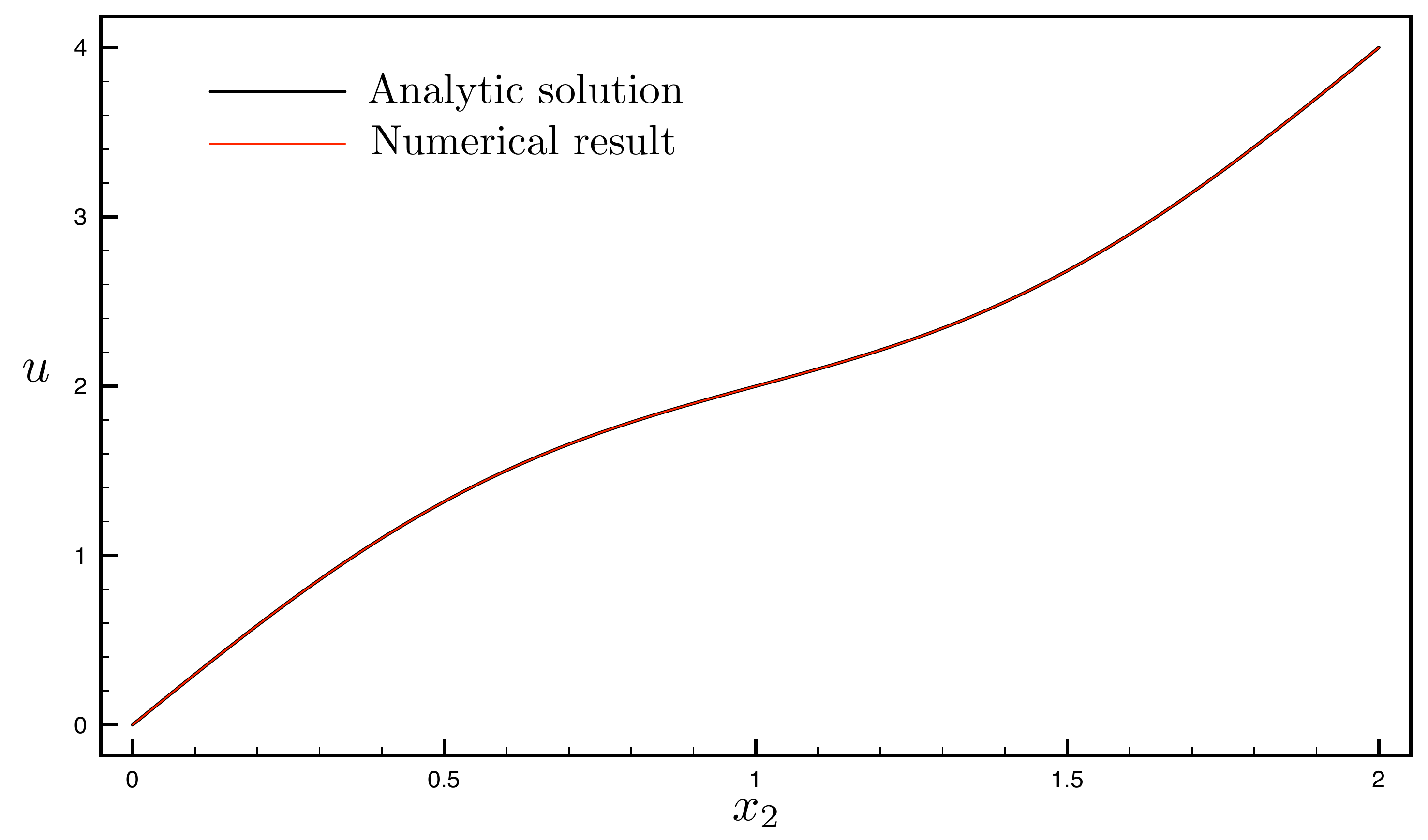}
\caption{Solution for Test 4.}
\label{fig:plate_f=pisinpiuy_u}
\end{subfigure}
\begin{subfigure}{0.4\linewidth}
\centering
\includegraphics[width=	\linewidth]{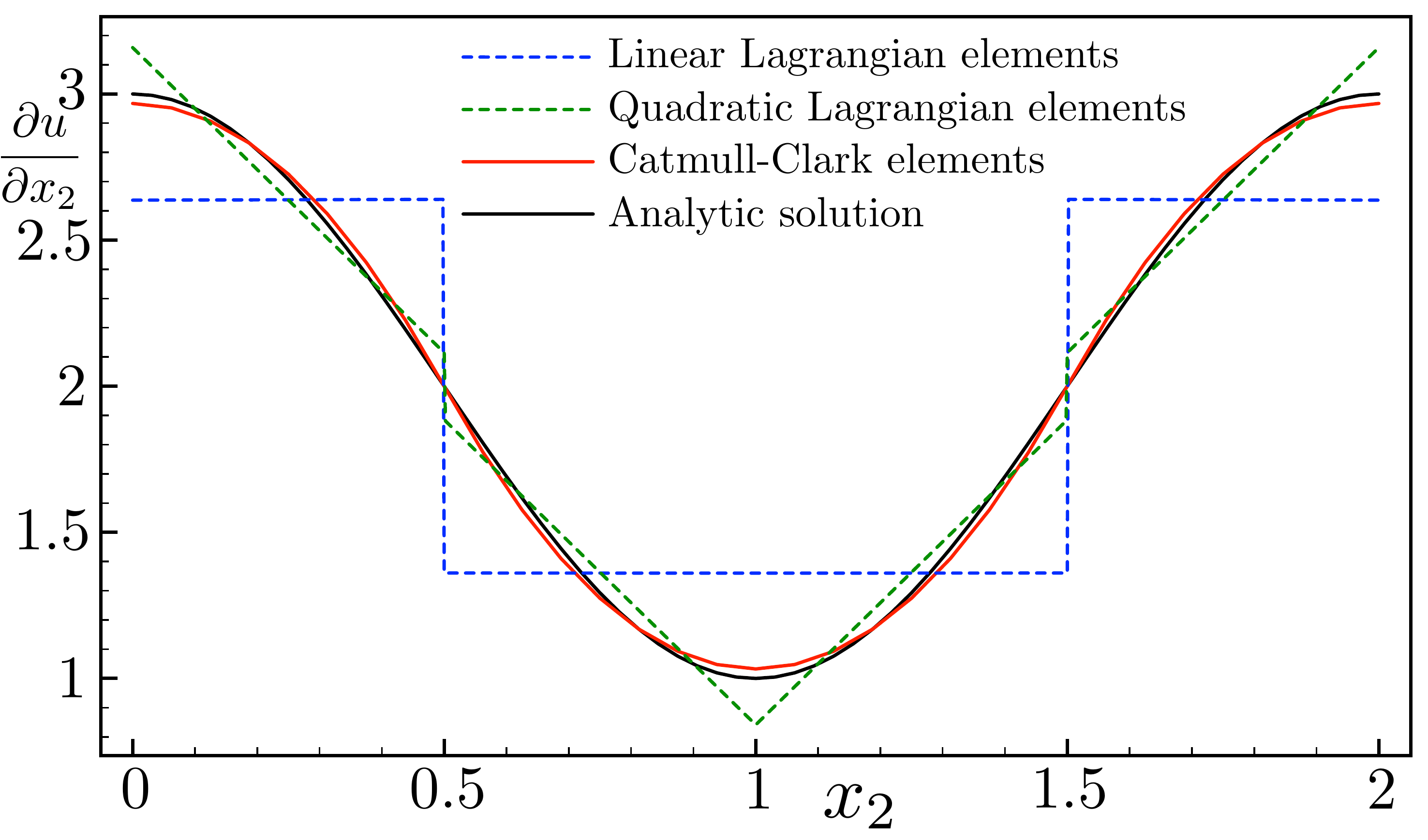}
\caption{Gradient for Test 4.}
\label{fig:plate_f=pisinpiy_plot}
\end{subfigure}
\caption{Solution $u$ and the gradient $\frac{\partial u}{\partial x_2}$ plotted along the line $x_1 = 1$ for the plate test. The numerical results are compared to the analytical solutions for tests 2, 3 and 4. {Mesh 1 is used for all tests.}}
\label{fig:gradient_plot}
\end{figure}

For Test 1, the right hand side $f = 0$ so that $\frac{\partial u}{\partial x_2} = 2$. Solving the equation using the proposed Catmull-Clark subdivision method, the numerical result $u_h$ is exactly $2$ everywhere as shown in Figure~\ref{fig:plate_f=0_plot}. Thus passes the consistency test and the eigenvalue of the system matrix are all positive and non-zero after application of the essential boundary conditions.
The gradient $\frac{\partial u}{\partial x_2}$ for Test 2 and 3 are linear and quadratic respectively. Recall that when interpolating functions in elements with edges on physical boundaries, the basis functions are {modified}, see Equations~\eqref{eq:mirroring} and~\eqref{eq:boundary_basis}. In other words, the gradients of the function $u$ are expected to be constant at boundaries. Figures~\ref{fig:plate_f=0_u},~\ref{fig:plate_f=1_u} and~\ref{fig:plate_f=y_u} show the numerical results for these tests. The results are smooth and capture the analytical solutions well.  Figures~\ref{fig:plate_f=1_plot} and~\ref{fig:plate_f=y_plot} compare the numerical results of $\frac{\partial u}{\partial x_2}$ to the analytical solution for Test 2 and 3. The Catmull-Clark subdivision method is also compared to linear and quadratic Lagrangian finite element methods. There is a substantial error in both boundary regions in Test 2 for Catmull-Clark subdivision method. This is because the method imposes the gradient to be constant at both boundaries. The numerical result of the Catmull-Clark subdivision method in Test 3 has a substantial error in the region close to the top boundary ($x_2=2$) but captures the gradient in the region close to the bottom boundary ($x_2=0$) well because the analytic solution for the gradient in the bottom boundary region is near-constant. These errors at the boundaries will pollute the numerical result in the interior of the domain, which will reduce the convergence rate. The gradients approximated by the linear and quadratic Lagrangian finite elements are piecewise constant and piecewise linear, respectively. The results of the Catmull-Clark subdivision methods for these two tests lies between the linear and quadratic Lagrangian elements. The gradient $\frac{\partial u}{\partial x_2}$ in Test 4 is a cosine function which is non-polynomial and it behaves as a constant in both boundary regions shown in Figure~\ref{fig:plate_f=pisinpiy_plot}. The Lagrangian elements only possess $C^0$ continuity across elements and their gradients hence have jumps between elements. The Catmull-Clark subdivision elements capture the gradients of the given function better as they are $C^1$ smooth.

\begin{figure}
\centering
\begin{subfigure}{0.49\linewidth}
\centering
	\includegraphics[width=\linewidth]{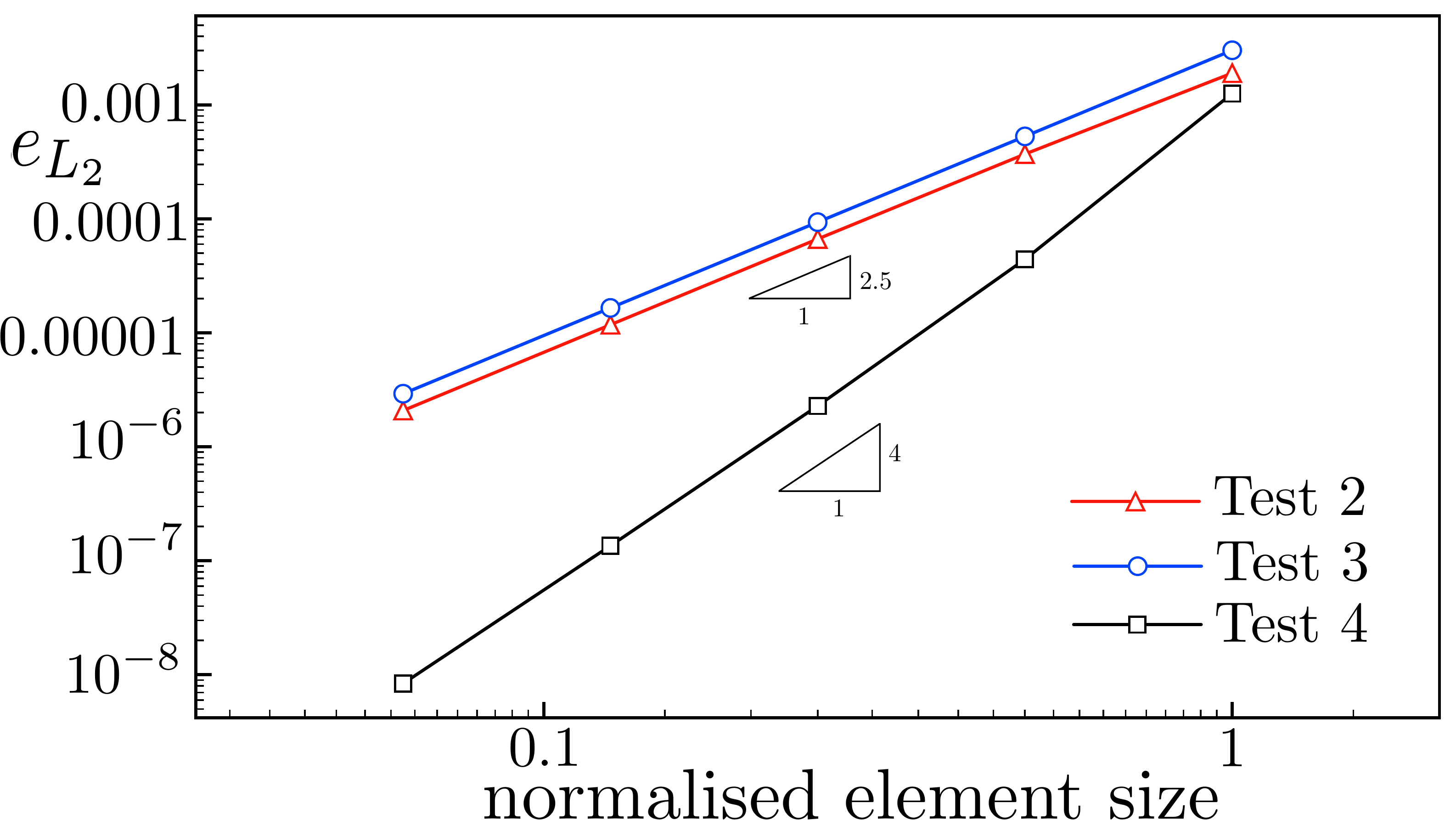}
	\caption{Normalised global $L_2$ error (Mesh 1).}
	\label{fig:plate_converge_l2}
\end{subfigure}
\begin{subfigure}{0.49\linewidth}
\centering
	\includegraphics[width=\linewidth]{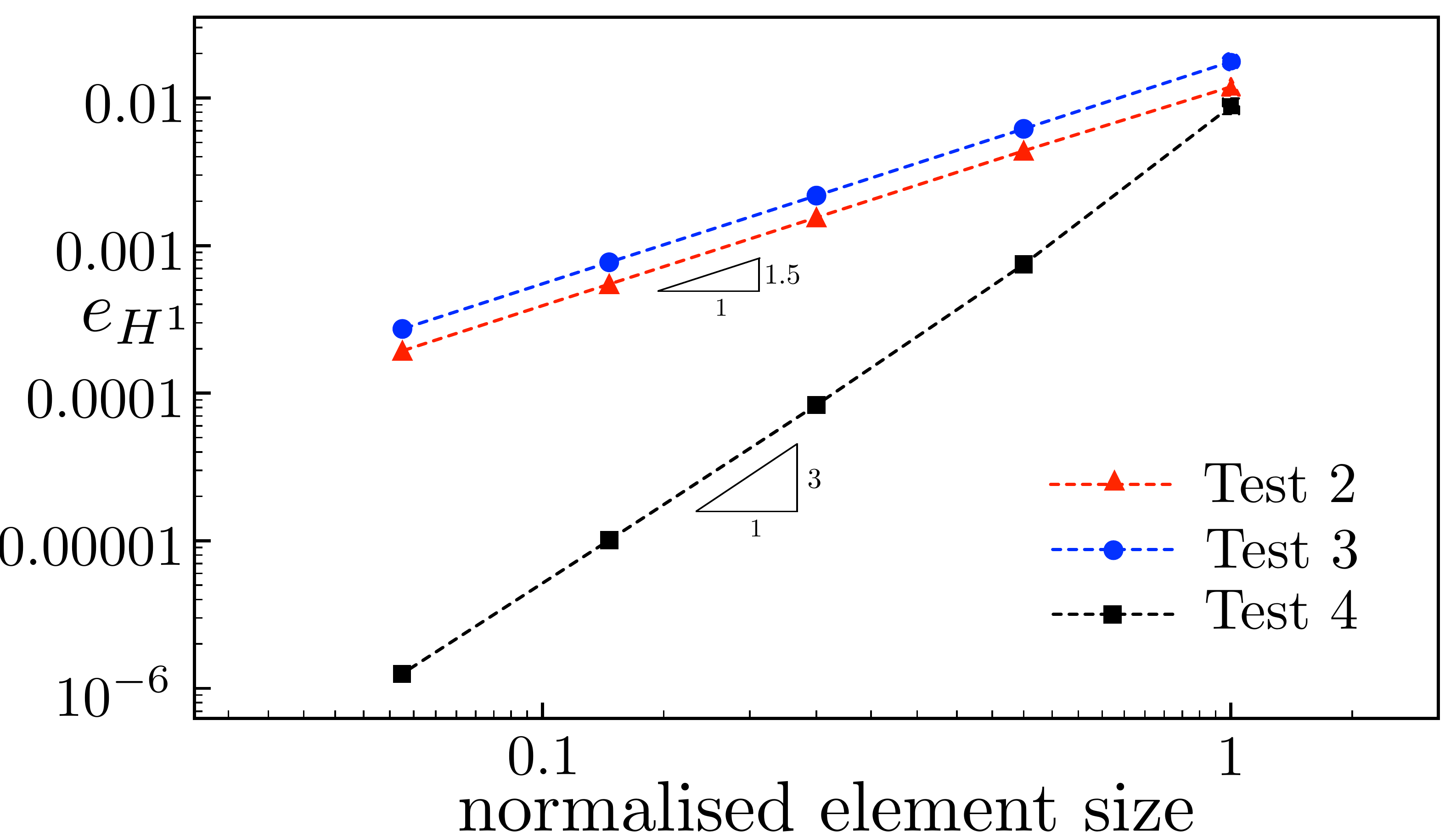}
	\caption{ Normalised global $H^1$ error (Mesh 1).}
	\label{fig:plate_converge_h1}
\end{subfigure}
\begin{subfigure}{0.49\linewidth}
\centering
	\includegraphics[width=\linewidth]{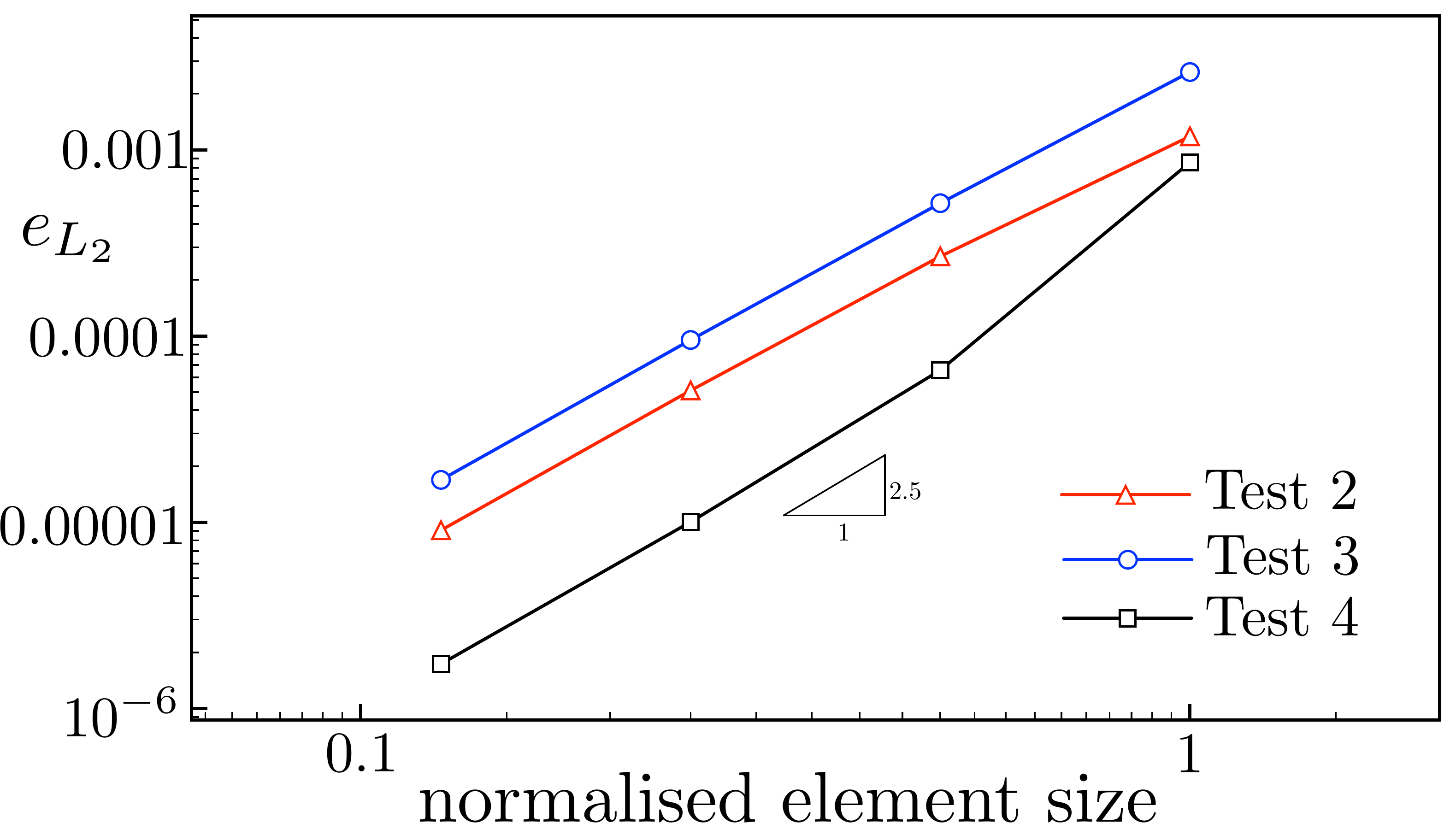}
	\caption{Normalised global $L_2$ error (Mesh 2).}
	\label{fig:plate_converge_l2_irreg}
\end{subfigure}
\begin{subfigure}{0.49\linewidth}
\centering
	\includegraphics[width=\linewidth]{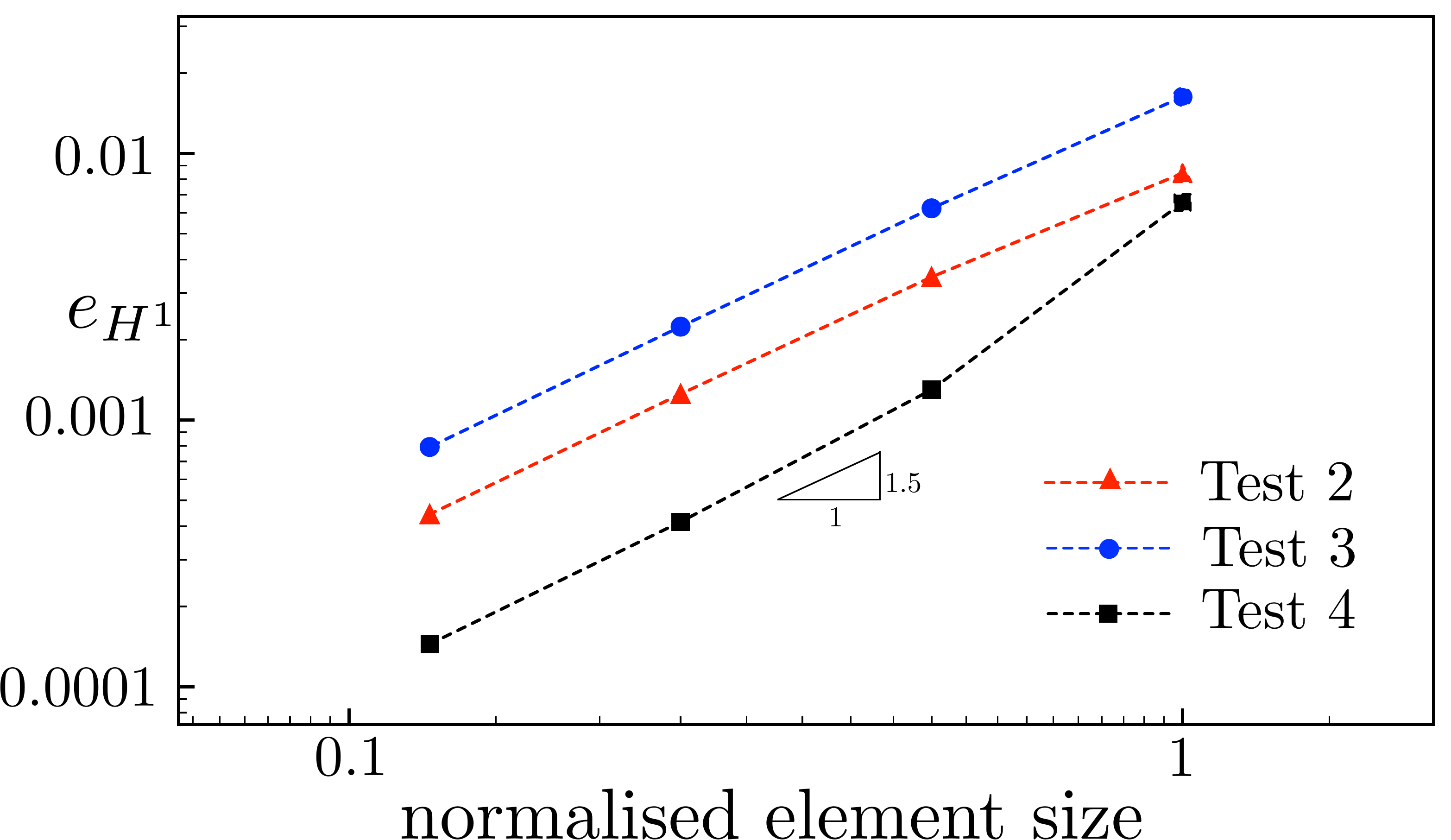}
	\caption{ Normalised global $H^1$ error (Mesh 2).}
	\label{fig:plate_converge_h1_irreg}
\end{subfigure}
\caption{Convergence study for Tests 2, 3 and 4 using the regular mesh (Mesh 1) and the mesh with an extraordinary vertex (Mesh 2).}
\label{fig:plate_converge}
\end{figure}

Figure~\ref{fig:plate_converge} shows the plots of normalised global $L_2$ and $H^1$ errors against the element size. The normalised global $L_2$ error is defined by
\begin{equation}
e_{L_2} = \frac{ \norm{u-u_h}_{L_2}}{ \norm{u}_{L_2}},
\end{equation}
where $\norm{\bullet}_{L_2}$ is the $L_2$ norm defined as $\norm{\bullet}_{L_2} = \sqrt{\int_\Gamma |\bullet|^2 \mathrm{d} \Gamma}$. The normalised global $H^1$ error is computed as
\begin{equation}
e_{H^1} = \frac{ \norm{u-u_h}_{H^1}}{ \norm{u}_{H^1}},
\end{equation}
where $\norm{\bullet}_{H^1}$ is the $H^1$ norm defined as $\norm{\bullet}_{H^1} =\sqrt{\int_\Gamma |\bullet|^2 \mathrm{d} \Gamma +  \int_\Gamma | \nabla(\bullet)|^2 \mathrm{d} \Gamma} $. We set the element size of the coarsest mesh as $1$. Then, the normalised element size for the refined meshes are $\frac{1}{2},  \frac{1}{4}, \cdots$. The convergence rate of Test 2 and 3 are sub-optimal at $2.5$ ($L_2$ error) and $1.5$ ($H^1$ error). The optimal convergence rate for cubic elements should be $p+1 = 4$ ($L_2$ error) and $p =3$ ($H^1$ error), where $p$ is the polynomial degree of the basis functions.  The numerical result captures the analytical solution well and the convergence rate for Test 4 is optimal. The same convergence study is now repeated starting from a mesh containing extraordinary vertices as shown in Figure~\ref{fig:plate_irreg_grid}. Figures~\ref{fig:plate_converge_l2} and~\ref{fig:plate_converge_h1} show the plots of normalised element sizes against the $L_2$ and $H^1$ errors, respectively. The same convergence rates are obtained for Tests 2 and 3. However, the convergence rate of Test 4 is also reduced to  $2.5$ ($L_2$ error) and $1.5$ ($H^1$ error). Figure~\ref{fig:plate_converge_l2_irreg} and~\ref{fig:plate_converge_h1_irreg} show the plots of normalised element sizes against $L_2$ and $H^1$ errors, respectively, for the mesh with an extraordinary vertex.

The Catmull-Clark subdivision method can pass the patch test when the function gradient is a constant but has difficulties to capture the gradients in boundary regions when they do not behave like a constant.  When the gradient behaves like a constant in the boundary regions, the optimal convergence rate can be obtained. If this is not the case, a reduction of the convergence rate is observed. The presence of the extraordinary vertex in the patch also reduces the convergence rate. It is also important to note that the Catmull-Clark subdivision elements have advantages in describing non-polynomial functions since their basis functions are cubic and $C^2$ continuous.

\begin{figure}
\centering
\begin{subfigure}{0.49\linewidth}
\centering
	\includegraphics[width=\linewidth]{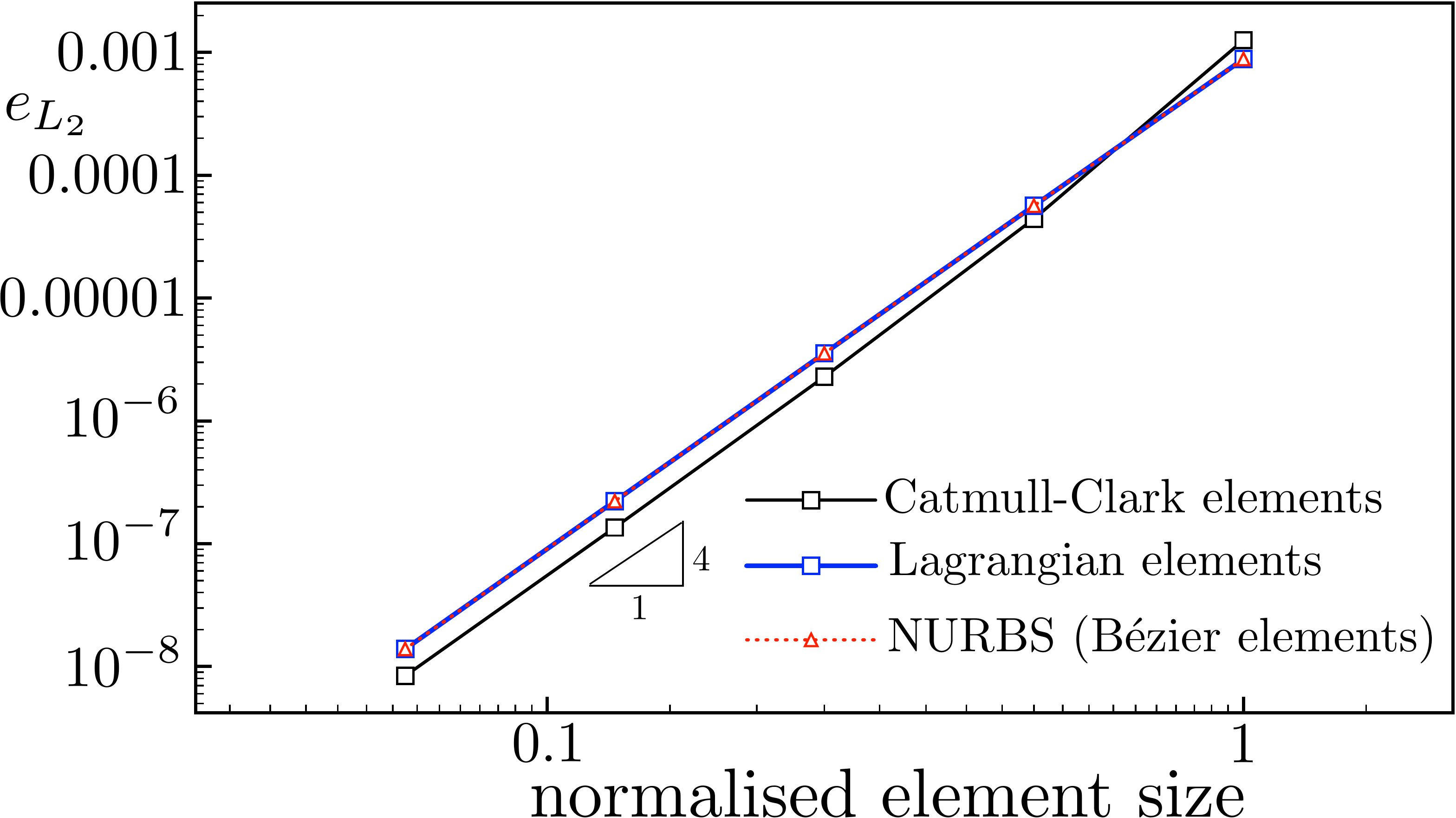}
	\caption{Normalised global $L^2$ error.}
	\label{fig:lagrangevsbersteinvsCC_l2}
\end{subfigure}
\begin{subfigure}{0.49\linewidth}
\centering
	\includegraphics[width=\linewidth]{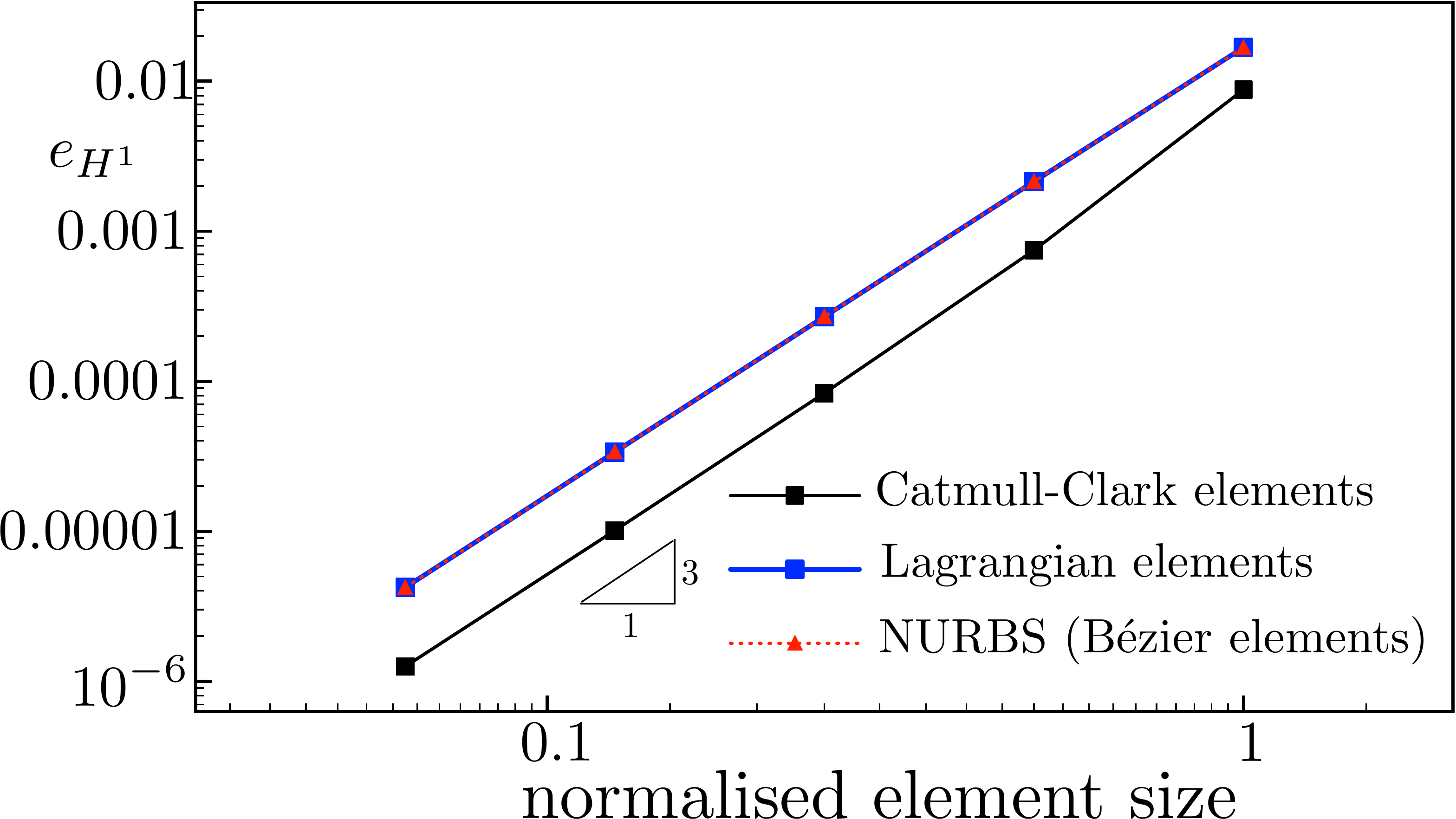}
	\caption{Normalised global $H^1$ error.}
	\label{fig:lagrangevsbersteinvsCC_h1}
\end{subfigure}
\caption{Convergence studies for Test 4 using Mesh 1. The Catmull-Clark elements are compared with Lagrangian elements and NURBS (B\'ezier elements). $p = 3$ for all cases.}
\label{fig:plate_converge_compare_elements}
\end{figure}
{
\paragraph*{Comparison with NURBS and Lagrangian elements \\}
We now compare the convergence rate associated with Catmull-Clark elements against conventional Lagrangian elements and NURBS.  B\'ezier extraction~\cite{borden2011isogeometric} is adopted to decompose a NURBS surfaces into $C_0$ B\'ezier elements to provide an element structure for the isogeometric Galerkin method. This is a widely-used method for isogeometric analysis using T-splines~\cite{scott2011isogeometric}. As the Lagrangian and B\'ezier elements can fully pass the `patch test', they both have no approximation error for Test 1, 2 and 3. Figure~\ref{fig:plate_converge_compare_elements} compares their behaviour in approximating non-polynomial solution in Test 4. Mesh 1 is used for all methods. All methods exhibit an optimal convergence rate. Since no geometry error is involved in the `patch test’, the B\'ezier element provides the same performance as the Lagrangian element without the advantages of exact geometry representation. The Catmull-Clark element is slightly more accurate than other two methods for this specific test.  }

\subsection{The Laplace-Beltrami equation}
The following sections will solve the Laplace-Beltrami equation~\eqref{eq:Poisson} on different two dimensional Catmull-Clark subdivision manifolds. An analytical solution of function $u^m$ is manufactured as
\begin{equation}
u^m (\mathbf{x}) =  \sin(\pi x_1) \cos(\pi x_2)e^{x_3},
\label{eq:solution}
\end{equation}
where $\mathbf{x}(x_1,x_2,x_3)$ is a point on the two dimensional manifold in three dimensional space. Applying the Laplacian operator on $u^m$ gives
\begin{equation}
\Delta u^u (\mathbf{x}) = -2\pi^2 \sin(\pi x_1) \cos(\pi x_2)e^{x_3} + \sin(\pi x_1) \cos(\pi x_2)e^{x_3}.
\label{eq:analytical_term}
\end{equation}
Then, the Hessian matrix $\nabla^2 u(\mathbf{x}) $  is computed as
\begin{equation}
\nabla^2 u^m(\mathbf{x}) =
\left[ 
\begin{array}{ccc}
-\pi^2 \sin(\pi x_1) \cos(\pi x_2)e^{x_3} & -\pi^2 \cos (\pi x_1) \sin(\pi x_2)e^{x_3} &  \pi \cos(\pi x_1)\cos(\pi x_2)e^{x_3}\\
-\pi^2\cos (\pi x_1) \sin(\pi x_2)e^{x_3} & -\pi^2 \sin(\pi x_1) \cos(\pi x_2)e^{x_3} &  -\pi \sin(\pi x_1) \sin(\pi x_2)e^{x_3}\\
  \pi \cos(\pi x_1)\cos(\pi x_2)e^{x_3}& -\pi \sin(\pi x_1) \sin(\pi x_2)e^{x_3}   & \sin(\pi x_1) \cos(\pi x_2)e^{x_3}
\end{array}
\right]
\end{equation}

The second term in~\eqref{eq:bt_rhs} includes the normal vector $\mathbf{n}(\mathbf{x})$ and its gradient which can not be computed analytically. In the present work, an $L^2$ projection is used to compute the coefficients of normal vectors associated with control points, denoted by $\hat{\mathbf{n}}$, in order to numerically interpolate the normal vector derivatives at any surface points, thus
\begin{equation}
\frac{\partial \mathbf{n}}{\partial \mathbf{x}} = \sum_{A = 1}^{n_b} \frac{\partial N_A}{\partial \mathbf{x}} \hat {\mathbf{n}}_A.
\end{equation}

\subsubsection{Cylindrical surface example}
\begin{figure}
\centering
\begin{subfigure}{0.49\linewidth}
\centering
  \includegraphics[width=0.9\linewidth]{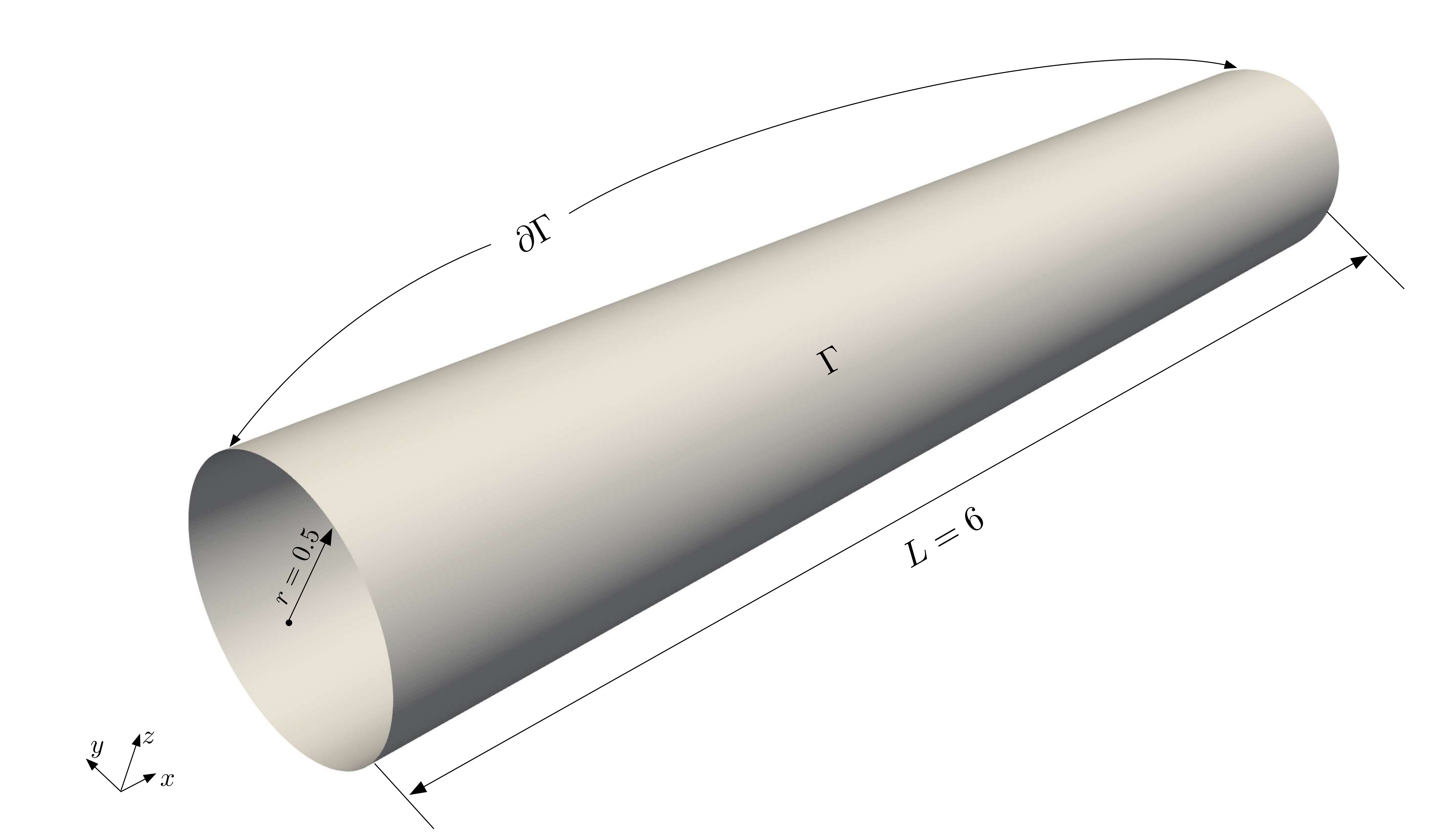}
\caption{Geometry}
\label{fig:cylinder_limit}
\end{subfigure}
\begin{subfigure}{0.49\linewidth}
\centering
  \includegraphics[width=0.9\linewidth]{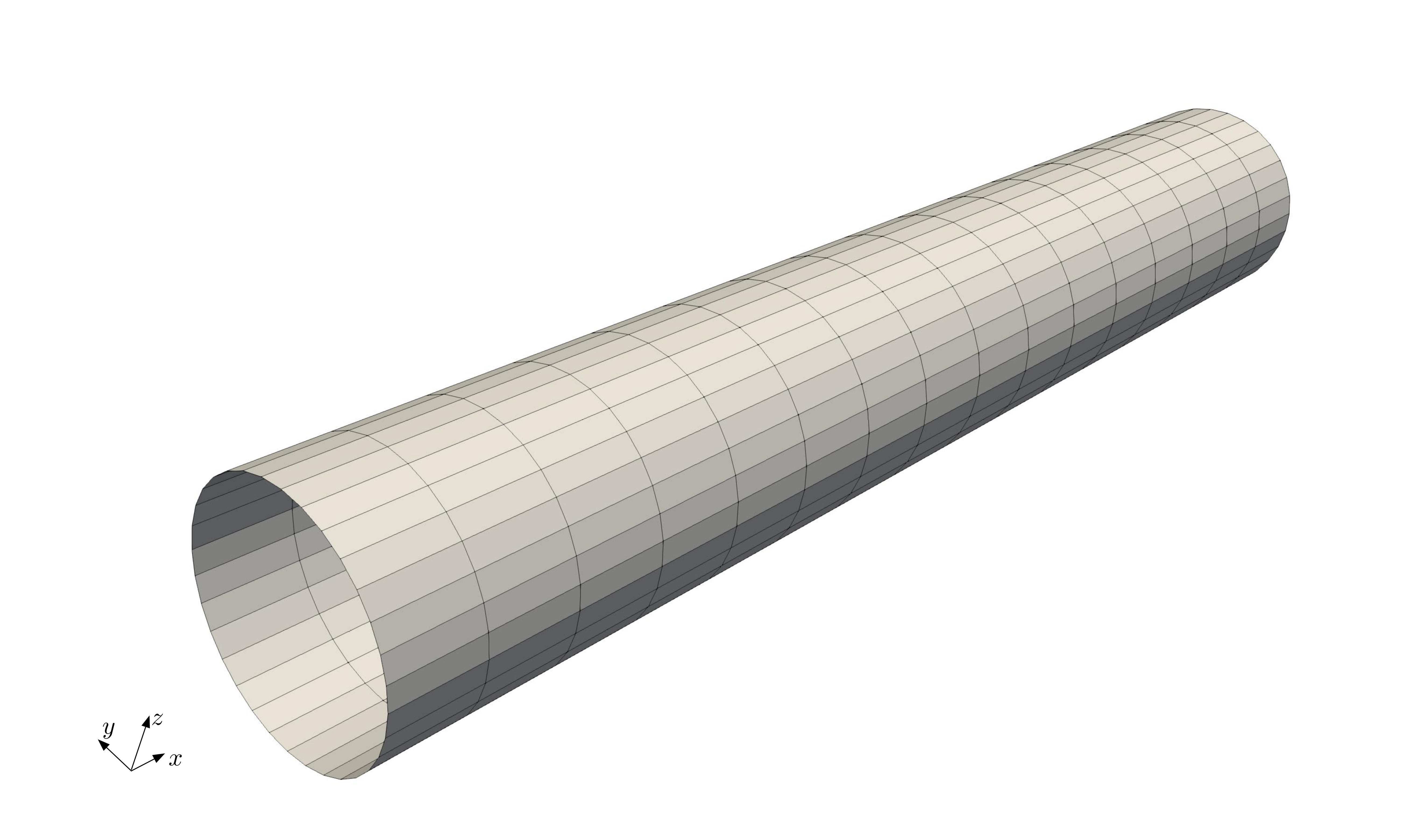}
\caption{Control mesh}
\label{fig:cylinder_init}
\end{subfigure}
\begin{subfigure}{0.7\linewidth}
\centering
\includegraphics[width=0.9	\linewidth]{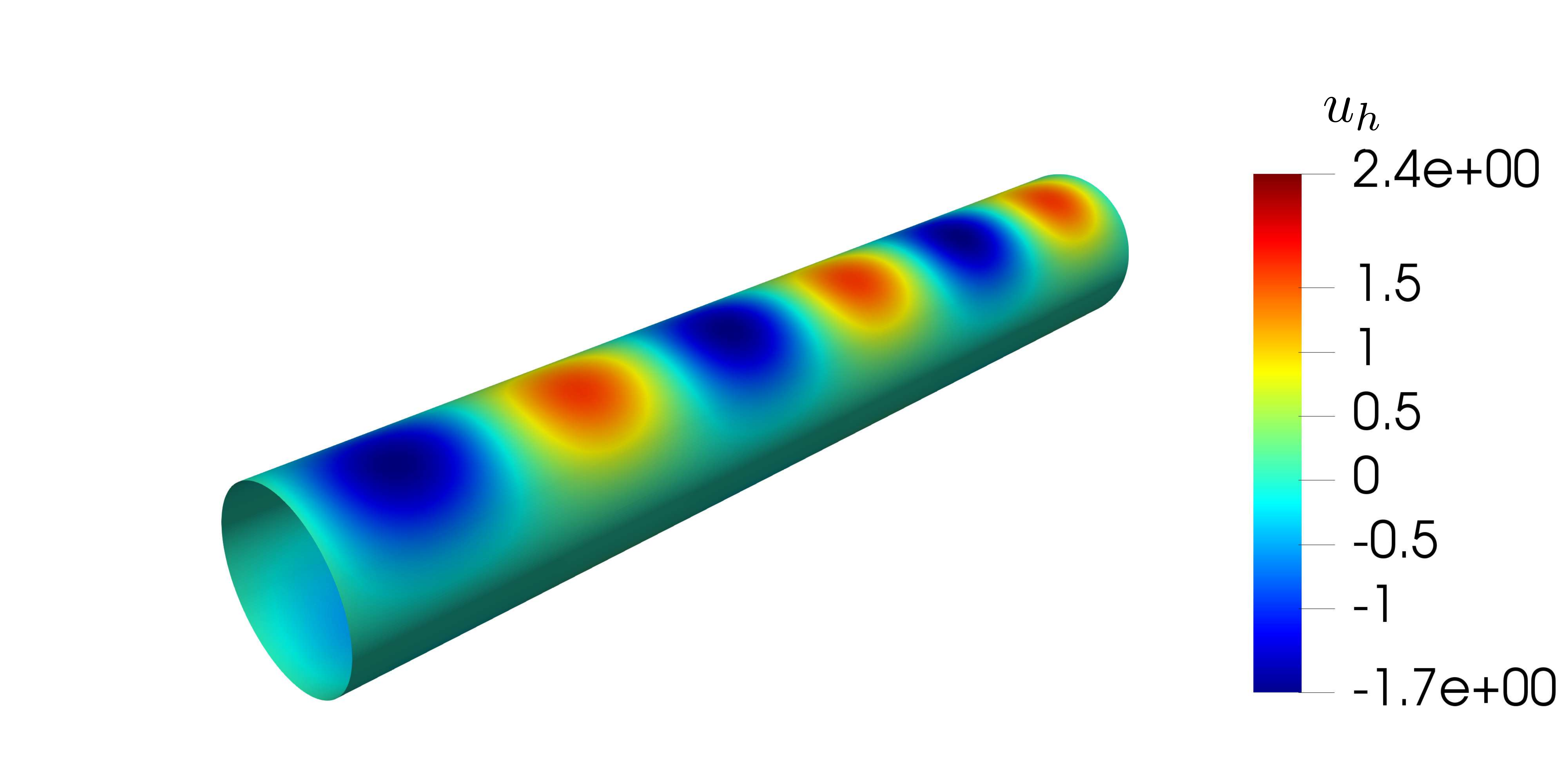}
\caption{Numerical result.}
\label{fig:cylinder_result}
\end{subfigure}
\caption{The geometry is given in Figure (a) and Figure (b) is the control mesh which constructs the best approximating Catmull-Clark subdivision surface of the given geometry. The control mesh is generated using least-squares fitting. Figure (c) shows the numerical result $u_h$ on the cylindrical surface.}
\end{figure}

The first numerical example considered is a cylindrical surface. The analysis domain of the problem is the cylindrical surface shown in Figure~\ref{fig:cylinder_limit}. Surfaces fitting methods are used to construct the control mesh, see Section~\ref{sec:geometry_fitting}. The first level control mesh is shown in Figure~\ref{fig:cylinder_init}. This has no extraordinary vertices. The Laplace-Beltrami problem on this manifold domain is solved using the Galerkin formulation presented in Section~\ref{sec:FE}. {Essential boundary conditions are applied on $\partial \Gamma$.} The right-hand side function $f$ is computed using the definition in Equation~\eqref{eq:bt_rhs}. Figure~\ref{fig:cylinder_result} shows the numerical result $u_h$ which matches the manufactured analytical solution~\eqref{eq:solution} very well.  

\begin{figure}
\centering
\begin{subfigure}{0.49\linewidth}
\centering
\includegraphics[width=0.9	\linewidth]{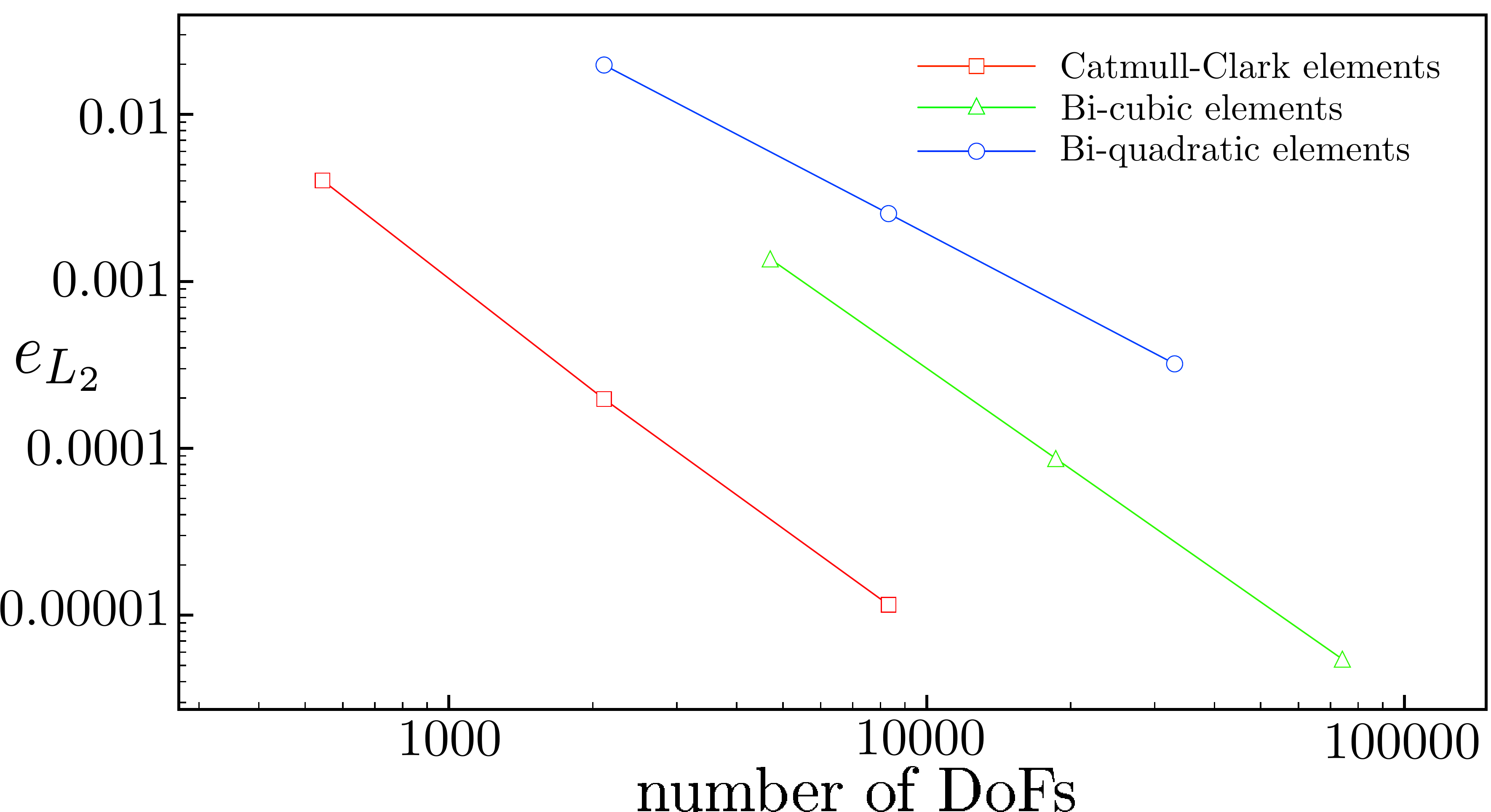}
\caption{Global $L_2$ error against number of degrees of freedom.}
\label{fig:cylinder_convergence_dofs}
\end{subfigure}
\begin{subfigure}{0.49\linewidth}
\centering
\includegraphics[width=0.9	\linewidth]{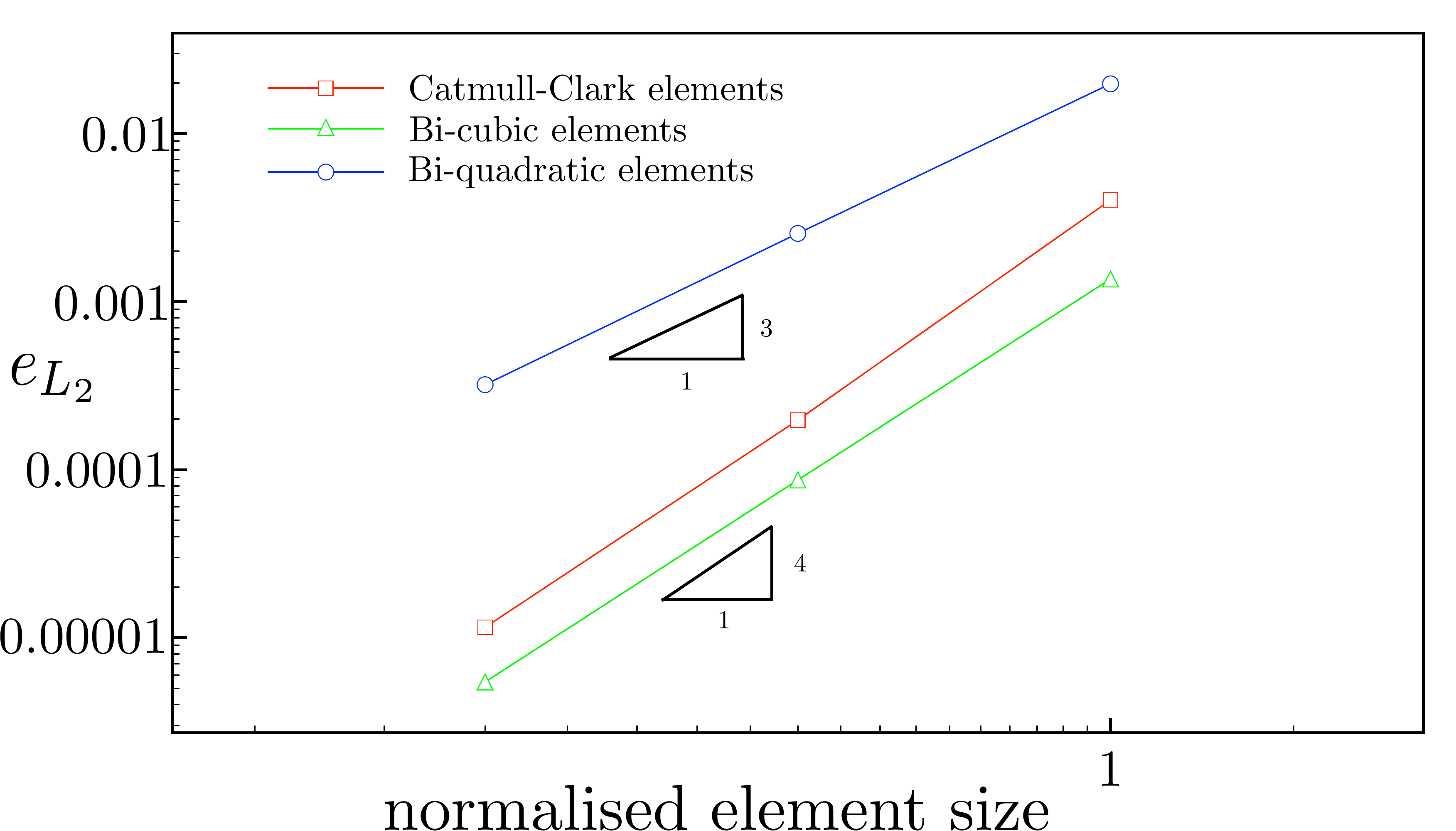}
\caption{Global $L_2$ error against normalised element size.}
\label{fig:cylinder_convergence_elements}
\end{subfigure}
\caption{Convergence study for cylindrical example: comparison of the Catmull-Clark elements to the quadratic and cubic Lagrangian elements.}
\label{fig:cylinder_convergence}
\end{figure}

A convergence study is now conducted for this geometry. The refined control meshes are constructed using the least-squares fitting method described in Section~\ref{sec:geometry_fitting}. Figure~\ref{fig:cylinder_convergence} compares the convergence rates between Catmull-Clark subdivision surfaces with two different order Lagrangian elements. In this example, the shortcoming caused by extraordinary vertices and boundary gradients are not present, and the Catmull-Clark subdivision surfaces have the same convergence rate $p+1$ as cubic Lagrangian elements. 

\subsubsection{Hemispherical surface example}
\begin{figure}
\centering
\begin{subfigure}{0.49\linewidth}
\centering
  \includegraphics[width=0.9\linewidth]{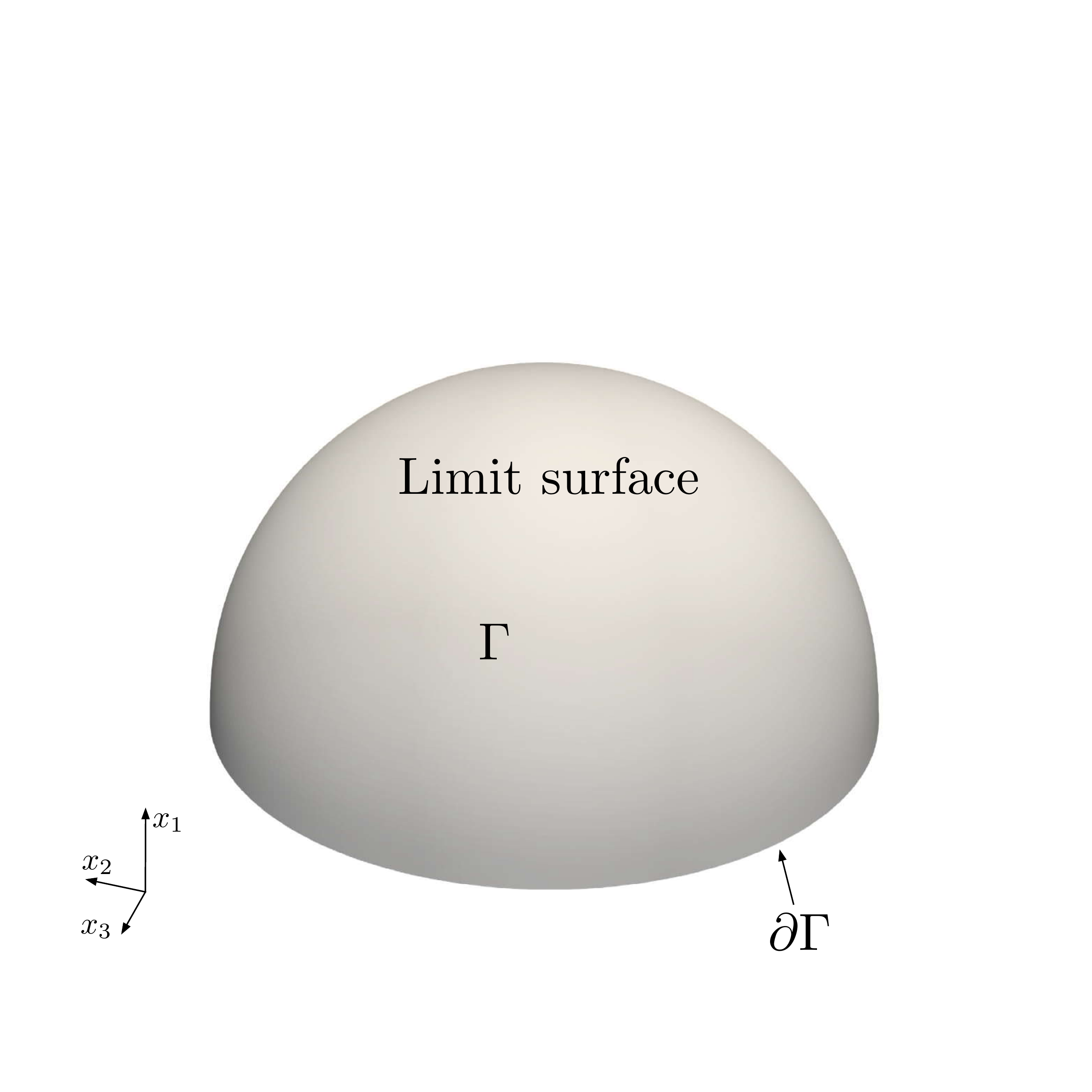}
\caption{}
\label{fig:sphere_geometry}
\end{subfigure}
\begin{subfigure}{0.49\linewidth}
\centering
  \includegraphics[width=0.81\linewidth]{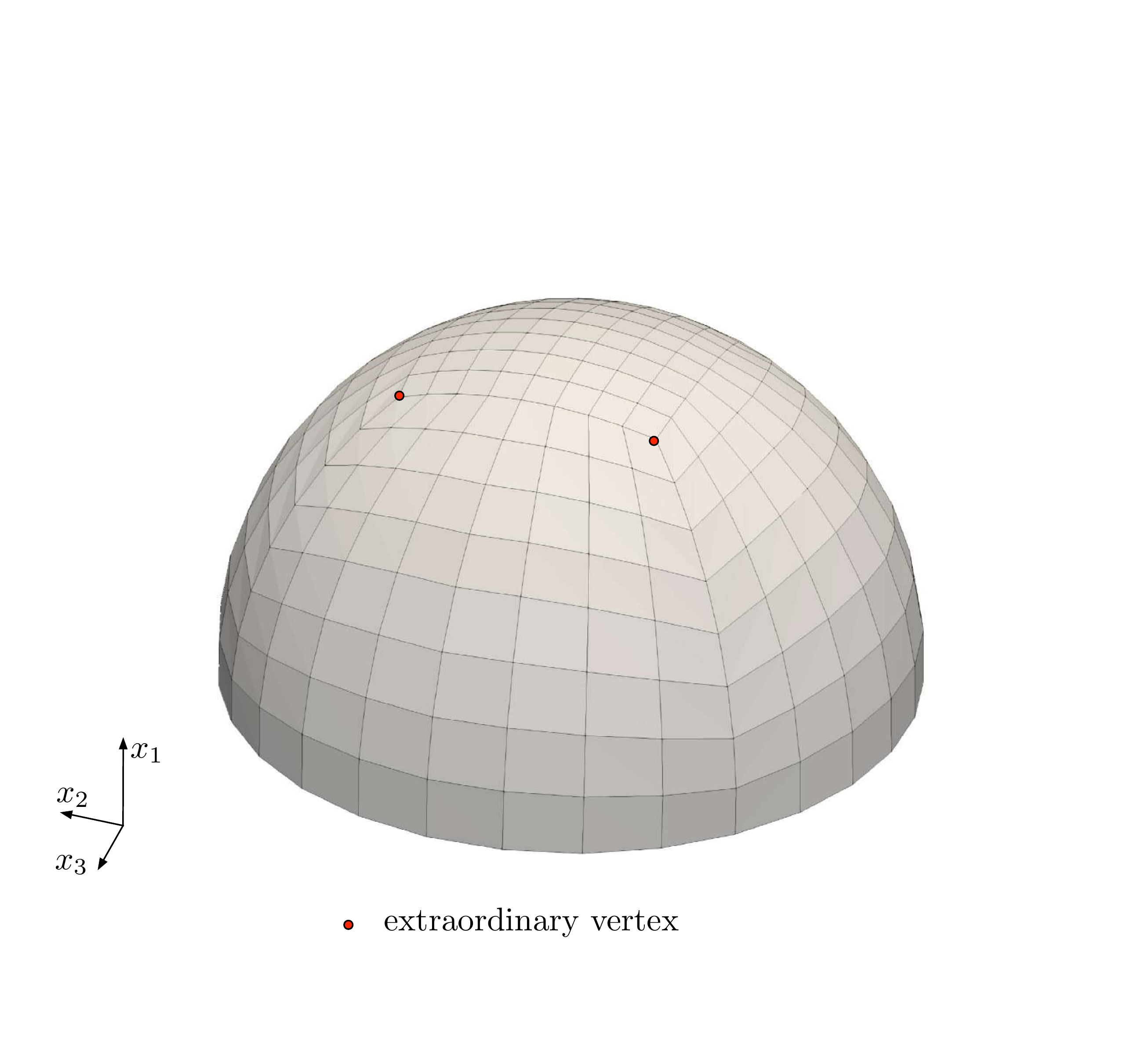}
\caption{}
\label{fig:sphere_grid_0}
\end{subfigure}
\begin{subfigure}{0.49\linewidth}
\centering
  \includegraphics[width=0.72\linewidth]{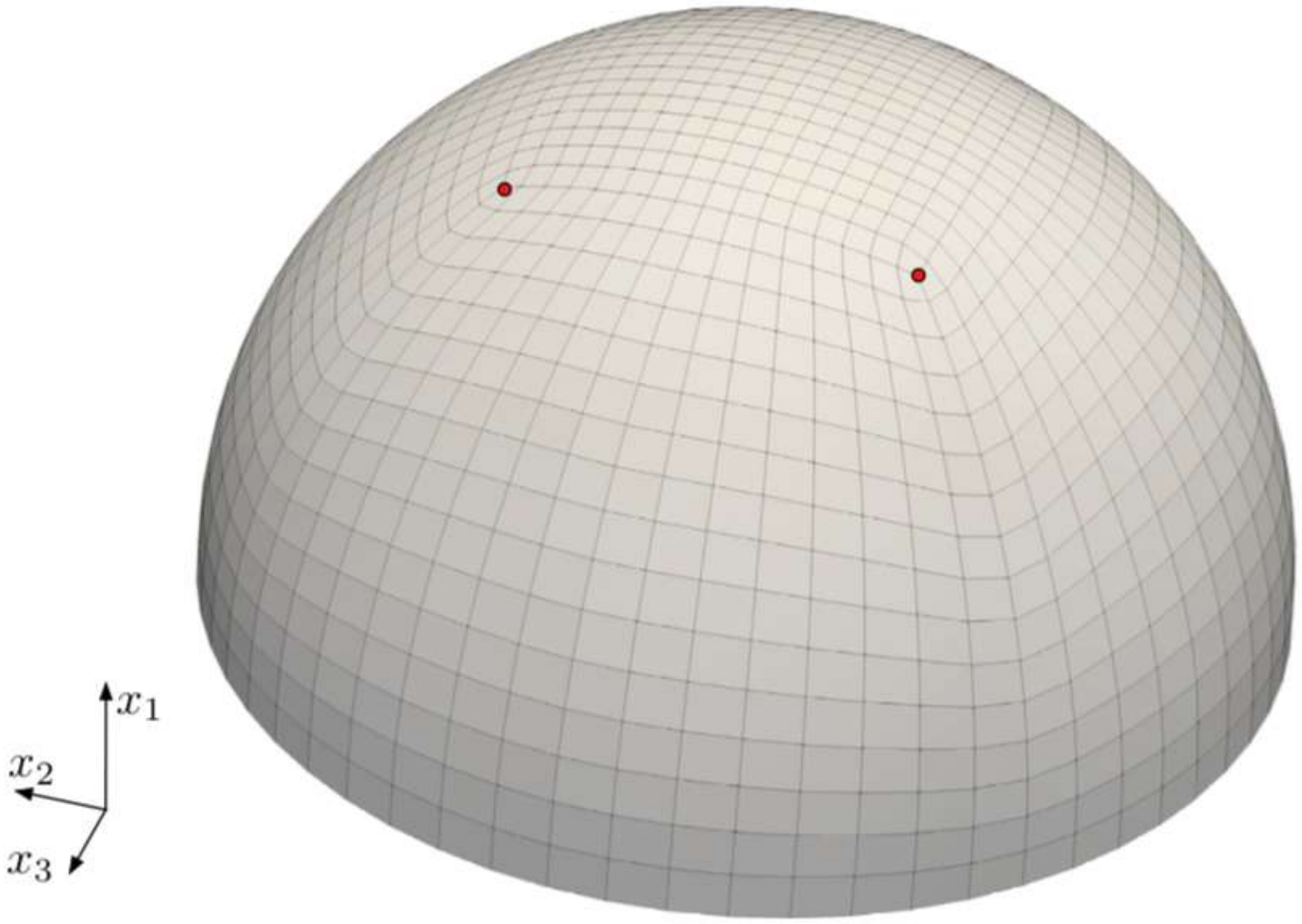}
\caption{}
\label{fig:sphere_grid_1}
\end{subfigure}
\begin{subfigure}{0.49\linewidth}
\centering
  \includegraphics[width=0.72\linewidth]{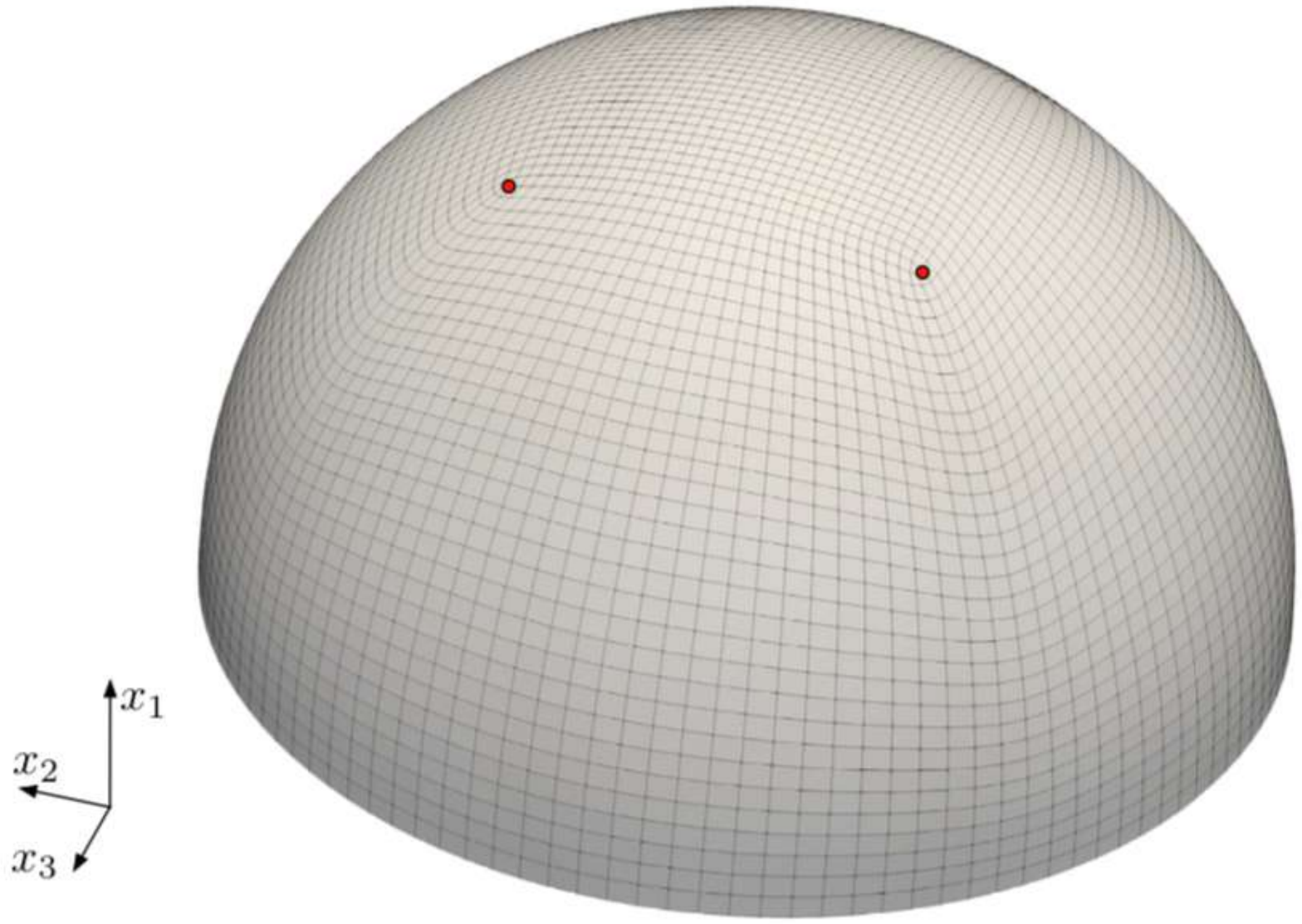}
\caption{}
\label{fig:sphere_grid_2}
\end{subfigure}
\begin{subfigure}{0.7\linewidth}
\centering
\includegraphics[width=0.7\linewidth]{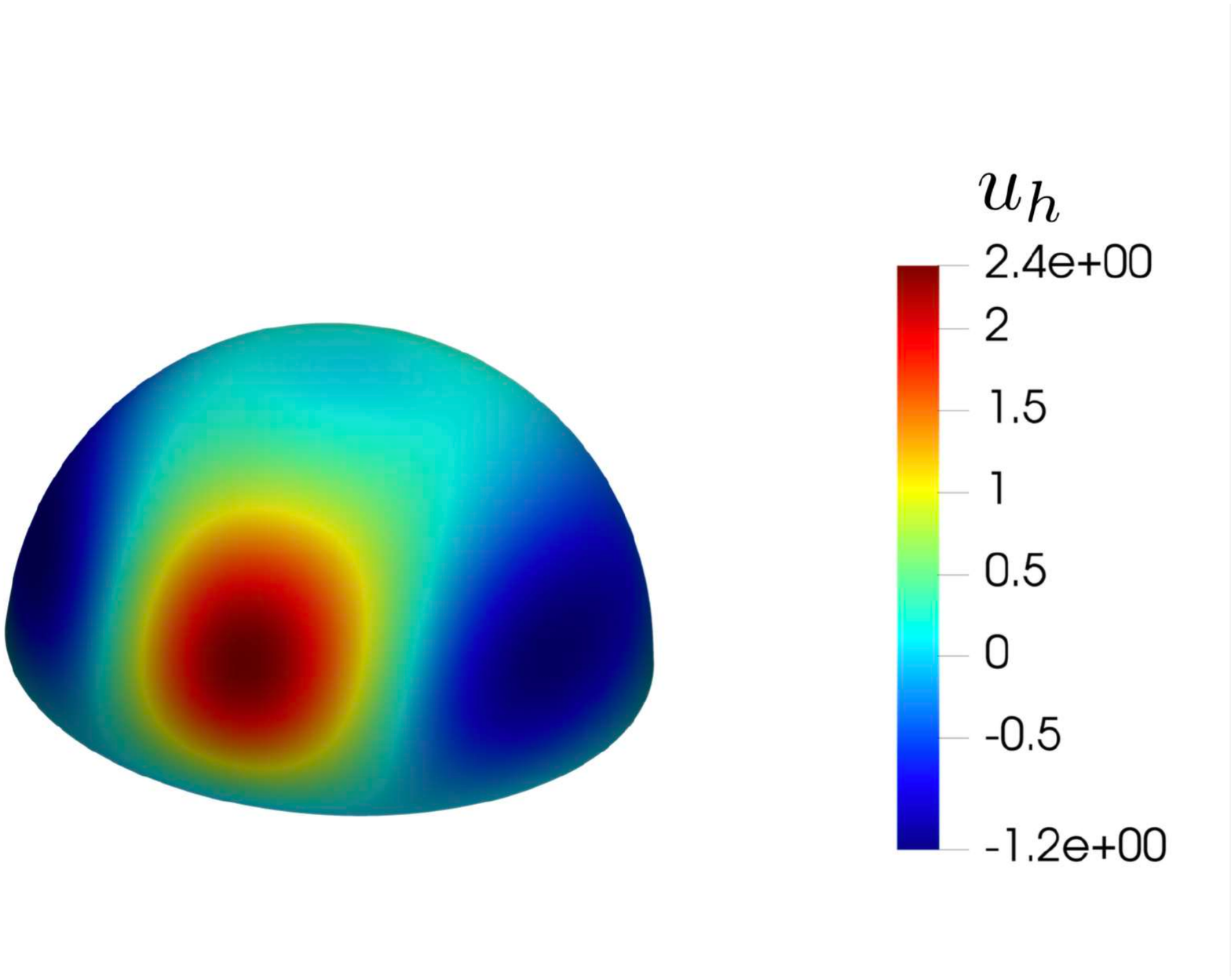}
\caption{}
\label{fig:sphere_result}
\end{subfigure}
\label{fig:sphere_refinements}
\caption{(a) is a hemispherical surface. (b) is the control mesh for constructing subdivision surfaces to fit the hemispherical surface.(c) is 1-Level refined mesh for the hemispherical surface. (d) is 2-level refined mesh for the hemispherical surface. (e) shows the numerical result $u_h$ on this surface.}
\end{figure}

The second geometry investigated is a hemispherical surface with radius equal to 1 as shown in Figure~\ref{fig:sphere_geometry}. We use the same strategy to fit the Catmull-Clark subdivision surfaces to the hemispherical surface. The control mesh shown in Figure~\ref{fig:sphere_grid_0} is generated to discretise the surface into a number of Catmull-Clark elements. The control mesh has four extraordinary vertices. Figure~\ref{fig:sphere_result} shows the solution $u_h$.
\paragraph*{Convergence study with an isogeometric approach \\}
In engineering, designers usually do not know the geometry of the product in advance. The geometry information is purely from the CAD model. Catmull-Clark subdivision surfaces, as a design tool, provide the geometry which is the design of the engineering product. In this case, engineers do not need to approximate the given geometry with Catmull-Clark elements. They can directly adopt the discretisation from the CAD model for analysis. For example, we adopt the control mesh shown in Figure~\ref{fig:sphere_grid_0} as the initial control mesh. It can be used to generate a limit surface approximating a hemisphere, as shown in Figure~\ref{fig:sphere_geometry}, with Catmull-Clark subdivision bases.  It is important to note the limit surface is not an exact hemisphere since it is evaluated using cubic basis spline functions. However, this surface is the domain of our problem and it will stay exact the same during the entire analysis (isogeometric) and $h$-refinement with subdivision algorithm will not change the geometry.  

The same problem is solved on the subdivision surfaces. A convergence study is done with another two levels of subdivision control mesh as shown in Figure~\ref{fig:sphere_grid_1} and~\ref{fig:sphere_grid_2}. Note, refinement does not change the number of extraordinary vertices. The two new meshes still have four extraordinary vertices. The two control meshes can be used to evaluate the same limit surface shown in Figure~\ref{fig:sphere_geometry}. 
The Catmull-Clark subdivision surfaces are compared with quadratic and cubic Lagrangian elements. Generally, Catmull-Clark subdivision elements can achieve higher accuracy per degree of freedom than Lagrangian elements.  From the initial to the second level of mesh refinement, the Catmull-Clark subdivision elements have a similar convergence rate to cubic Lagrangian elements. After that, the convergence rate is equivalent to quadratic Lagrangian elements. 
\begin{figure}
\centering

\begin{subfigure}{.49\linewidth}
\centering
\includegraphics[width=0.9	\linewidth]{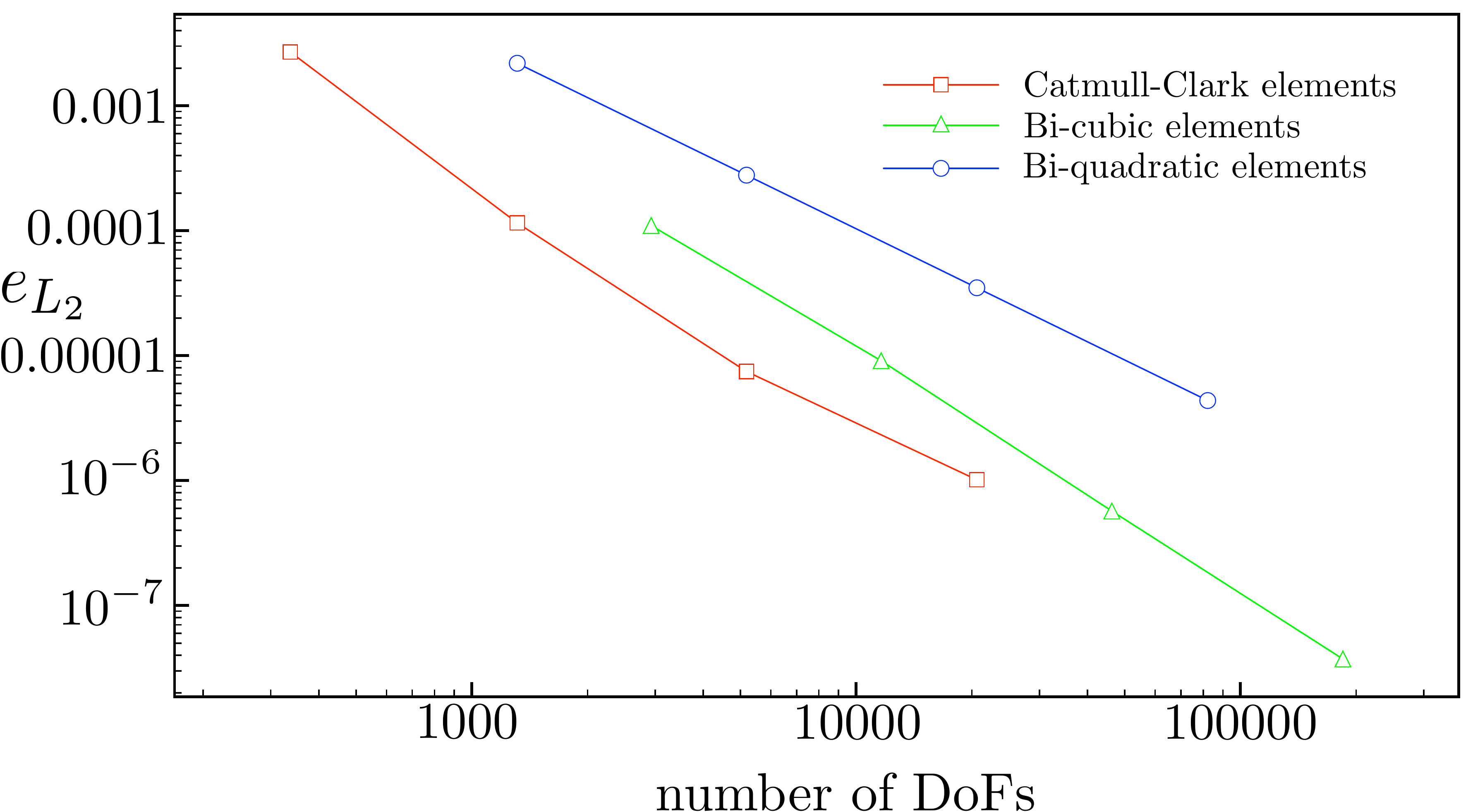}
\caption{$e_{L_2}$ against number of degrees of freedom.}
\end{subfigure}
\begin{subfigure}{.49\linewidth}
\centering
\includegraphics[width=0.9\linewidth]{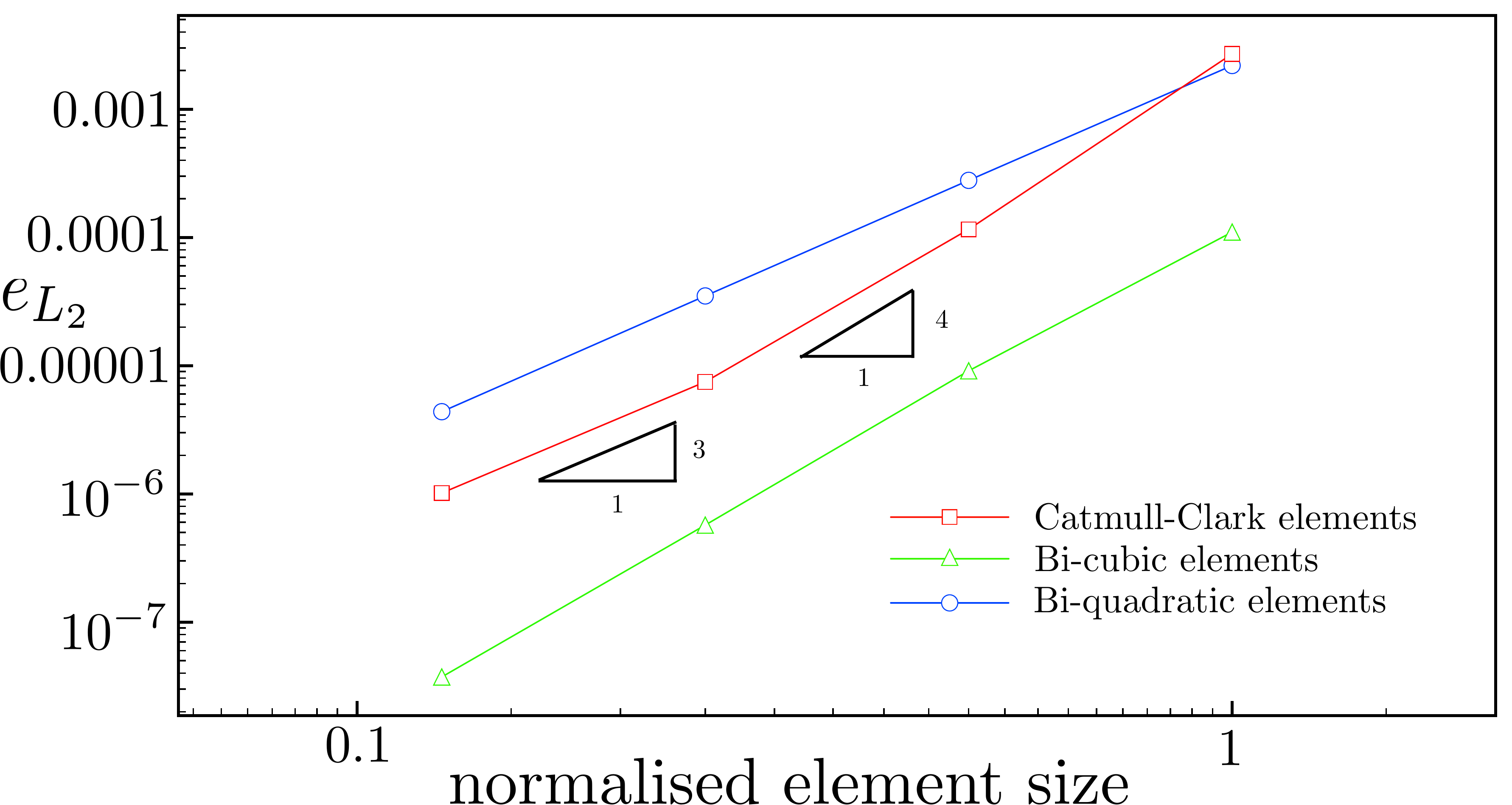}
\caption{$e_{L_2}$ against normalised element sizes.}

\end{subfigure}

\label{fig:convergence}
\caption{Convergence study: comparison of the Catmull-Clark elements with the quadratic and cubic Lagrangian elements.}
\end{figure}

\paragraph*{Sparsity patterns\\}

\begin{figure}
\centering
\begin{subfigure}{0.32\linewidth}
\includegraphics[width=\linewidth]{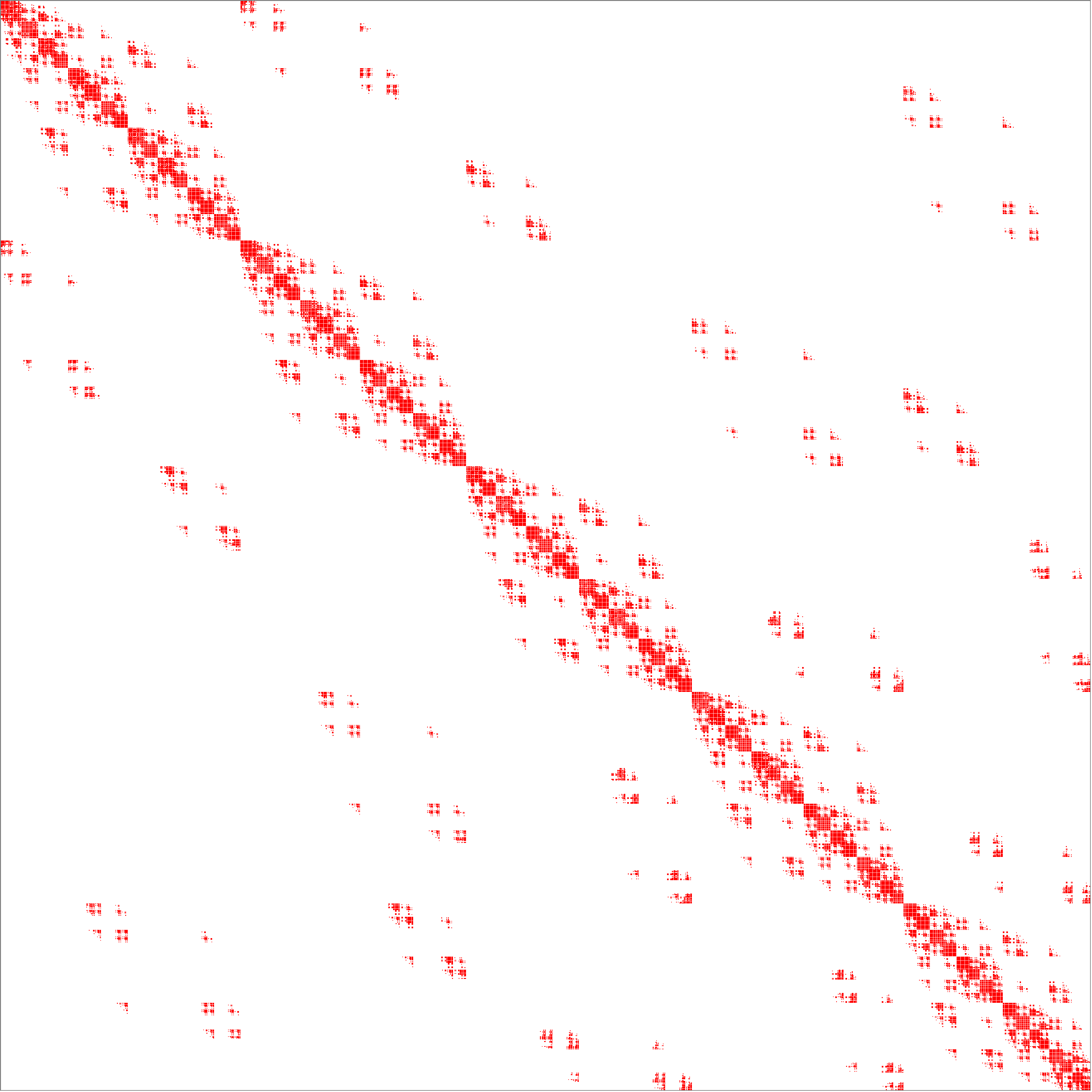}
\caption{Catmull-Clark subdivision elements (1280 elements, 1313 DoFs).}
\label{fig:dsp_sparsity_pattern_cc}
\end{subfigure}
\begin{subfigure}{0.32\linewidth}
\includegraphics[width=\linewidth]{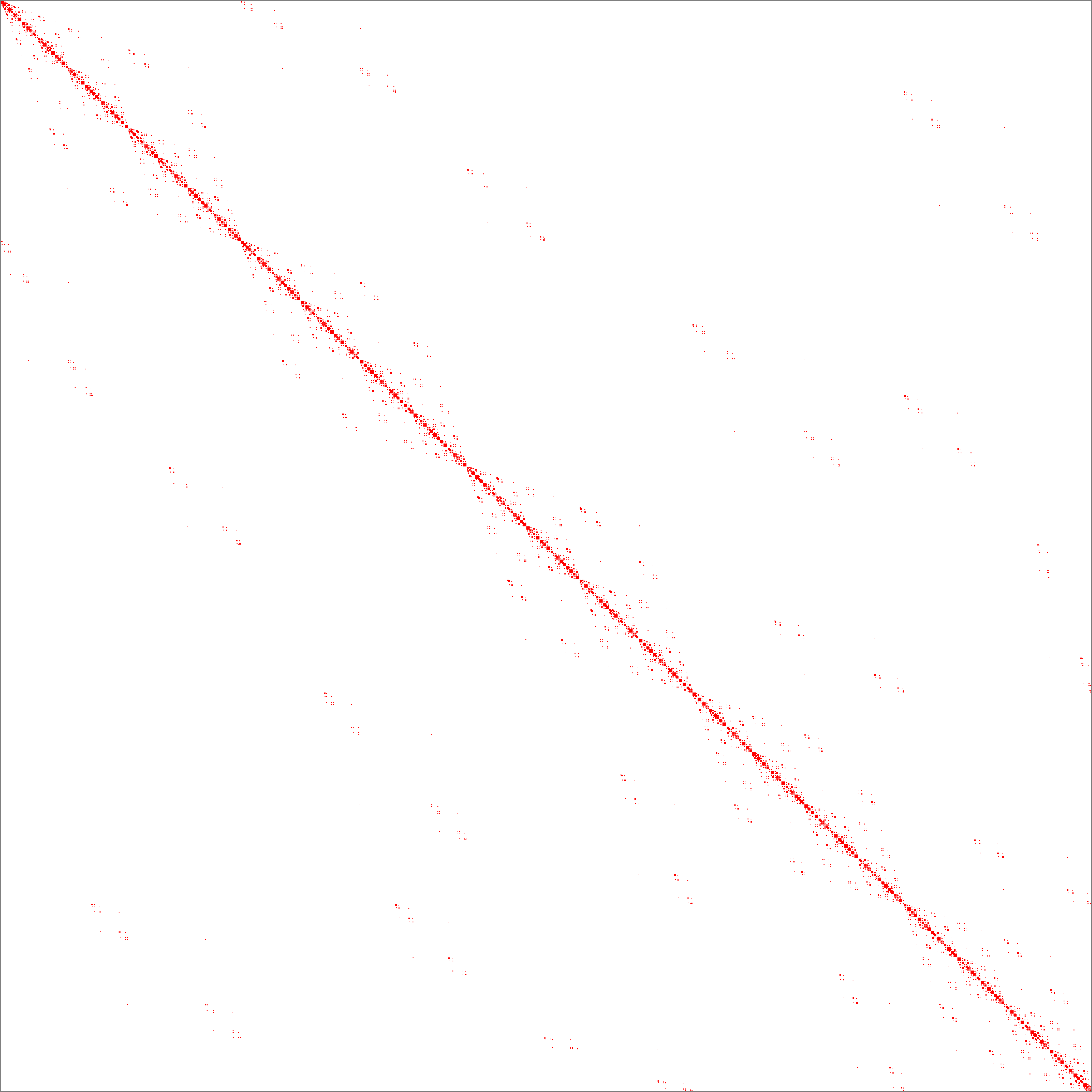}
\caption{Linear Lagrangian elements (1280 elements, 1313 DoFs).}
\label{fig:dsp_sparsity_pattern_linear}
\end{subfigure}
\begin{subfigure}{0.32\linewidth}
\includegraphics[width=\linewidth]{dsp_sparsity_pattern_cubic.pdf}
\caption{Cubic Lagrangian elements (1280 elements, 11617 DoFs).}
\label{fig:dsp_sparsity_pattern_cubic}
\end{subfigure}
\caption{Comparison of sparsity patterns between the Catmull-Clark elements and the Lagrangian elements}
\end{figure}
Figure~\ref{fig:dsp_sparsity_pattern_cc} shows the sparsity pattern of the system matrix $\ary{K}$ for the Catmull-Clark subdivision discretisation. The size of the matrix is the same as the system matrix assembled using a linear Lagrange discretisation. However, because the Catmull-Clark subdivision discretisation uses cubic basis functions with non-local support and there are 16 shape functions in a subdivision element with no extraordinary vertex, the number of non-zero entries in columns and rows is more than the linear Lagrange discretisation  (i.e. the sparsity is decreased and the bandwidth increased). Thus, the system matrix of a Catmull-Clark subdivision discretisation has the same size but is denser than the linear Lagrange discretisation shown in~\ref{fig:dsp_sparsity_pattern_linear}. Figure~\ref{fig:dsp_sparsity_pattern_cubic} is the sparsity patterns of cubic Lagrange discretisations. $p$-refinement increase the number of degrees of freedom as well as the number of non-zero entries in rows and columns. Thus there is no significant change in the density of the system matrices. The Catmull-Clark subdivision discretisation has the same number of non-zero entries in {each row or column} as the cubic Lagrangian discretisation but has a much smaller size.

\subsection{Investigation of extraordinary vertices}
\begin{figure}
\centering
\captionsetup{justification=centering}
\begin{subfigure}{.5\linewidth}
  \centering
  \includegraphics[width=\linewidth]{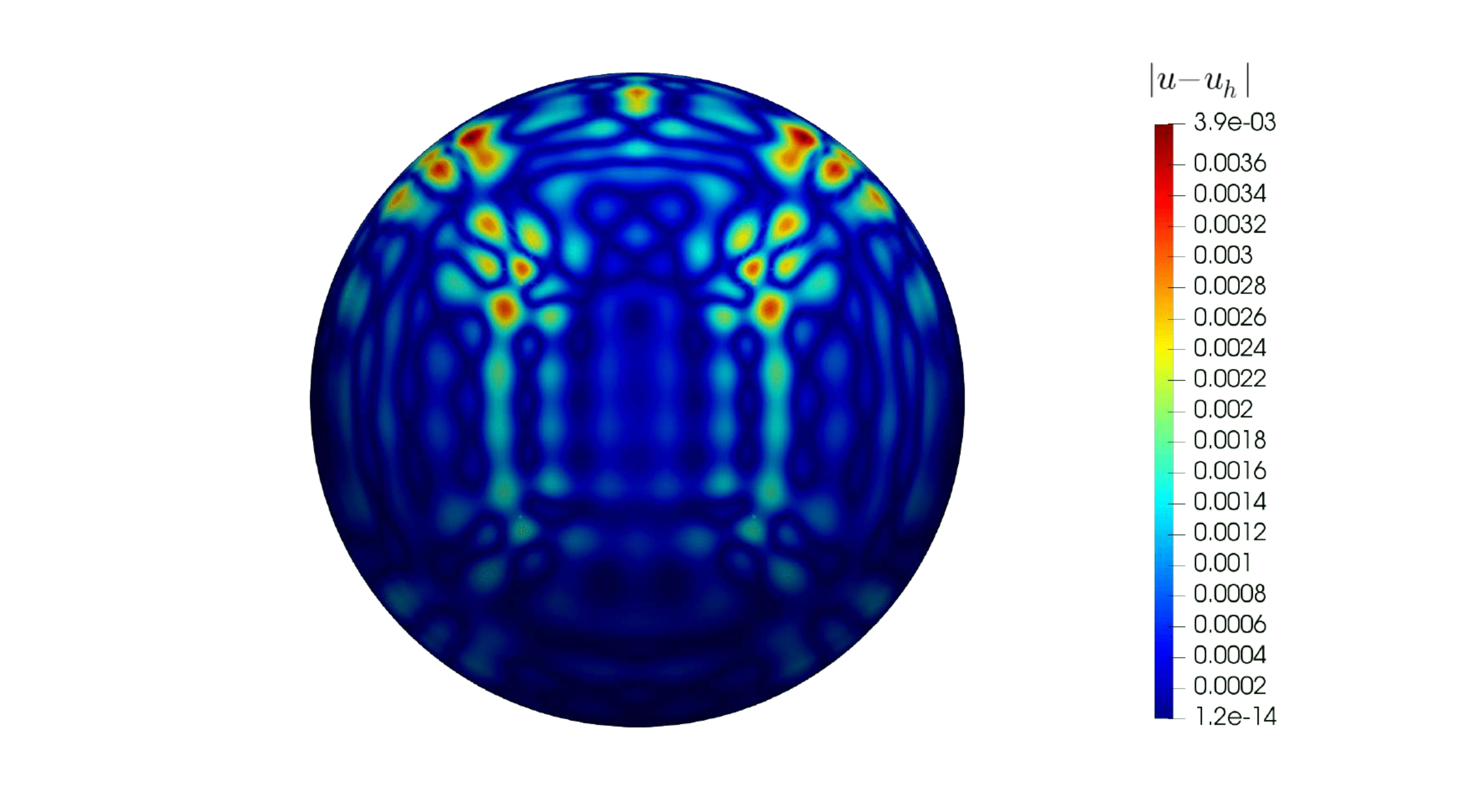}
  \caption{Initial mesh with standard Gauss quadrature.}
  \label{fig:sphere_error_plot_0}
\end{subfigure}%
\begin{subfigure}{.5\linewidth}
  \centering
  \includegraphics[width=\linewidth]{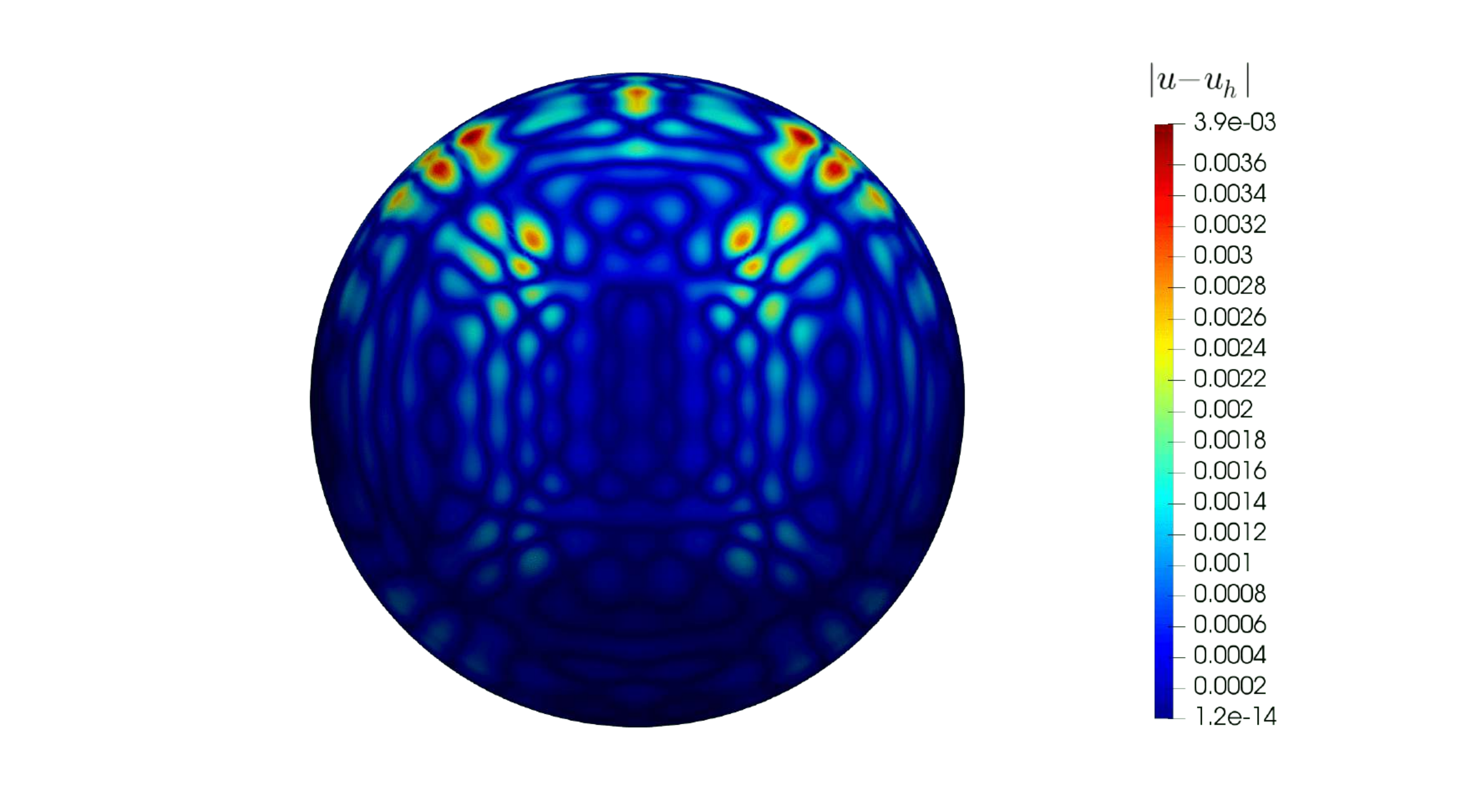}
  \caption{Initial mesh with adaptive Gauss quadrature.}
  \label{fig:sphere_error_plot_aq_0}
\end{subfigure}

\vspace{2ex}

\begin{subfigure}{.5\linewidth}
  \centering
  \includegraphics[width=\linewidth]{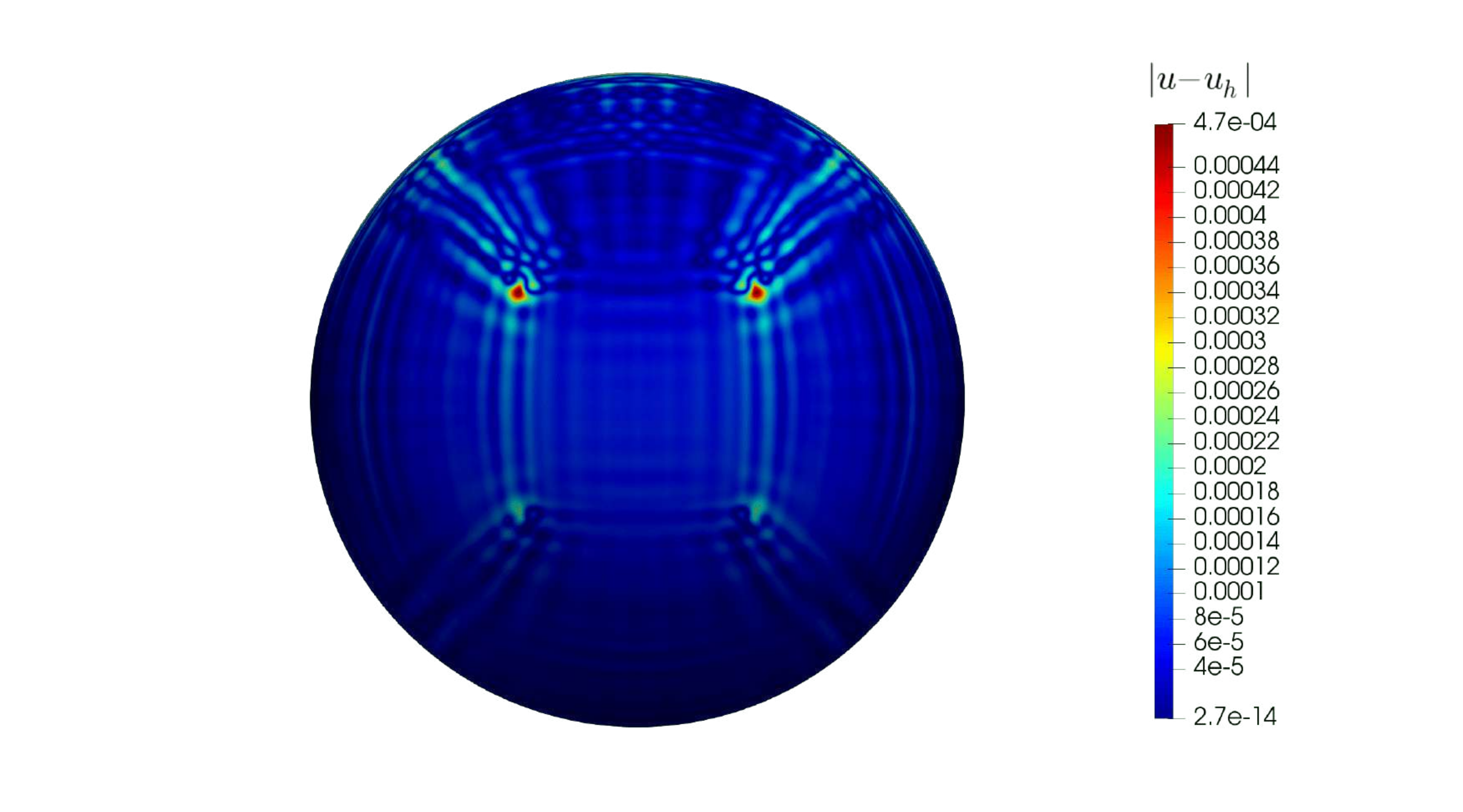}
  \caption{$1^{\text{st}}$ level subdivision mesh with standard Gauss quadrature.}
  \label{fig:sphere_error_plot_1}
\end{subfigure}%
\begin{subfigure}{.5\linewidth}
  \centering
  \includegraphics[width=\linewidth]{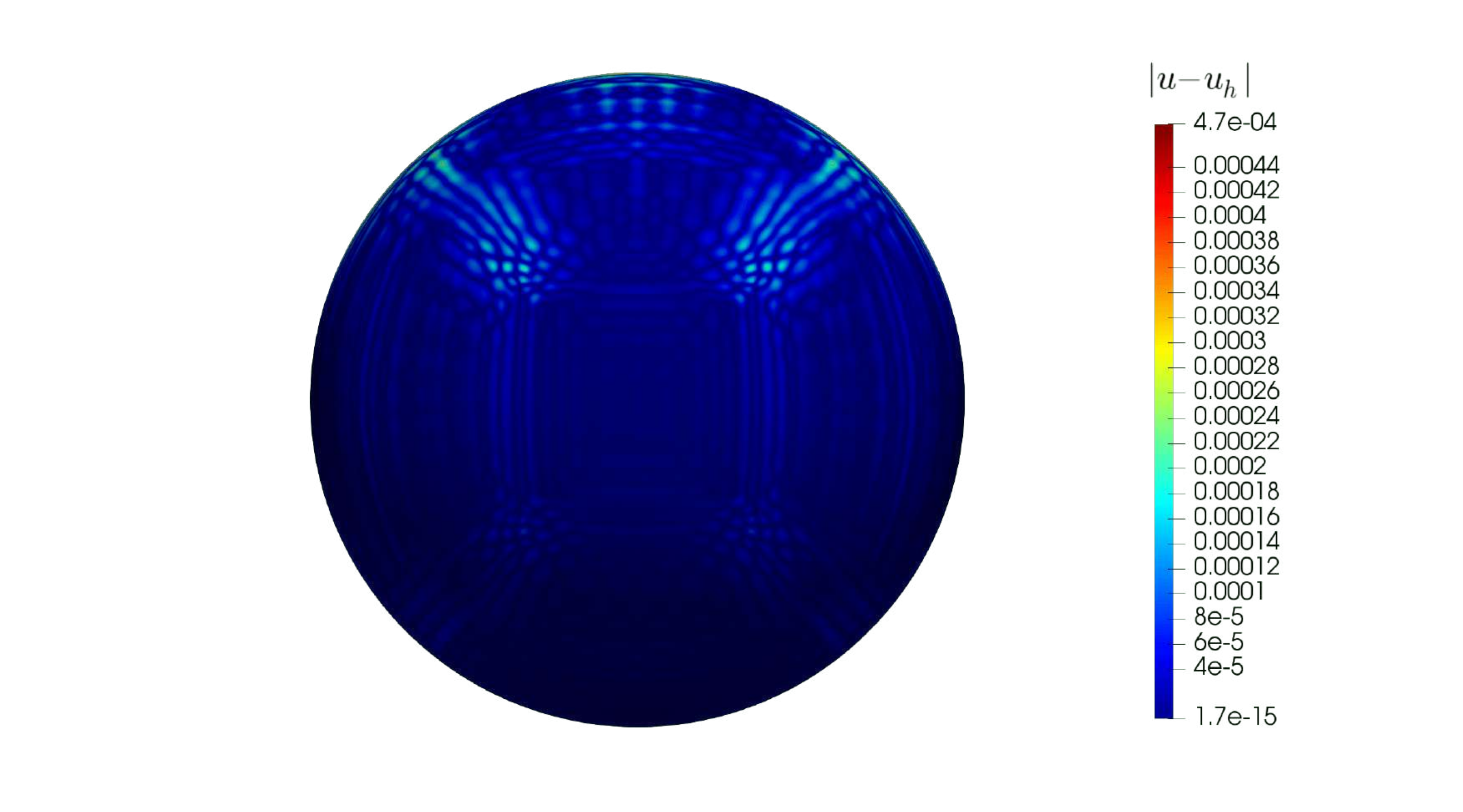}
  \caption{$1^{\text{st}}$ level subdivision mesh with adaptive Gauss quadrature.}
  \label{fig:sphere_error_plot_aq_1}
\end{subfigure}

\vspace{2ex}

\begin{subfigure}{.5\linewidth}
  \centering
  \includegraphics[width=\linewidth]{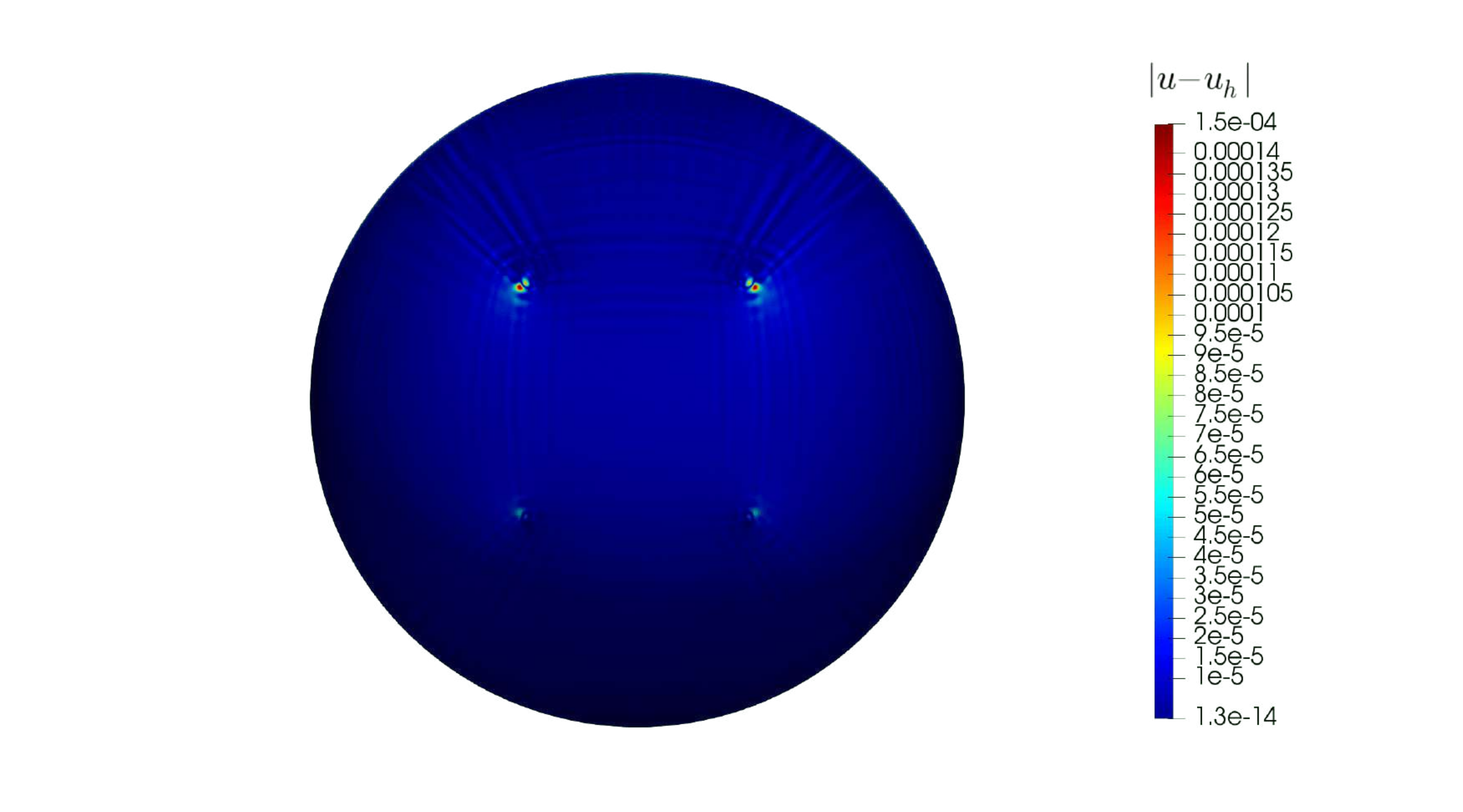}
  \caption{$2^{\text{nd}}$ level subdivision mesh with standard Gauss quadrature.}
  \label{fig:sphere_error_plot_2}
\end{subfigure}%
\begin{subfigure}{.5\linewidth}
  \centering
  \includegraphics[width=\linewidth]{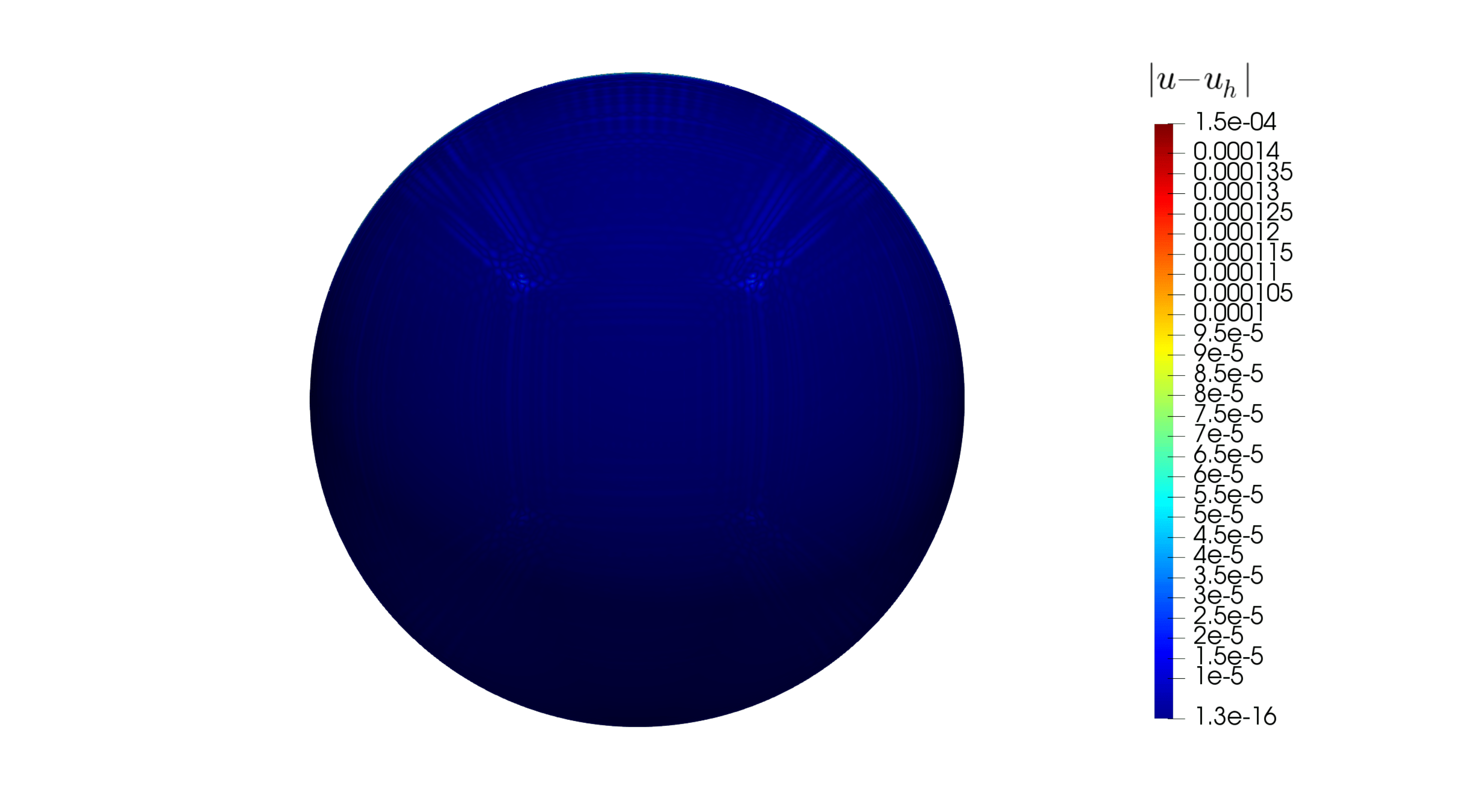}
  \caption{$2^{\text{nd}}$ level subdivision mesh with adaptive Gauss quadrature.}
  \label{fig:sphere_error_plot_aq_2}
\end{subfigure}

\caption{Point-wise error $|u- u_h|$ plots over spherical surfaces.}
\label{fig:test}
\end{figure}

\paragraph*{Quadrature error\\}
The presence of extraordinary vertices leads to difficulties in integration as described in Section~\ref{sec:aq_ev}. Figures~\ref{fig:sphere_error_plot_0},~\ref{fig:sphere_error_plot_1} and~\ref{fig:sphere_error_plot_2} show the point-wise errors at surface points for three levels of mesh refinement using the standard Gauss quadrature rule.  {The number of extraordinary vertices remains 4 after refinement.} For the analysis using the initial mesh, the error in the regions around extraordinary vertices have similar magnitudes to the other regions. However, after a level of refinement, the error in the other regions is reduced more than the area around the four extraordinary vertices. After the second refinement, the error is concentrated in the areas around the four extraordinary vertices. Figures~\ref{fig:sphere_error_plot_aq_0},~\ref{fig:sphere_error_plot_aq_1} and~\ref{fig:sphere_error_plot_aq_2} plot the point-wise errors on the same mesh analysed with the adaptive quadrature rule shown in Section~\ref{sec:aq_ev}. The errors around extraordinary vertices are now decreased. 

\begin{figure}
\centering
\captionsetup{justification=centering}

\begin{subfigure}[t]{0.5\linewidth}
\centering
\includegraphics[width=0.8\linewidth]{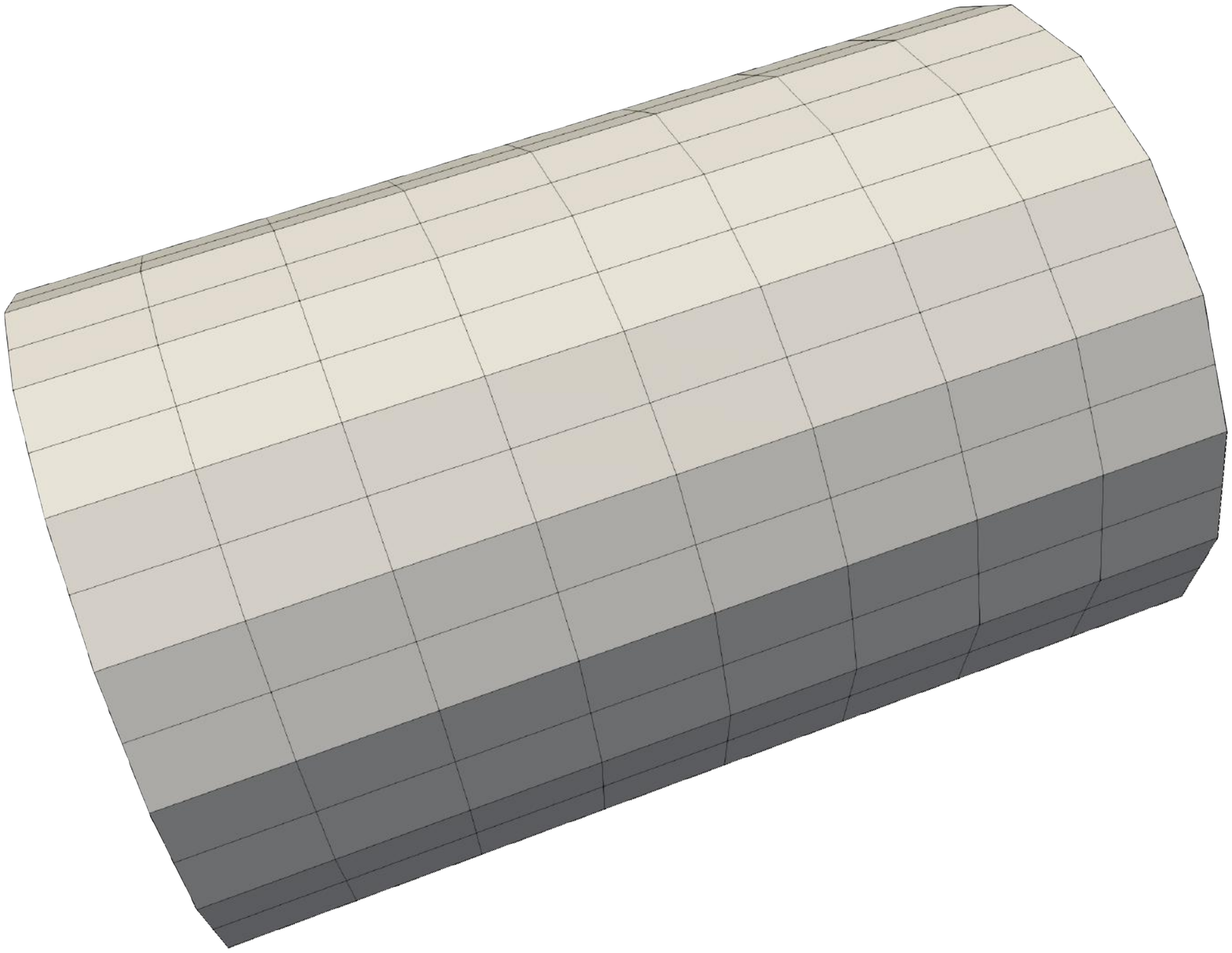}
\caption{Regular mesh has no extraordinary vertex.}
\label{fig:cylinder_mesh_exv0}
\end{subfigure}%
\begin{subfigure}[t]{0.5\linewidth}\centering
\includegraphics[width=0.8\linewidth]{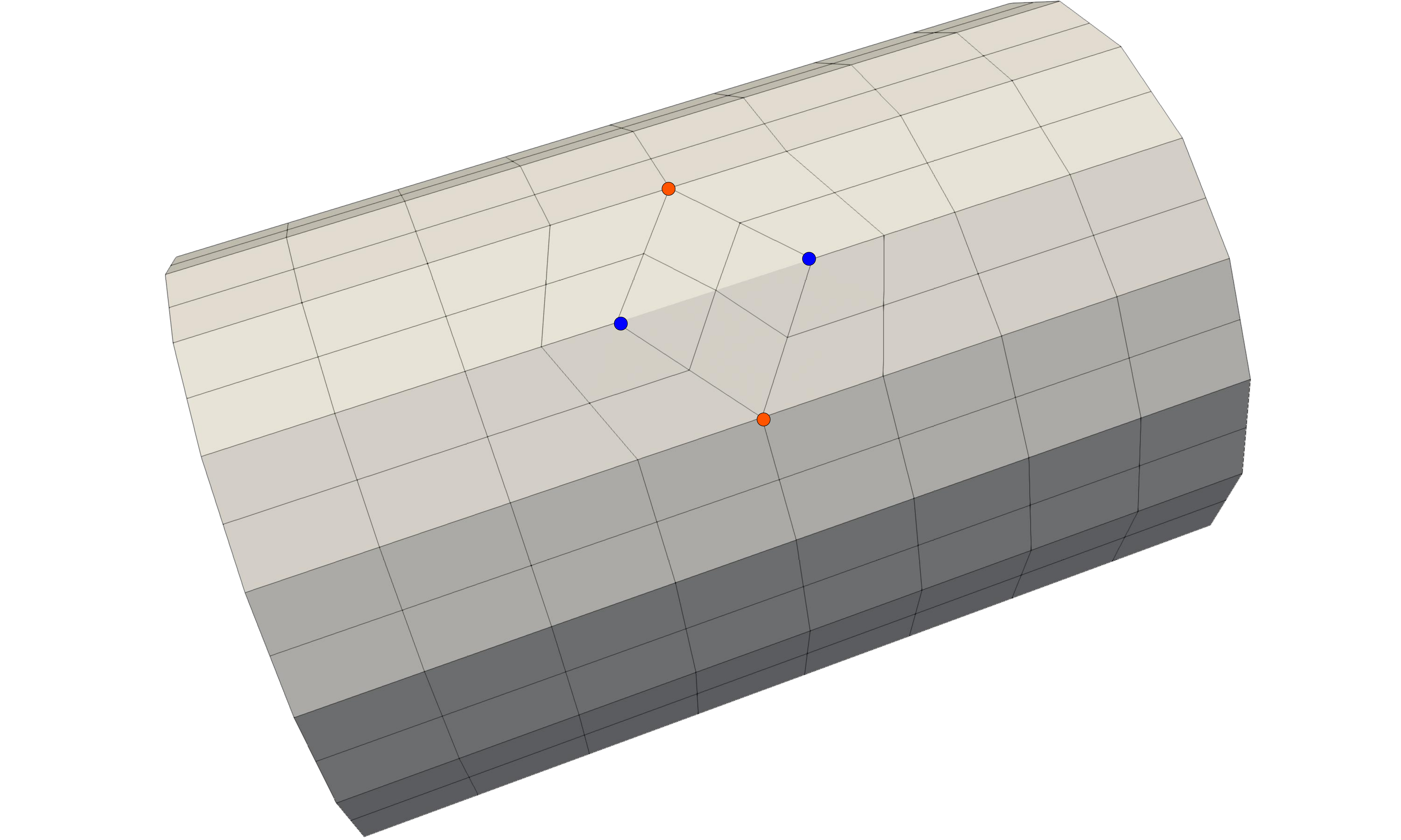}
\caption{Mesh has four extraordinary vertices (two valence 5 vertices and two valence 3 vertices).}
\label{fig:cylinder_mesh_exv4}
\end{subfigure}

\begin{subfigure}{0.7\linewidth}\centering
\includegraphics[width=\linewidth]{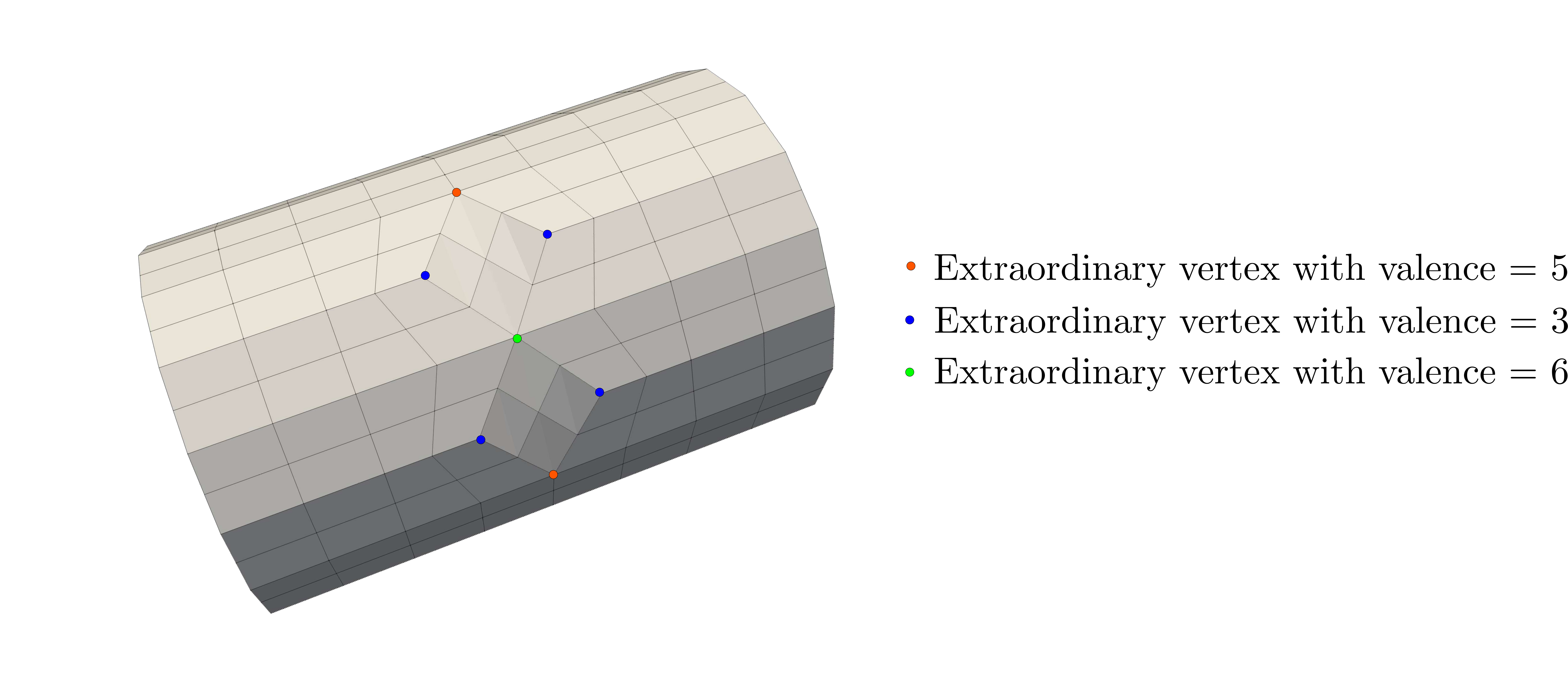}
\caption{Mesh has seven extraordinary vertices (one valence 6 vertex, two valence 5 vertices and four valence 3 vertices).}
\label{fig:cylinder_mesh_exv7}
\end{subfigure}

\label{fig:cylinder_ev_meshes}
\caption{Three control meshes with different number of extraordinary vertices.}
\end{figure}

\begin{figure}
\centering
\captionsetup{justification=centering}

\begin{subfigure}[b]{0.5\linewidth}
\centering
\includegraphics[width=\linewidth]{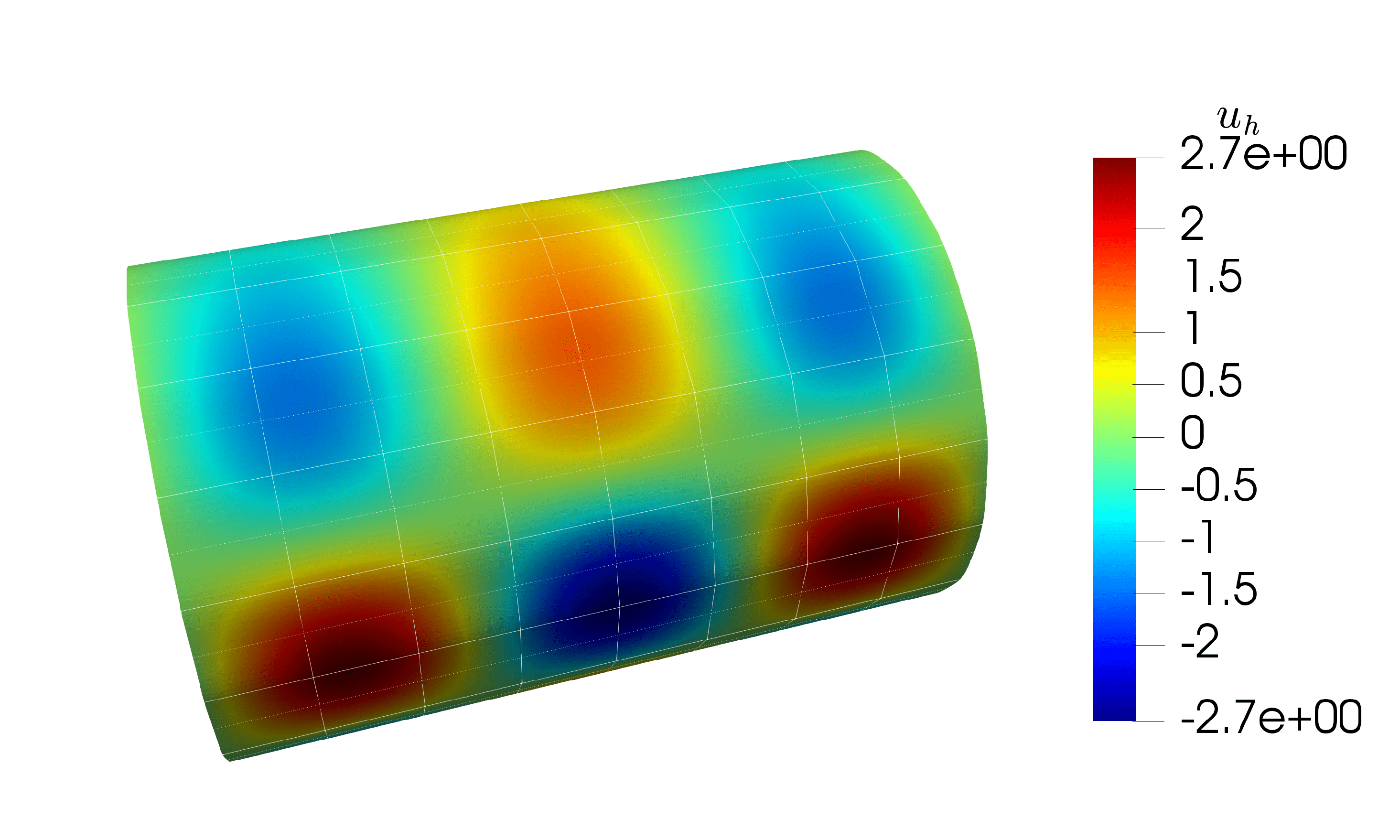}
\caption{Numerical result $u_h$ on the cylindrical surface with no extraordinary vertex.}
\label{fig:sphere_result_ev_0}
\end{subfigure}%
\begin{subfigure}[b]{0.5\linewidth}
\centering
\includegraphics[width=\linewidth]{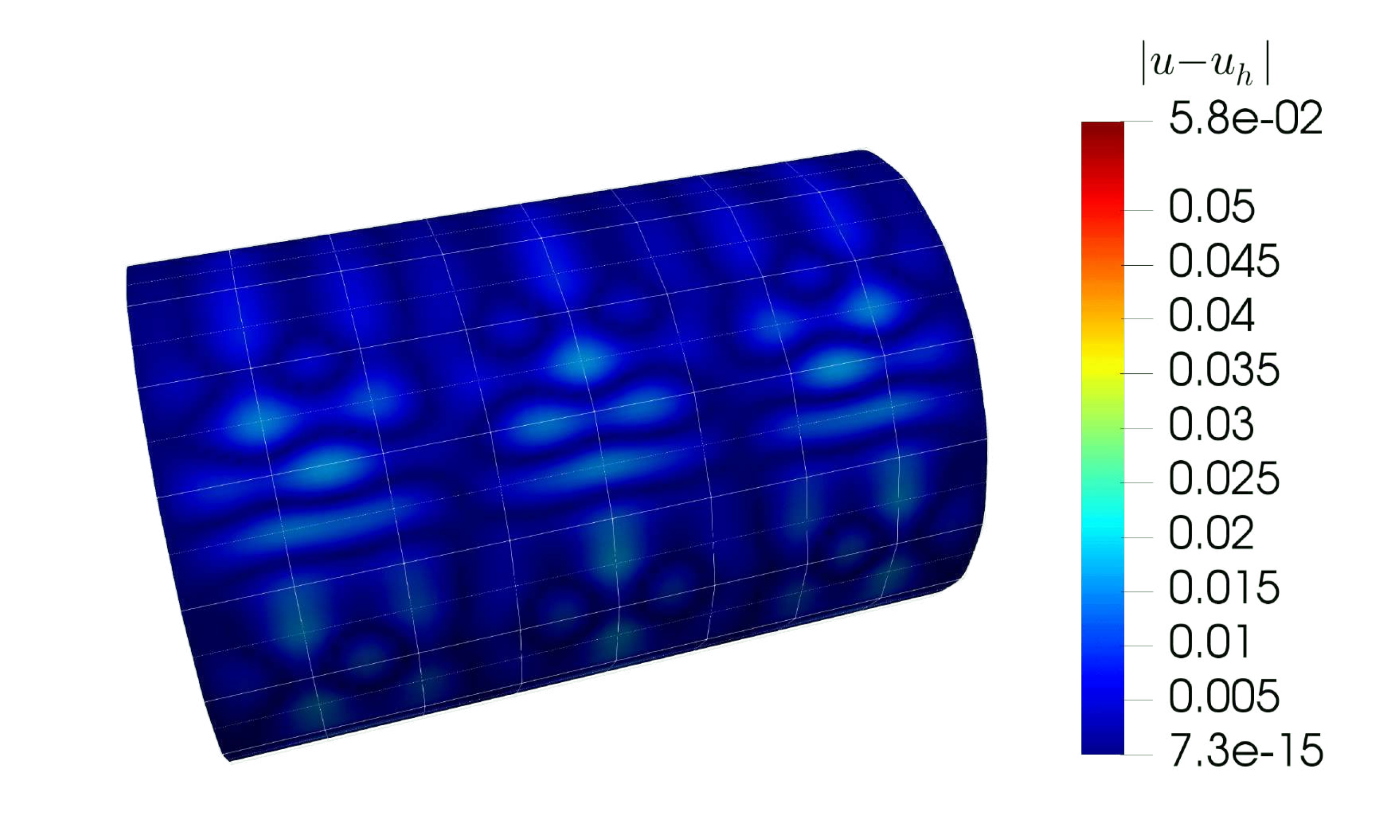}
\caption{Point-wise error on the cylindrical surface with no extraordinary vertex.}
\label{fig:cylinder_pe_exv_0}
\end{subfigure}
\begin{subfigure}[b]{0.5\linewidth}
\centering
\includegraphics[width=\linewidth]{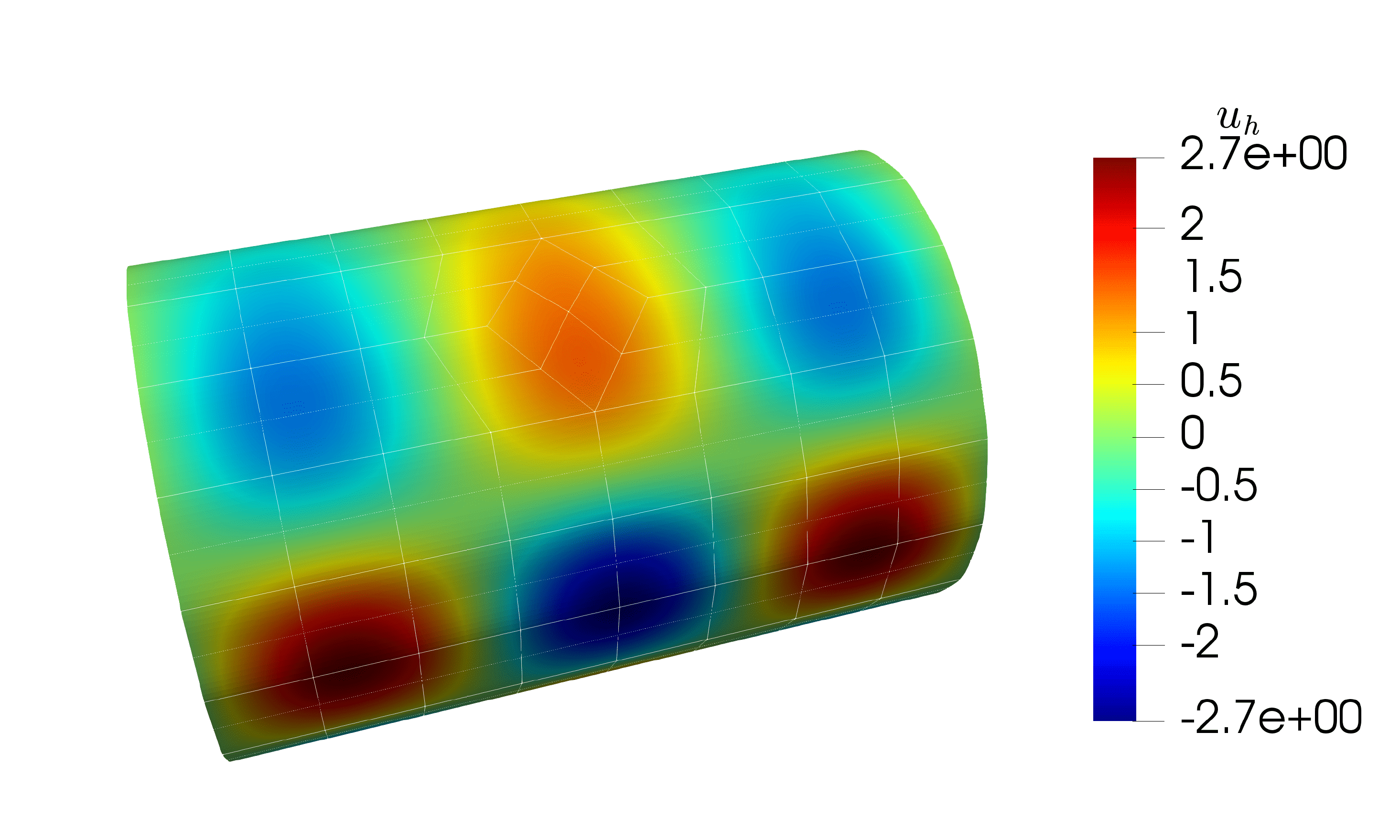}
\caption{Numerical result $u_h$ on the cylindrical surface with 4 extraordinary vertices.}
\label{fig:sphere_result_ev_4}
\end{subfigure}%
\begin{subfigure}[b]{0.5\linewidth}
\centering
\includegraphics[width=\linewidth]{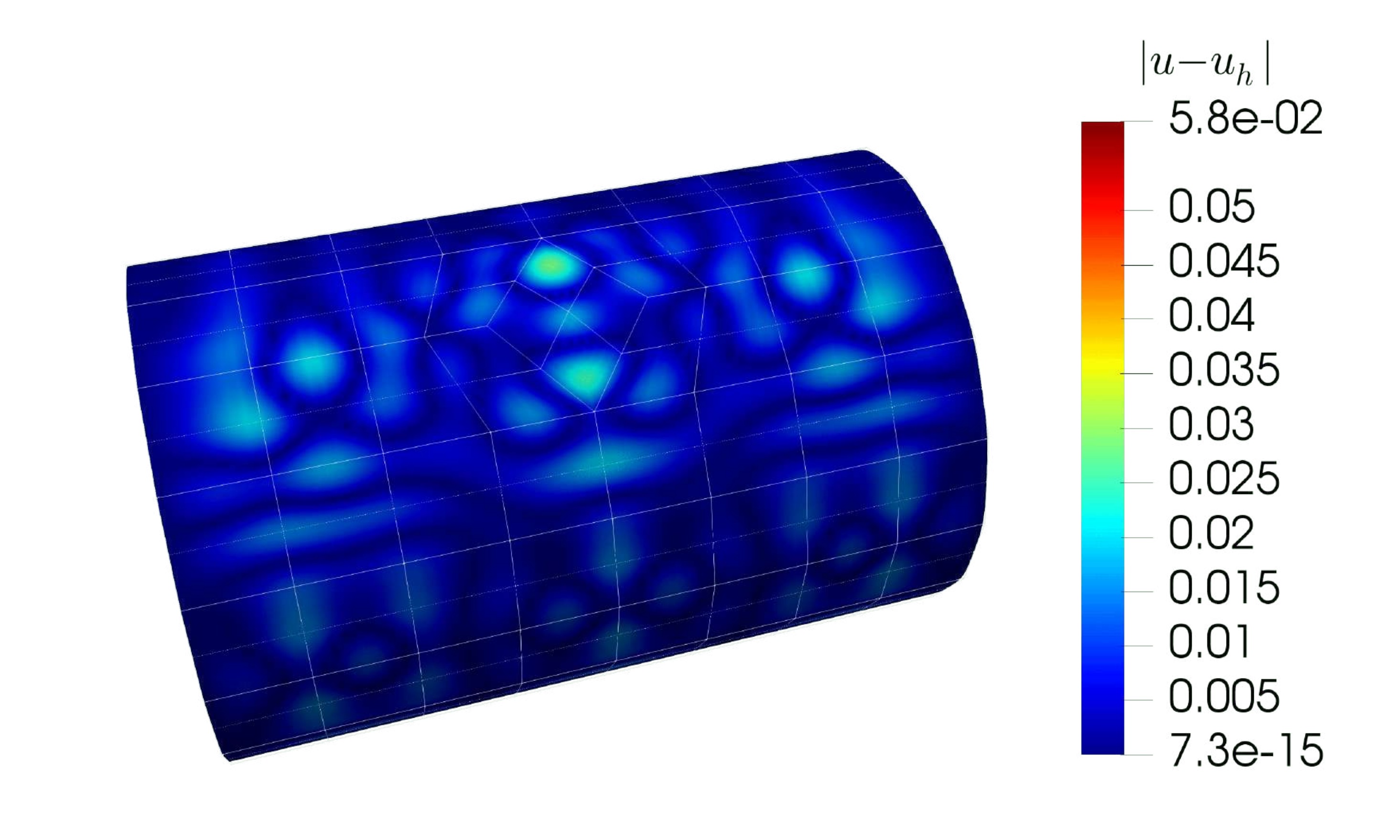}
\caption{Point-wise error on the cylindrical surface with 4 extraordinary vertex.}
\label{fig:cylinder_pe_exv_4}
\end{subfigure}
\begin{subfigure}[b]{0.5\linewidth}
\centering
\includegraphics[width=	\linewidth]{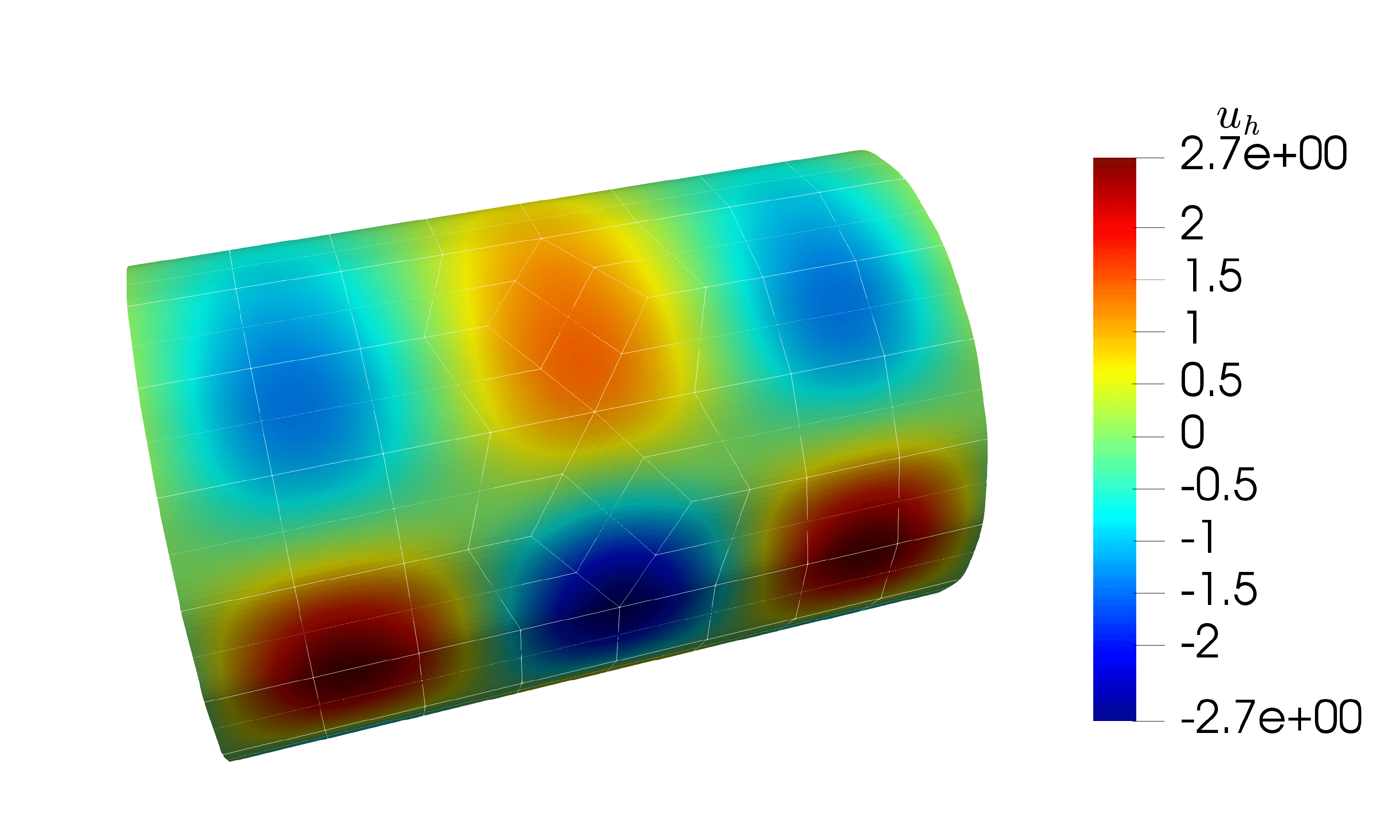}
\caption{Numerical result $u_h$ on the cylindrical surface with 7 extraordinary vertices.}
\label{fig:sphere_result_ev_7}
\end{subfigure}%
\begin{subfigure}[b]{0.5\linewidth}
\centering
\includegraphics[width=\linewidth]{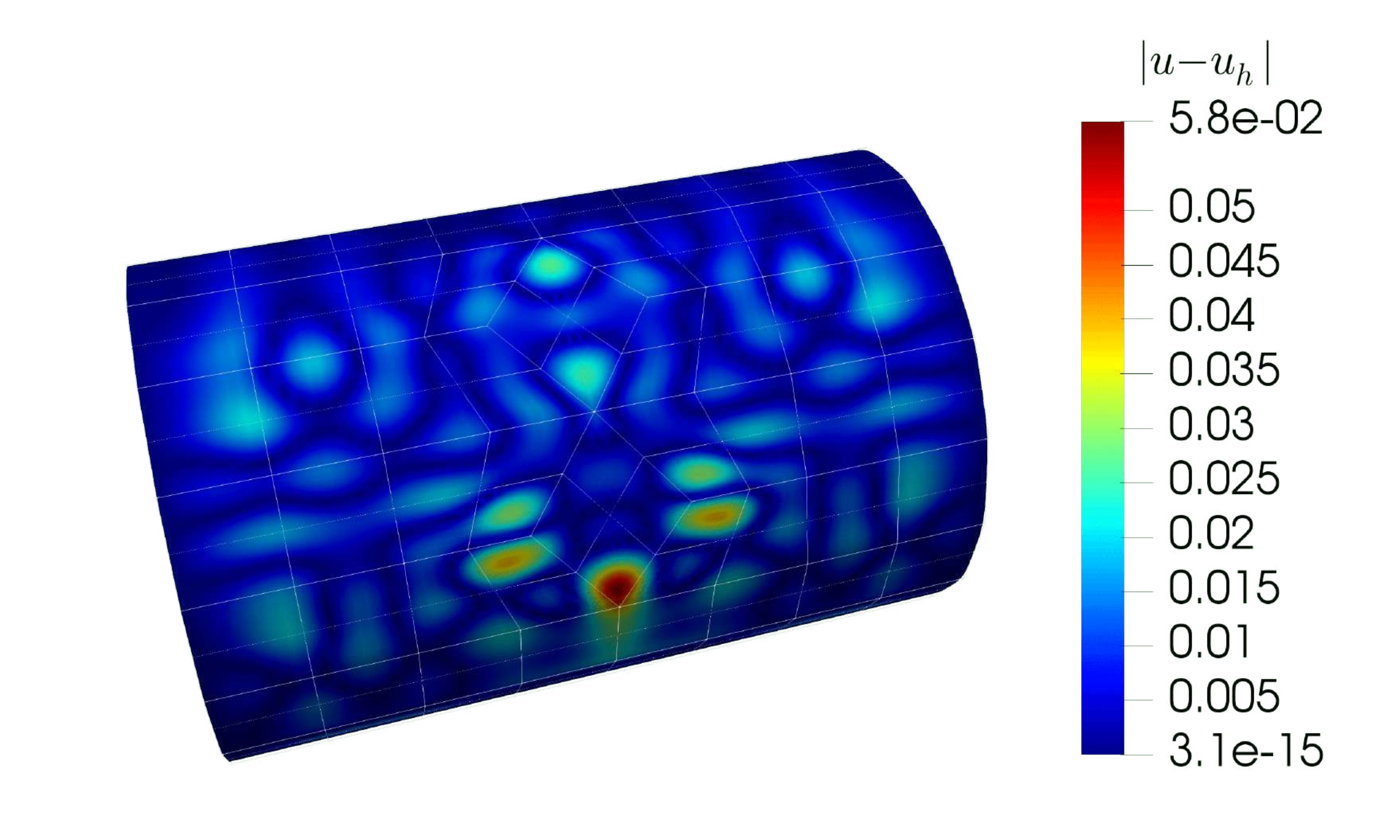}
\caption{Point-wise error on the cylindrical surface with 7 extraordinary vertex.}
\label{fig:cylinder_pe_exv_7}
\end{subfigure}
\label{fig:cylinder_ev_meshes_errors}
\caption{Numerical results and point-wise errors on the cylindrical surfaces constructed by thee different meshes. The white grids are used to indicate the locations of extraordinary vertices.}\end{figure}

\begin{figure}
\centering
\includegraphics[width=0.9\linewidth]{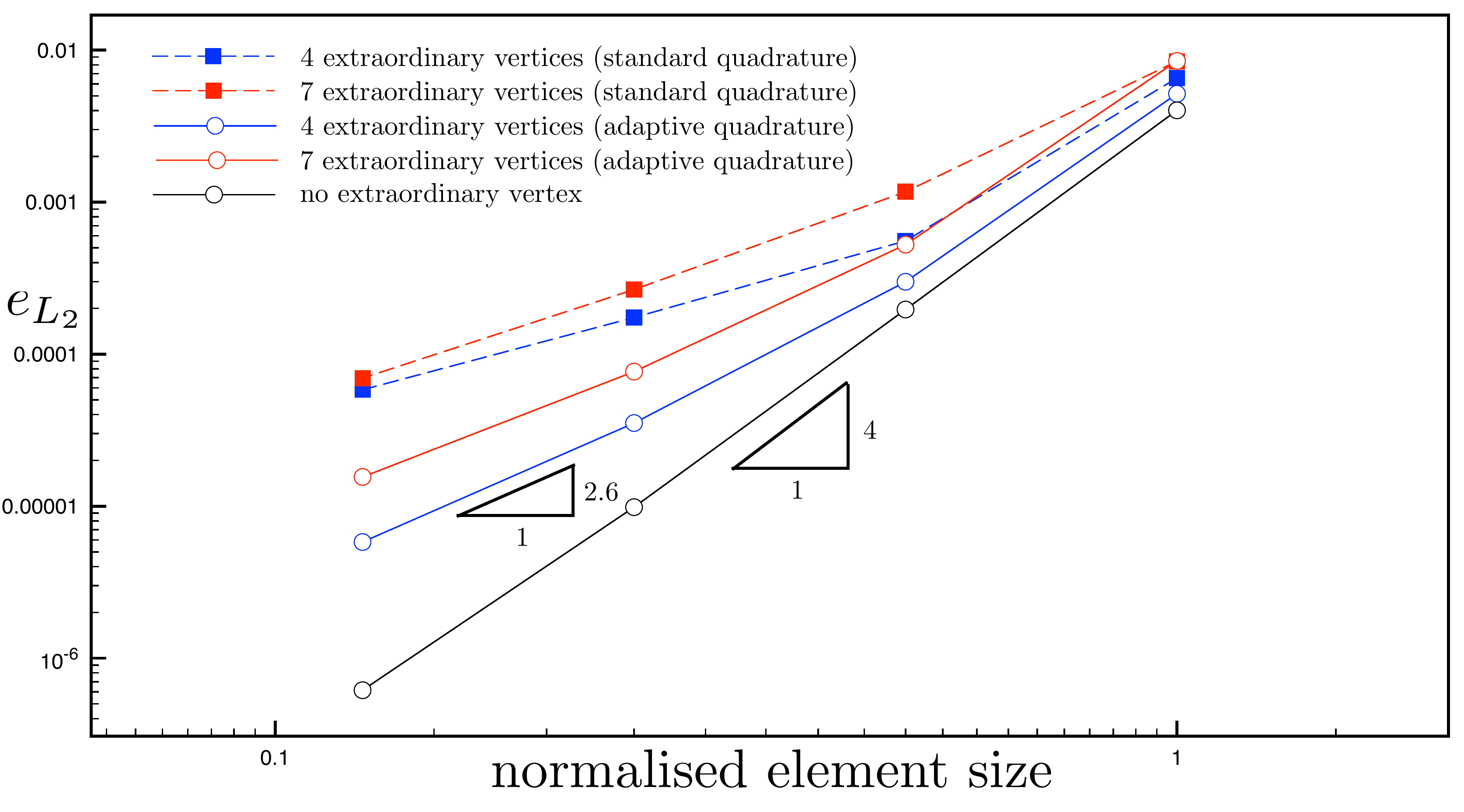}
\caption{Convergence study: comparison of the Catmull-Clark subdivisions with different number of extraordinary vertices. The adaptive Gauss quadrature achieves better convergence rates over the standard Gauss quadrature rule.}
\label{fig:cylinder_convergence_exv}
\end{figure}

\paragraph*{Approximation error\\}
{The presence of extraordinary vertices introduces approximation errors. Then we investigate the effect of the number and valence of extraordinary vertices on numerical accuracy. } Figures~\ref{fig:cylinder_mesh_exv0},~\ref{fig:cylinder_mesh_exv4} and~\ref{fig:cylinder_mesh_exv7} are three control meshes with different numbers of extraordinary vertices. Figure~\ref{fig:cylinder_mesh_exv0} shows a control mesh without an extraordinary vertex. Figure~\ref{fig:cylinder_mesh_exv4} shows a control mesh with four extraordinary vertices, including two vertices with a valence of 3 and two vertices with a valence of 5. The control mesh in Figure~\ref{fig:cylinder_mesh_exv7} has seven extraordinary vertices, including four vertices with a valence of 4, two vertices with a valence of 5 and one vertex with a valence of 6. It is important to note the three different control meshes construct different but similar geometries.  The Laplace-Beltrami problem is solved using the Galerkin formulation with the same right-hand side function $f$ computed in~\eqref{eq:bt_rhs}. Both standard and adaptive Gauss quadrature rules are used for all cases. Figures~\ref{fig:sphere_result_ev_0}, \ref{fig:sphere_result_ev_4} and~\ref{fig:sphere_result_ev_7} show the solution of $u$ on the surfaces constructed using the three meshes. Because of the similarity of the geometries and solutions, the three cases are used to investigate the effects of the number of extraordinary vertices in numerical results. 
The point-wise errors on the three surfaces are shown in Figures~\ref{fig:cylinder_pe_exv_0}, ~\ref{fig:cylinder_pe_exv_4} and~\ref{fig:cylinder_pe_exv_7}. Meshes with extraordinary vertices have larger maximum point-wise errors close to the extraordinary vertices, while the mesh without extraordinary vertices has increased uniform point-wise error. 
Figure~\ref{fig:cylinder_convergence_exv} shows the convergence rates for the three cases. Meshes without extraordinary vertices can achieve the optimal $p+1$ convergence rate and $p=3$. In general, the more extraordinary vertices a mesh contains, the more error results. The extraordinary vertices increase the global errors in the results and reduce the convergence rate. {Since the global errors also include quadrature errors, the adaptive quadrature rule serves to reduce the quadrature errors. With the adaptive quadrature rule, the convergence rates are improved for the 4 and 7 extraordinary vertices cases but the results still agree with our assumption that increasing the number and valence of extraordinary vertices will produce higher error. } 
\paragraph*{Computational cost\\}
{Table~\ref{tab:computational_time} compares the computational cost for assembling the system matrix for the standard and adaptive quadrature rules. Because the number of extraordinary vertices remains constant after subdivision, the difference in computational time between the standard and adaptive quadrature schemes diminishes. }

{
\begin{table}[]
\centering
\begin{tabular}{@{}crrrr@{}}
\toprule
Number of  &   &   &   &   \\extraordinary vertices                & \multirow{-2}{*}{$n_e$} & \multirow{-2}{*}{$n_d$}  & \multirow{-2}{*}{$L_2$ error} & \multirow{-2}{*}{Assembly time}                  \\ \midrule
                                 &                    &      0              & {6.55e-3}     & {2.46s}   \\
                                 &                    &       3                & 6.71e-3                            & 3.98s                          \\
                                 &    \multirow{-3}{*}{{260}}                 & 7  & {6.72e-3}  & {6.39s}   \\ \cmidrule(l){2-5} 
                                 &                    &        0            & {5.93e-4}     & {9.39s}   \\
                                 &                    &        3              & 3.48e-4                            & 10.16s                         \\
                                 & \multirow{-3}{*}{{1040}}                   & 7  & {3.49e-4} & {13.87s}  \\ \cmidrule(l){2-5} 
                                 &                    & {0}                       & {1.80e-4}     & {34.9s}   \\
                                 &                    & {3}                       & {5.24e-5}     & {35.51s}  \\
\multirow{-9}{*}{4}              & \multirow{-3}{*}{{4160}}                   &  7 & {5.31e-5}     & {37.02s}  \\ \midrule
                                 & {264}          &    7      & {8.49e-3}     & {9.59s}  \\
                                 & {1056}         &    7      & { 5.54e-4}     & {16.03s} \\
\multirow{-3}{*}{7}               & {4224}    & {7}               & {8.60e-5}     & {40.24s} \\ \bottomrule
\end{tabular}
\caption{The computational times for assembling the system matrix are recorded for tests on cylindrical meshes with different number of extraordinary vertices. Standard quadrature ($n_d = 0$) is compared with adaptive quadratures with $n_d = 3$ and $7$. }
\label{tab:computational_time}
\end{table}

\subsection{Complex geometry}
This final example considers the ability of the Catmull-Clark method to provide high-order discretisations of complex geometry. The model considered is that of a racing car from CAD and imported into Autodesk Maya~\cite{maya} for removal of extraneous geometry and the generation of the surface mesh shown in Figure~\ref{fig:car_control_mesh}. Modelling such geometry using NURBS surfaces would require a number of patches to be spliced together. A model based on Catmull-Clark subdivision surface can directly evaluate the smooth limit surface in Figure~\ref{fig:car_limit_surfaces} using the control mesh containing extraordinary vertices. The minimum bounding box for this model is defined by $[x_i^{\min}, x_i^{\max}]^3 = [-1.047,1.122]\times[0.097,0.692]\times[-0.460,0.460]$. The control mesh has 9154 vertices. The physical domain is naturally discretised into a number of elements expressed as
\begin{align}
\Gamma = \bigcup_{e=1}^{n_{e}}\Gamma_e,
\end{align}
where $n_e = 9152$ for this example. Figure~\ref{fig:car_bc} indicates the domain where essential (Dirichlet) and natural boundary conditions are applied. The essential boundary $\Gamma_d$  is composed of two parts as 
\begin{align}
\Gamma_d = \Gamma^1_d + \Gamma^1_d,
\end{align}
where
\begin{align}
\Gamma^1_d &= \{\Gamma_e \,|\, x_2 < -0.9 \quad \forall \mathbf{x}  \in \Gamma_e\},\\
\Gamma^2_d &= \{\Gamma_e\, | \,x_2 > 0.9 \quad \forall \mathbf{x} \in \Gamma_e\}.
\end{align}
The natural boundary conditions is applied to the rest of the domain $\Gamma_n = \Gamma \backslash \Gamma_d$. The numerical result matches the analytical solution well as shown in Figure~\ref{fig:car_solution}. Figure~\ref{fig:car_xy} shows the results on $\Gamma_n$ and a maximum point-wise error $2.8\%$ is observed in Figure~\ref{fig:car_pointwise_error}.
\begin{figure}
\centering
\begin{subfigure}{0.49\linewidth}
\includegraphics[width=\linewidth]{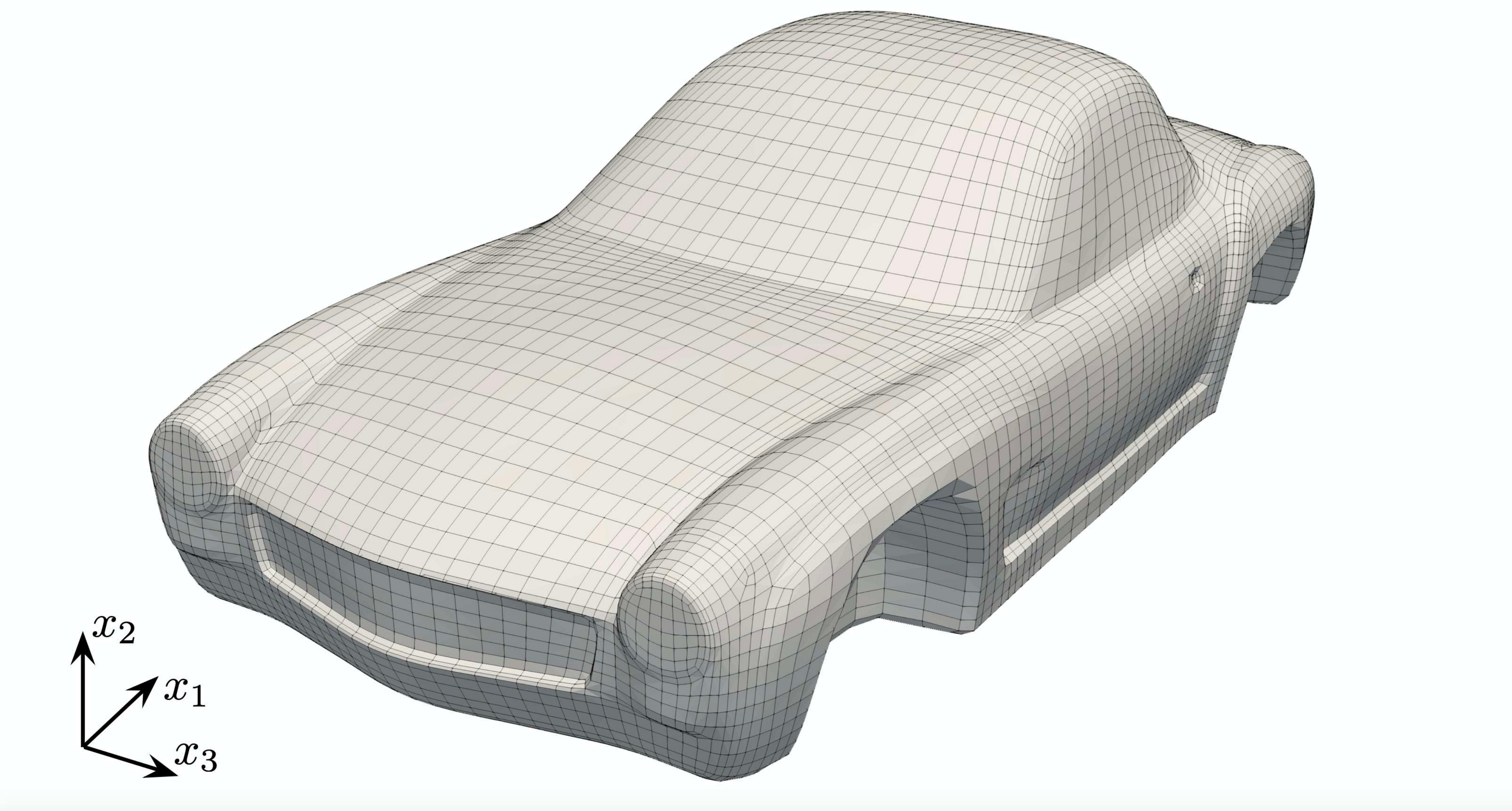}
\caption{A control mesh with 9154 vertices.}
\label{fig:car_control_mesh}
\end{subfigure}
\begin{subfigure}{0.49\linewidth}
\includegraphics[width=\linewidth]{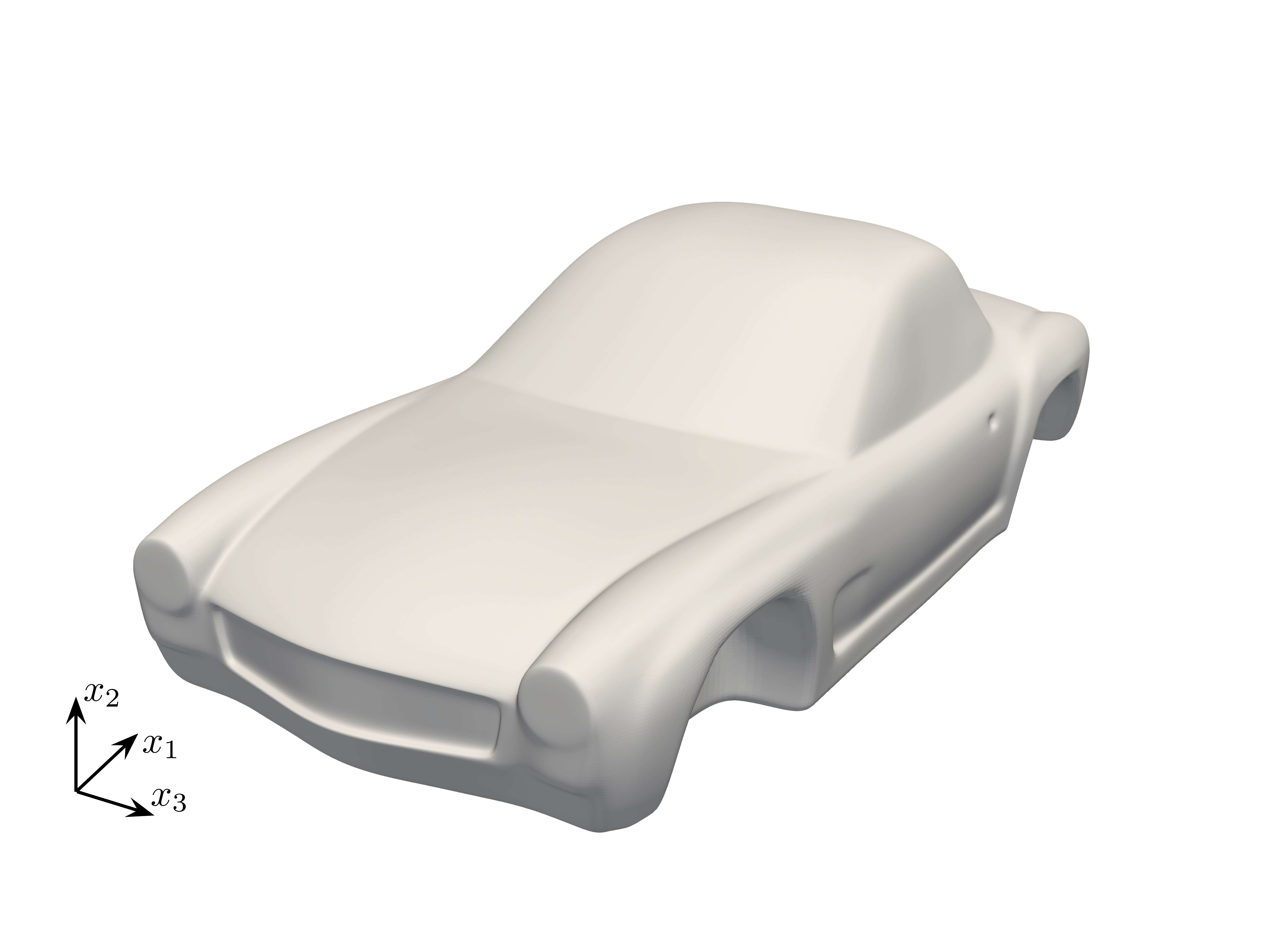}
\caption{Limit surface.}
\label{fig:car_limit_surfaces}
\end{subfigure}
\begin{subfigure}{0.49\linewidth}
\includegraphics[width=\linewidth]{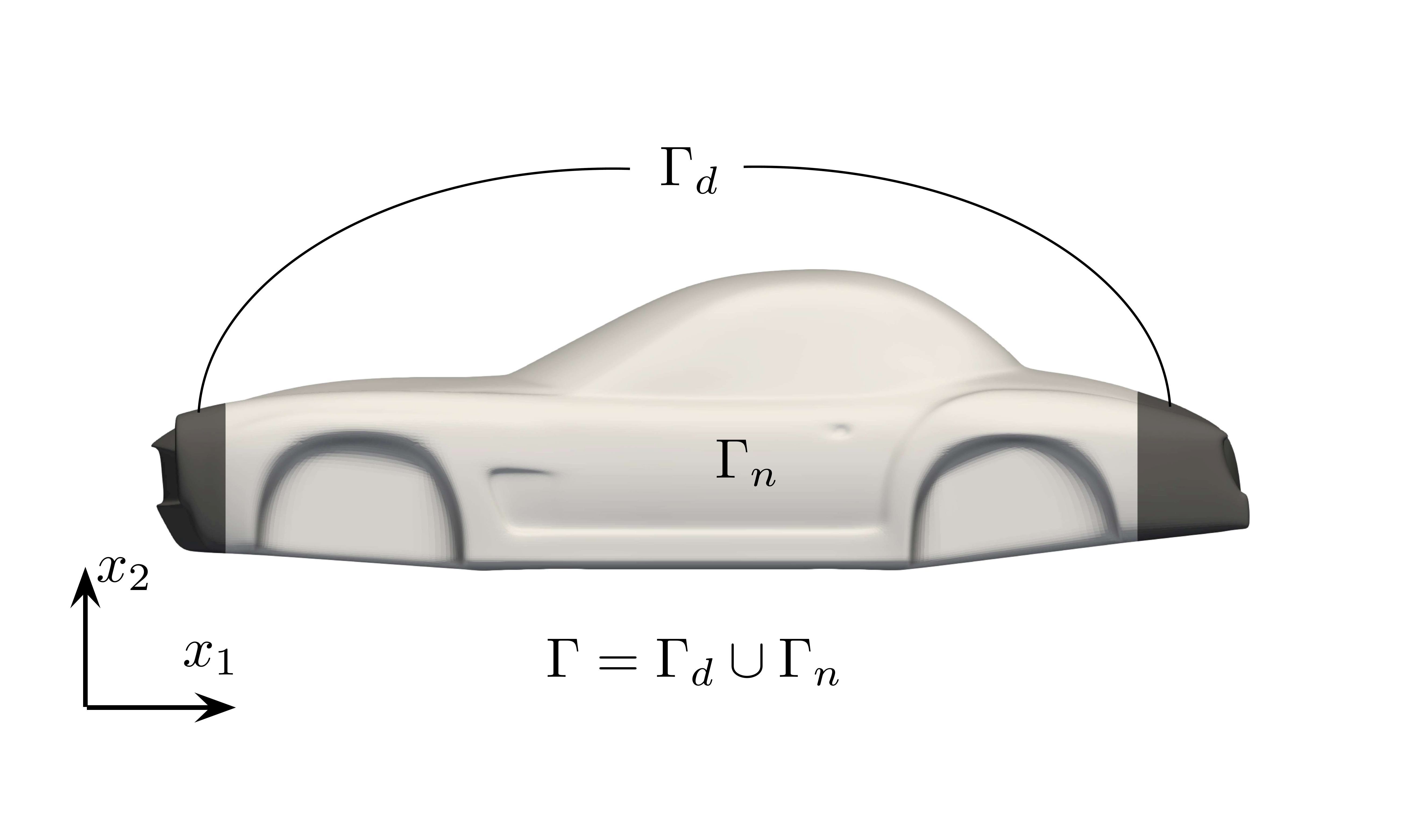}
\caption{Domains defined for applying boundary conditions.}
\label{fig:car_bc}
\end{subfigure}
\begin{subfigure}{0.49\linewidth}
\includegraphics[width=\linewidth]{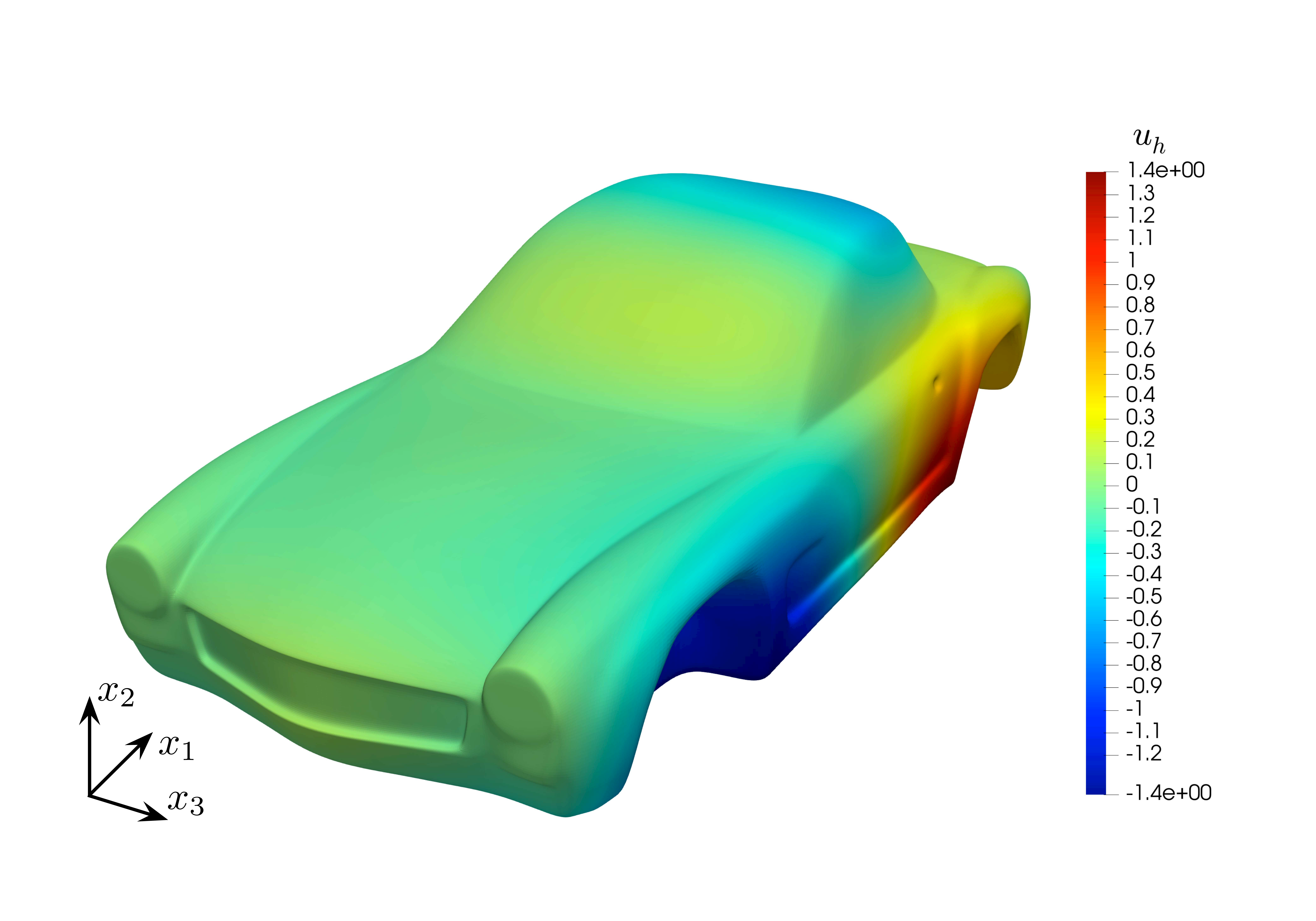}
\caption{Numerical results on the car surface.}
\label{fig:car_solution}
\end{subfigure}
\begin{subfigure}{0.49\linewidth}
\includegraphics[width=\linewidth]{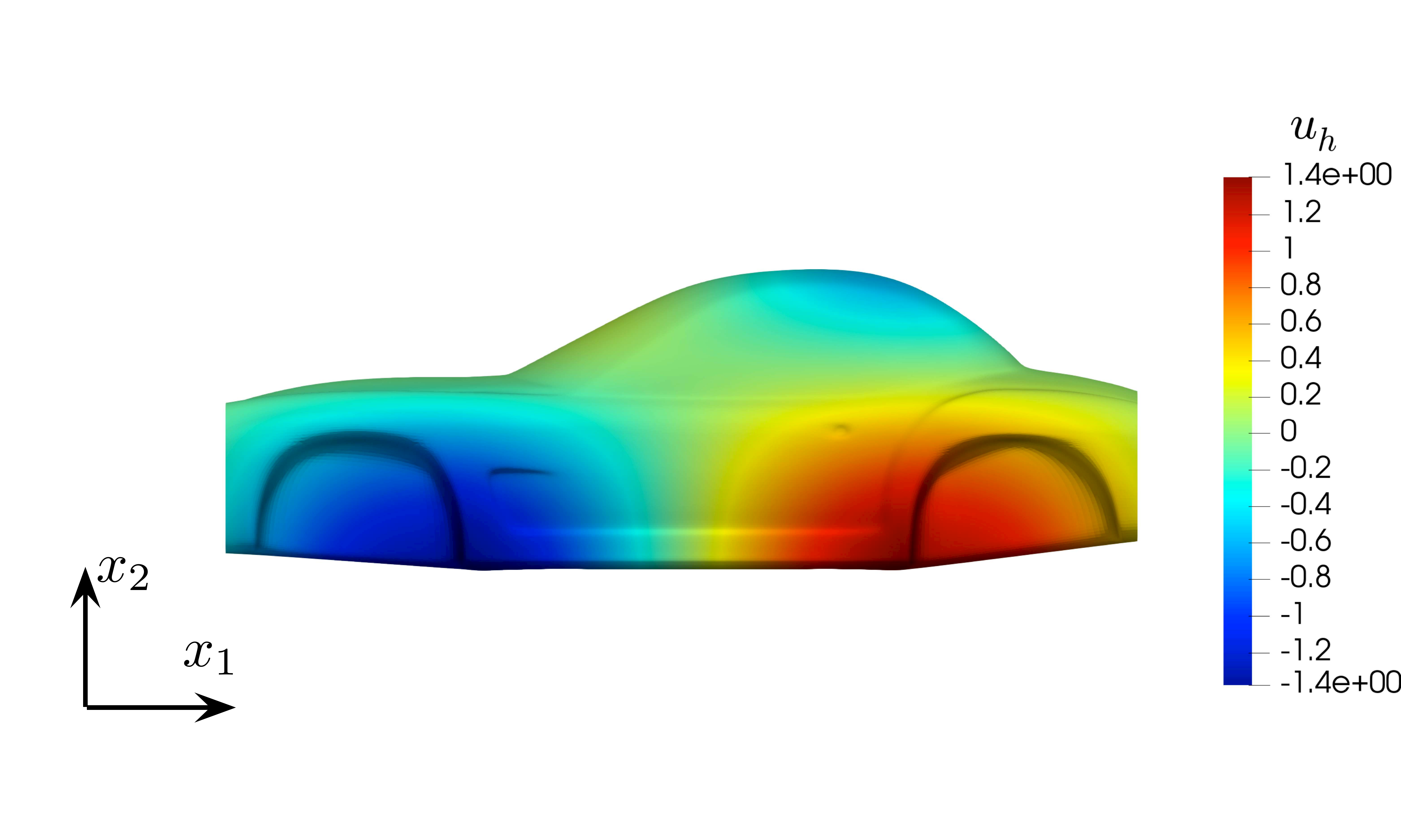}
\caption{Numerical results on $\Gamma_n$.}
\label{fig:car_xy}
\end{subfigure}
\begin{subfigure}{0.49\linewidth}
\includegraphics[width=\linewidth]{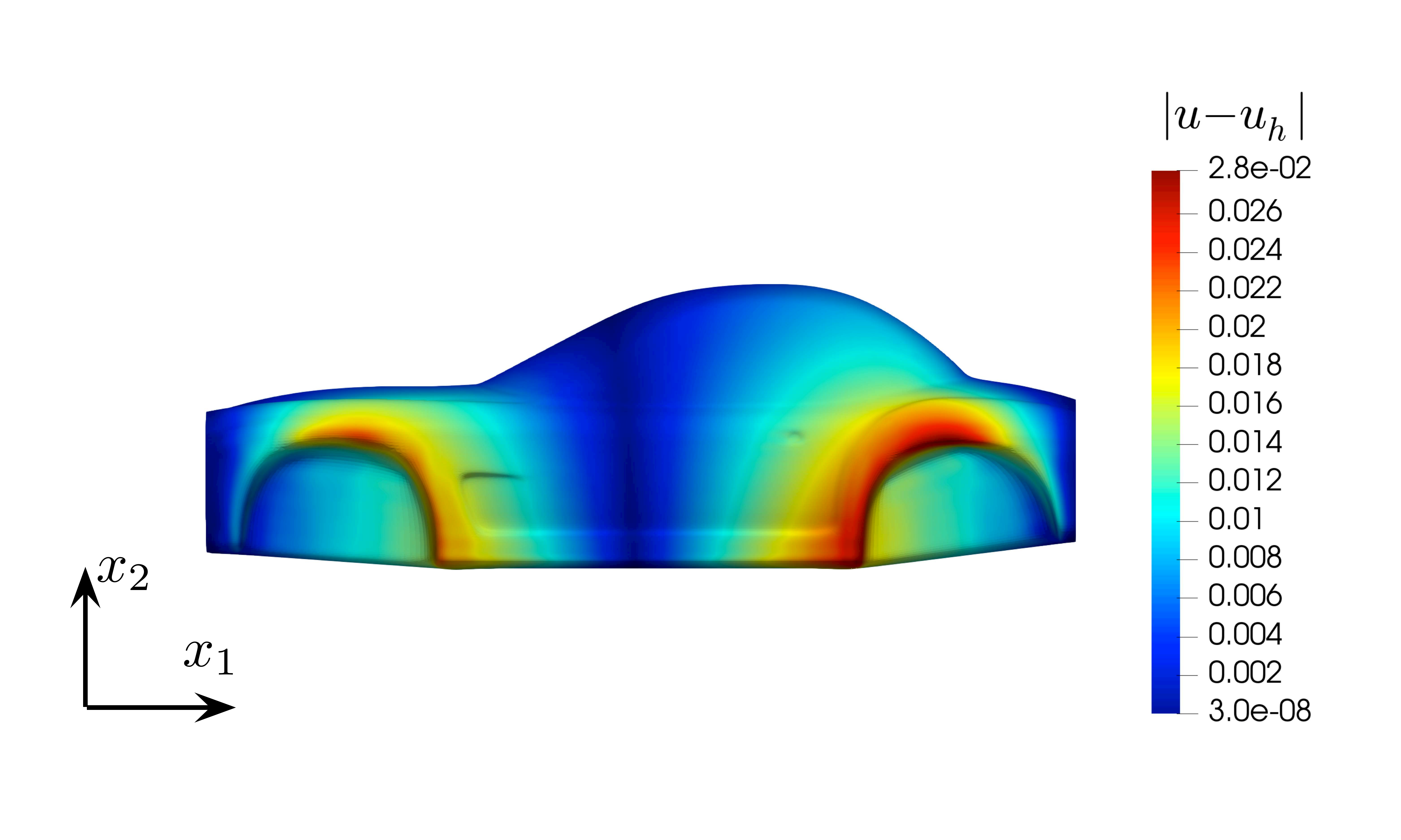}
\caption{Pointwise error on $\Gamma_n$.}
\label{fig:car_pointwise_error}
\end{subfigure}
\caption{The Laplace-Beltrami problem is sloved on a complex car geometry.}
\end{figure}
}
\section{Conclusions}
{ A thorough study of the isogeometric Galerkin method with Catmull-Clark subdivision surfaces has been presented. The same bases have been used for both geometry and the Galerkin discretisation. The method has been used to solve the Laplace-Beltrami equation on curved two-dimensional manifold embedded in three dimensional space using the Catmull-Clark subdivision surfaces. An approach to fit given geometries using Catmull-Clark subdivision scheme has been outlined.  A method to model open boundary geometries without involving `ghost' control vertices, but involving errors in function gradients close to boundary regions, has also been described. The penalty method has been adopted to impose the Dirichlet boundary conditions. The optimal convergence rate of $p+1$ has been obtained when using a cylindrical control mesh without extraordinary vertices. A reduction of convergence rates has been observed when the function gradients at the boundaries do not behave like constant, or control meshes contain extraordinary vertices. The adaptive quadrature scheme significantly improves the accuracy. The effect of the number and valence of the extraordinary vertices in convergence rates has been investigated and an adaptive quadrature rule implemented. This successfully improved the convergence rates for the proposed method. The convergence rate of the proposed method is not worse than  $2.5$ ($L_2$ error) and $1.5$ ($H^1$ error).}  

In future work, this method will be investigated with problems requiring $C^1$ continuity such as the deformations of thin shells.

\section*{Acknowledgements}
This work was supported by the UK Engineering and Physical Sciences Research Council grant EP/R008531/1 for the Glasgow Computational Engineering Centre.\\
\bibliographystyle{spbasic}      



\bibliography{cas-refs}
\appendix
\section{Appendix}
\label{ap:A}
\begin{figure}
\centering
	\includegraphics[width=0.7\linewidth]{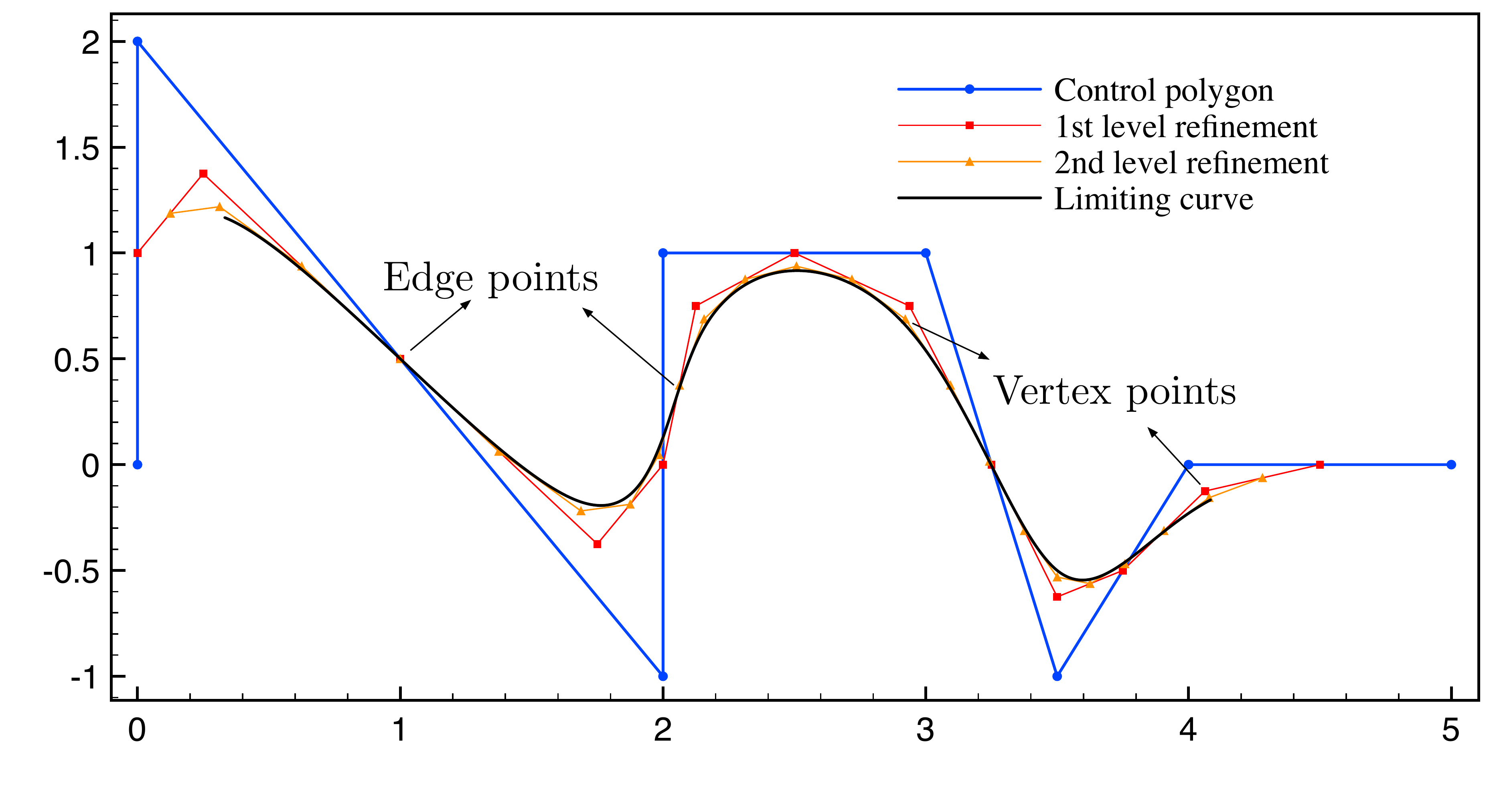}
	\caption{Schematics of the Lane-Riesenfeld Algorithm}
	\label{fig:Catmull_Clark_algorithm}
\end{figure}

{\subsection{Lane-Riesenfeld subdivision algorithm for curves}}
\label{ap:A.1}
 {The Lane-Riesenfeld algorithm} successively refines a curve starting from an initial control polygon. After a number of subdivisions, the curve is limited to a B-spline.
\begin{figure}
\centering
\includegraphics[width=0.8\linewidth]{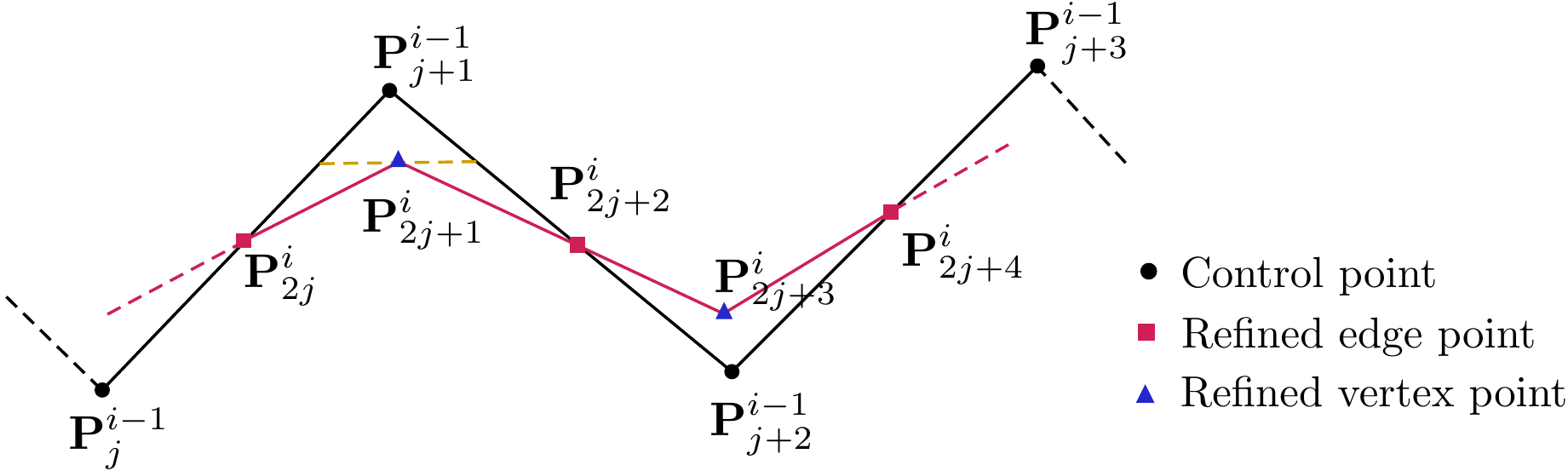}
\caption{Computing new control points through the Lane-Riesenfeld algorithm.}
\label{fig:catmull_clark_example}
\end{figure}
Figure~\ref{fig:catmull_clark_example} illustrates {a special case of this subdivision algorithm.} The control point $\mathbf{P}^i_{2j}$ in the $i^\text{th}$ level of refinement is computed from the upper level control points as
\begin{equation}
\mathbf{P}^i_{2j} = \frac{1}{2} \mathbf{P}^{i-1}_{j} + \frac{1}{2} \mathbf{P}^{i-1}_{j+1}.
\end{equation}
Point $\mathbf{P}^i_{2j}$ is the mid- point of $\mathbf{P}^{i-1}_{j}\mbox{--}\mathbf{P}^{i-1}_{j+1}$, and is called an `edge point'. The control point $\mathbf{P}^i_{2j+1}$ is computed as
\begin{equation}
\mathbf{P}^i_{2j+1} = \frac{1}{8} \mathbf{P}^{i-1}_{j} + \frac{3}{4} \mathbf{P}^{i-1}_{j+1} + \frac{1}{8} \mathbf{P}^{i-1}_{j+2}.
\end{equation}
To compute this point, one needs to connect the mid-points of $\mathbf{P}^{i}_{2j}\mbox{--}\mathbf{P}^{i-1}_{j+1}$ and $\mathbf{P}^{i-1}_{j+1}\mbox{--}\mathbf{P}^{i}_{2j+2}$. The point $\mathbf{P}^i_{2j+1}$ is the mid-point of the connecting line. This type of point is called `vertex point'. Each `vertex point' is associated with an upper level control point. Figure~\ref{fig:Catmull_Clark_algorithm} shows two levels of refinements using the {Lane-Riesenfeld algorithm} and the limiting result which is a {cubic B-spline curve}.

\subsection{Catmull-Clark subdivision algorithm for surfaces}
\label{ap:A.2}
\begin{figure}
\centering
	\includegraphics[width=0.7\linewidth]{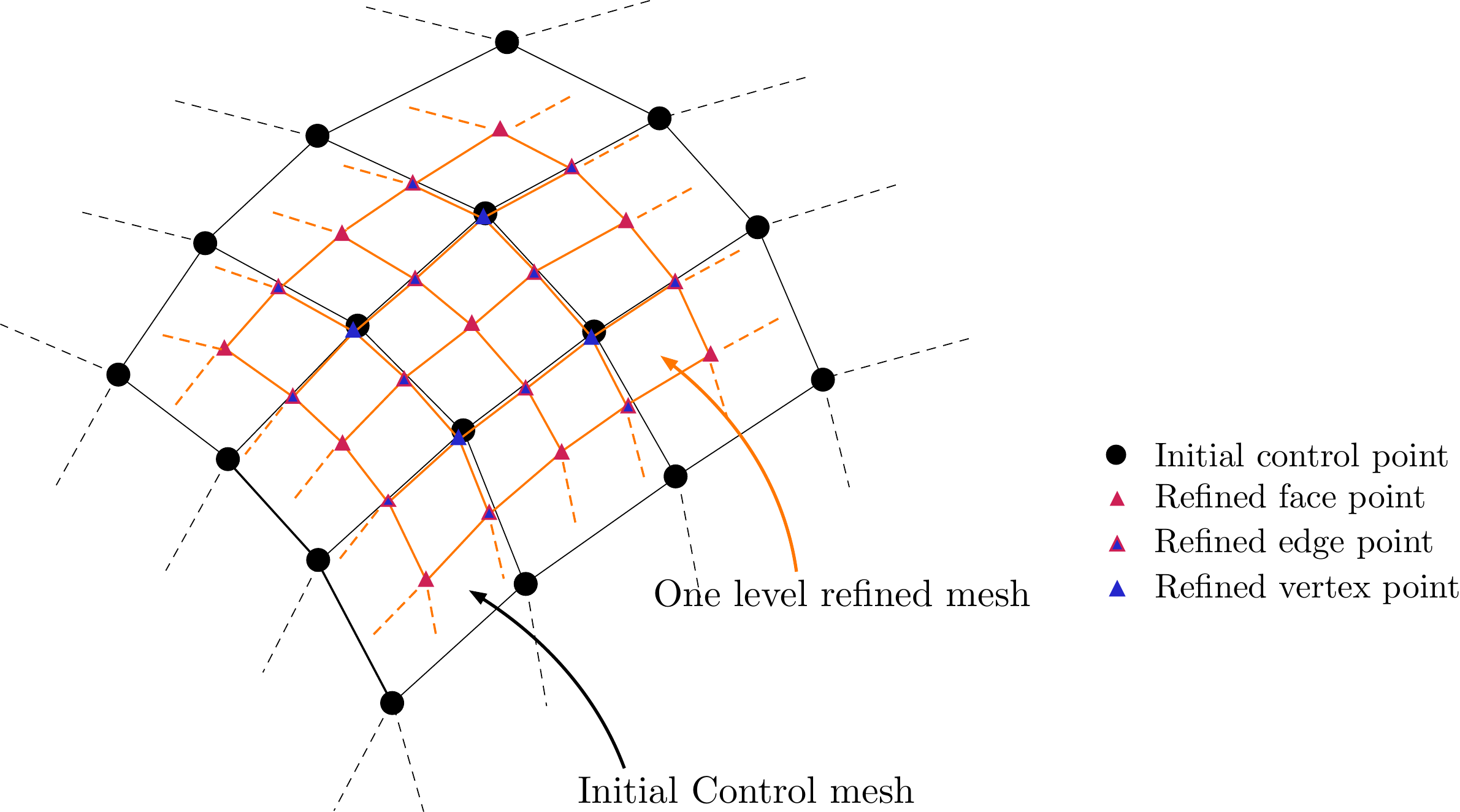}
	\caption{Catmull-Clark subdivision algorithm for surfaces.}
	\label{fig:CC_surface}
\end{figure}

The application of the subdivision algorithm to surfaces follows in a similar manner to curves. One face in the original control mesh is split into four new faces. For a closed surface, the numbers of faces and control vertices are doubled. Figure~\ref{fig:CC_surface} shows an example of generating a new control mesh through the Catmull-Clark algorithm.

\begin{figure}[]
\centering
	\begin{subfigure}[b]{0.25\linewidth}
		\includegraphics[width=\linewidth]{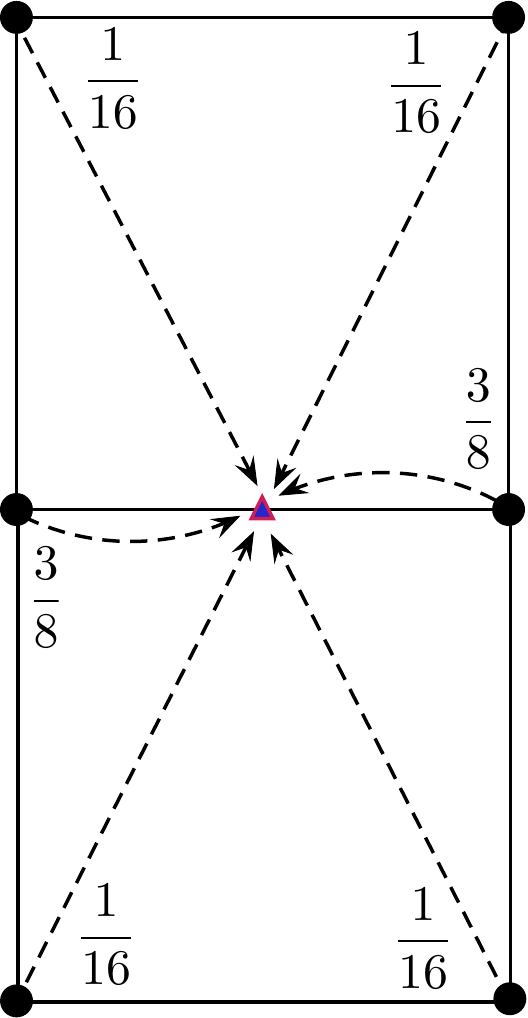}
		\caption{Edge point case 1.}
		\label{fig:CC_edge_point_1}
	\end{subfigure}
	\begin{subfigure}[b]{0.48\linewidth}
		\includegraphics[width=\linewidth]{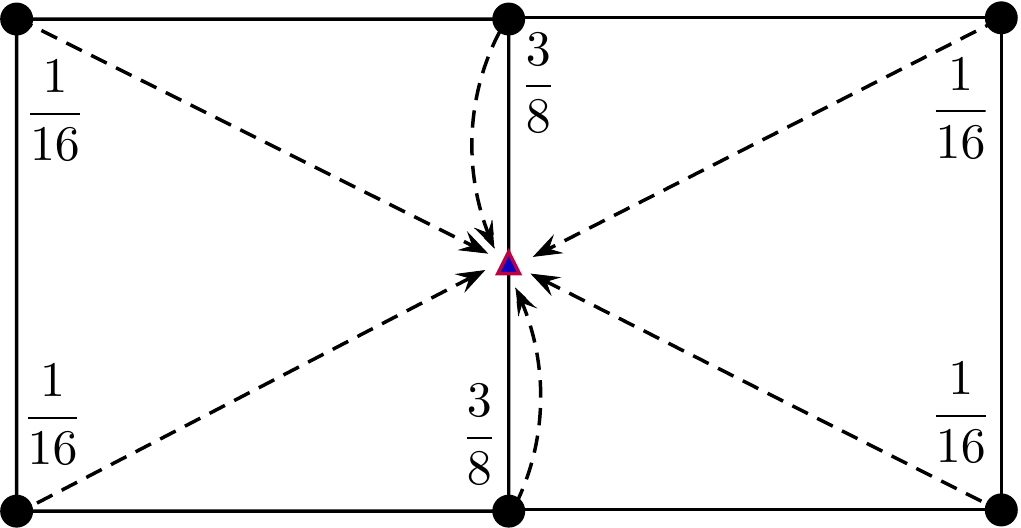}
		\caption{Edge point case 2.}
		\label{fig:CC_edge_point_2}
	\end{subfigure}
	\vspace{2ex}
		\begin{subfigure}[b]{0.48\linewidth}
		\includegraphics[width=\linewidth]{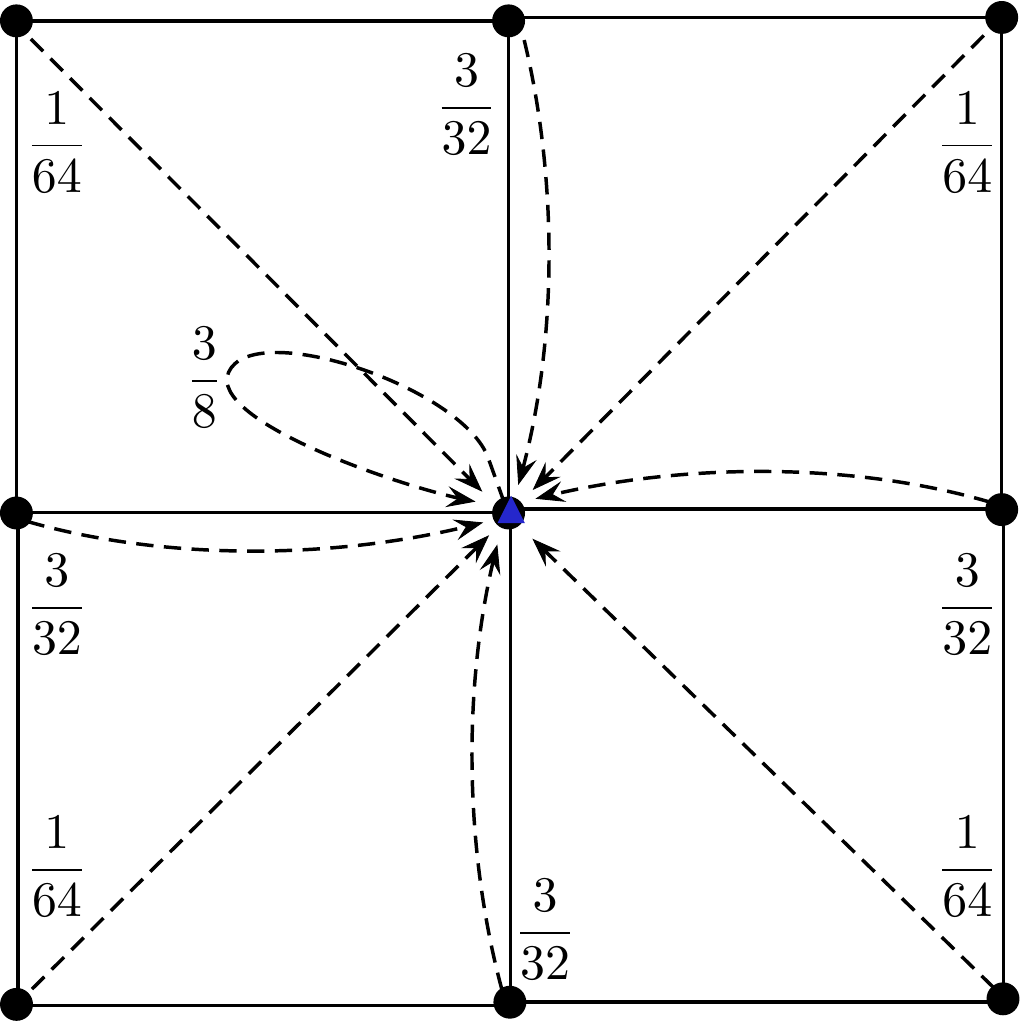}
		\caption{Vertex point.}
		\label{fig:CC_vertex_point}
	\end{subfigure}
	\begin{subfigure}[b]{0.25\linewidth}
		\includegraphics[width=\linewidth]{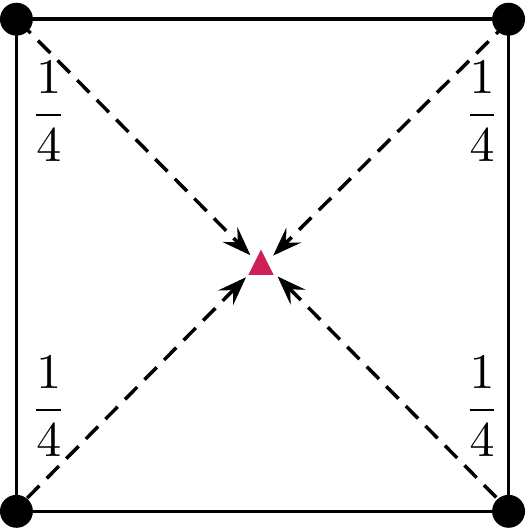}
		\caption{Face point.}
		\label{fig:CC_face_point}
	\end{subfigure}
	\caption{The weight distribution for computing different types of new control points.}
	\label{fig:CC_surface_weights}
\end{figure}

\begin{figure}
\centering
\includegraphics[width=0.65\linewidth]{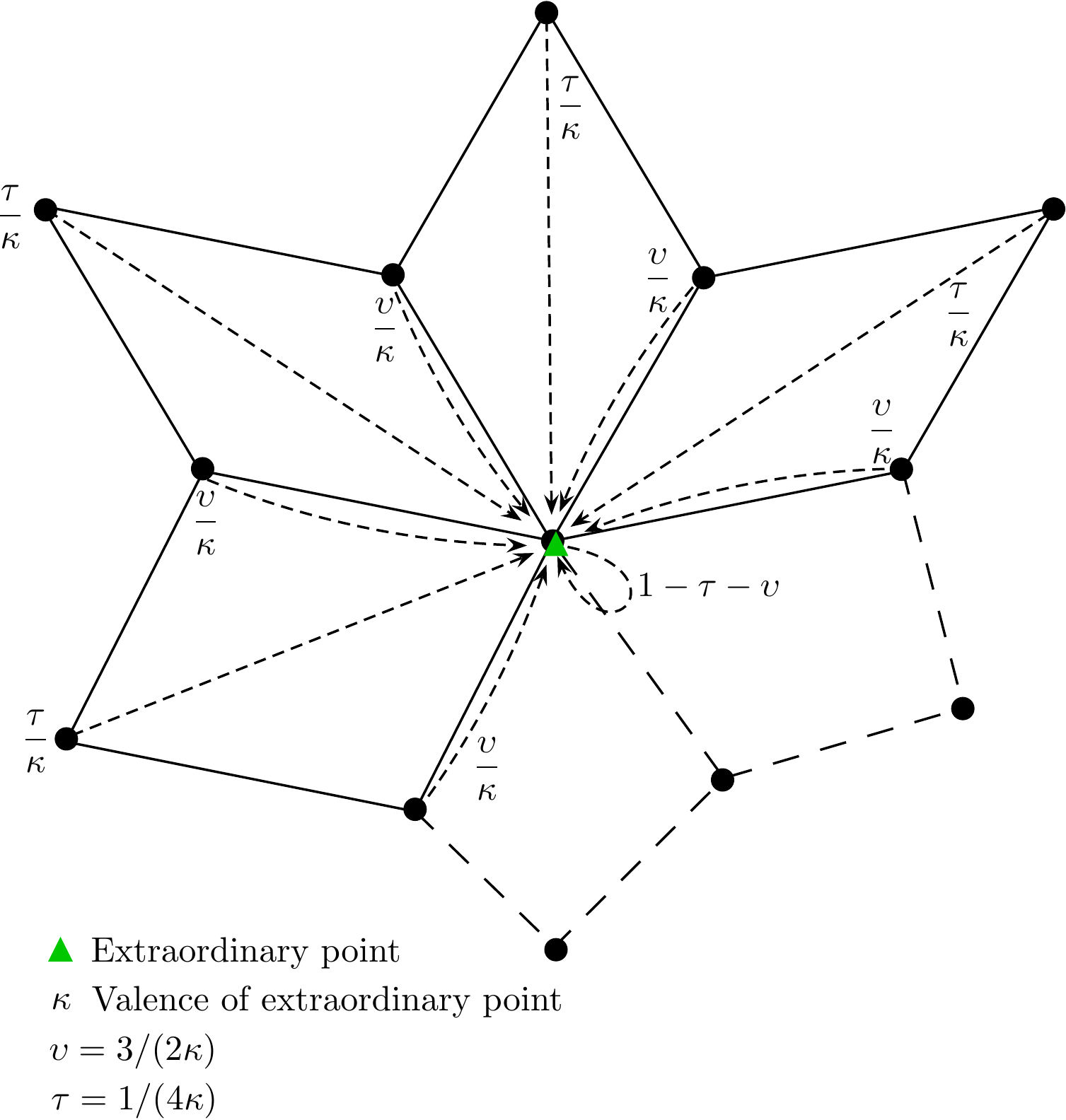}
\caption{Weight distributions for computing an extraordinary point with valence $\kappa$.}
\label{fig:CC_extraordinary_point}
\end{figure}

Similar to the one-dimensional curve, the new refined control points can be classified into three types: face points, edge points and vertex points. The face points in the $i^\text{th}$ refinement are computed as
\begin{equation}
\mathbf{P}^i_{2j,2k} = \frac{1}{4}\left[\mathbf{P}^{i-1}_{j,k}+\mathbf{P}^{i-1}_{j,k+1} + \mathbf{P}^{i-1}_{j+1,k} + \mathbf{P}^{i-1}_{j+1,k+1}\right],
\end{equation}
where $j$ and $k$ are indices of control points for orthogonal directions. The face point is the central point of the original face. The edge point is computed as
\begin{equation}
\mathbf{P}^i_{2j+1,2k} = \frac{1}{16}\left[\mathbf{P}^{i-1}_{j,k} + 6\mathbf{P}^{i-1}_{j,k+1} +\mathbf{P}^{i-1}_{j,k+2} +\mathbf{P}^{i-1}_{j+1,k} + 6\mathbf{P}^{i-1}_{j+1,k+1} + \mathbf{P}^{i-1}_{j+1,k+2}\right],
\end{equation}
or likewise
\begin{equation}
\mathbf{P}^i_{2j,2k+1} = \frac{1}{16}\left[\mathbf{P}^{i-1}_{j,k} + \mathbf{P}^{i-1}_{j,k+1} + 6\mathbf{P}^{i-1}_{j+1,k} +6\mathbf{P}^{i-1}_{j+1,k+1} + \mathbf{P}^{i-1}_{j+2,k} + \mathbf{P}^{i-1}_{j+2,k+1}\right].
\end{equation}
The `vertex point' is computed as
\begin{equation}
\begin{split}
\mathbf{P}^i_{2j+1,2k+1} = \frac{1}{64}  \Big[ &\mathbf{P}^{i-1}_{j,k} + 6\mathbf{P}^{i-1}_{j,k+1} +\mathbf{P}^{i-1}_{j,k+2} + 6\mathbf{P}^{i-1}_{j+1,k} +36\mathbf{P}^{i-1}_{j+1,k+1} + 6\mathbf{P}^{i-1}_{j+1,k+2}\\
&+ \mathbf{P}^{i-1}_{j+2,k} + 6\mathbf{P}^{i-1}_{j+2,k+1} + \mathbf{P}^{i-1}_{j+2,k+2}\Big].\\
\end{split}
\end{equation}
Equipped with these formulae, the new control points on the $i^\text{th}$ level of refinement $\mathscr {P}^i$ can be computed as:
\begin{equation}
\mathscr {P}^i = \ary {S} \mathscr {P}^{i-1},
\label{eq:subd_operator}
\end{equation}
$ \ary {S}$ is a subdivision operator -- a matrix consisting of a set of weights. Each weight is associated with a control point in $\mathscr P^{i-1}$. The weight distributions for different types of control points are shown in Figure~\ref{fig:CC_surface_weights}. The weight distributions for extraordinary point are shown in Figure~\ref{fig:CC_extraordinary_point}. After successive levels of refinements, a smooth B-spline surfaces is obtained. 

%
\subsection{Computing control point set for sub-elements}
\label{ap:bases_ev}
We denote the control points of an irregular patch in Figure~\ref{fig:Catmull_clark_ex_mesh} as a set $\mathscr{P}$. The initial control points of the patch are expressed as
\begin{equation}
\mathscr{P}_0 = \left\{\mathbf{P}^0_{0}, \mathbf{P}^0_{1}, \dots, \mathbf{P}^0_{{2\kappa+6}},  \mathbf{P}^0_{{2\kappa+7}}\right\}.
\end{equation}
Through one level of subdivision we generate $2\kappa + 17$ new control points ($\kappa$ is the valence), as shown in Figure~\ref{fig:catmull_clark_ex}, denoted by
\begin{equation}
\mathscr{P}_1 = \left\{\mathbf{P}^1_{0}, \mathbf{P}^1_{1}, \dots, \mathbf{P}^1_{{2\kappa+15}},\mathbf{P}^1_{{2\kappa+16}}\right\}.
\end{equation}
The subdivision step is represented as
\begin{equation}
\mathscr{P}_1 = \ary{A}\mathscr{P}_0,
\end{equation}
where $\ary A$ is the subdivision operator given by
\begin{equation}
\ary{A} = \left[
\begin{array}{cc}
\ary{S} & \ary{0}\\
\ary{S}_{11} & \ary{S}_{12}\\
\ary{S}_{21} & \ary{S}_{22}\\
\end{array}
\right].
\end{equation}
The terms $\ary{S}$, $\ary{S}_{11}$, $\ary{S}_{12}$, $\ary{S}_{21}$ and $\ary{S}_{22}$ are defined in~\cite{stam1998exact} and $\ary{S}$ is given in Equation~\ref{eq:subd_operator}. To evaluate the sub-element $\Omega_1$, $\Omega_2$ and $\Omega_3$ in Figure~\ref{fig:catmull_clark_ex}, one needs to pick $2\kappa+8$ control points out of the new $2\kappa+17$ control point patch. A selection operator $\ary D+\kappa$ for sub-element $\Omega_k$ and $ k = 1, 2, 3$ is used to select the necessary control points from $\mathscr{P}_1$, that is
\begin{equation}
\mathscr{P}_{1,k} = \ary{D}_k\mathscr{P}_1.
\end{equation}
Then a surface point can be evaluated with the cubic spline basis functions as
\begin{equation}
\mathbf{x}(\boldsymbol{\xi}) =  \sum_{A=0}^{15} N_A(\boldsymbol{\xi})\, \mathbf{P}^{1,k}_{A}.
\end{equation}
As shown in Figure~\ref{fig:CC_domain_subdivision}, after successive subdivisions, the non-evaluable element can be limited to a negligible region. 

Assume the target point has parametric coordinates $\boldsymbol{\xi} = (\xi, \eta)$. One first determines how many subdivisions are required for this point by:
\begin{equation}
n = \left \lfloor{\min \left(-\log_2(\xi), -\log_2(\eta)\right)+1}\right \rfloor.
\label{eq:n}
\end{equation}
The sub-element index $k$ is determined as
\begin{equation}
k = \left\{
\begin{split}
 &1 \quad  \text{if} \quad \boldsymbol{\xi} \in \left[\frac{1}{2^n}, \frac{1}{2^{n-1}}\right]\times \left[0, \frac{1}{2^n}\right],\\
 &2 \quad  \text{if} \quad \boldsymbol{\xi} \in \left[\frac{1}{2^n}, \frac{1}{2^{n-1}}\right]\times \left[\frac{1}{2^n}, \frac{1}{2^{n-1}}\right],\\
 &3 \quad  \text{if} \quad \boldsymbol{\xi} \in \left[0, \frac{1}{2^{n}}\right]\times \left[\frac{1}{2^n}, \frac{1}{2^{n-1}}\right].\\
\end{split}
\right.
\end{equation}
The surface point $\mathbf{x}$ is located in the regular sub-element $k$ after the $n^\text{th}$ refinement. The patch for this element is picked with the selection operator $\ary{D}_k$ as
\begin{equation}
\mathscr P_{n,k} = \ary D_k \mathscr P_n.
\end{equation}
The enlarged set $\mathscr P_n$ contains $2\kappa + 17$ control vertices, which is generated from the subdivision of ${\mathscr P^\ast}_{n-1}$ as
\begin{equation}
\mathscr P_n = \ary A {\mathscr P^\ast}_{n-1}.
\end{equation}
The set ${\mathscr P^\ast}_{n-1}$ has $2\kappa + 8$ control vertices. It is successively refined from the initial set $\mathscr P_0$ as
\begin{equation}
{\mathscr P^\ast}_{n-1} = \bar{\ary A}^{n-1} \mathscr P_0,
\end{equation}
where $\bar{\ary A}$ is a square matrix operator which subdivides the patch for computing the new patch for the irregular element, and is defined by
\begin{equation}
\bar{\ary{A}} = \left[
\begin{array}{cc}
\ary{S} & \ary{0}\\
\ary{S}_{11} & \ary{S}_{12}\\
\end{array}
\right].
\end{equation}

\end{document}